\documentclass[12pt,maitrise]{dms}

\usepackage[latin1]{inputenc}
\usepackage[T1]{fontenc}

\usepackage{amsfonts}
\usepackage{amssymb}
\usepackage{amsmath}
\usepackage{mathrsfs}
\usepackage{verbatim}

\newif\ifPDF
\ifx\pdfoutput\undefined\PDFfalse

\else\ifnum\pdfoutput > 0\PDFtrue
  \else\PDFfalse
    \fi
\fi

\ifPDF 
        \usepackage[pdftex]{color,graphicx}
        \usepackage[pdfauthor={Claude Gravel},%
                    pdftitle={Géométrie spectrale sur le disque : loi de Weyl et ensembles nodaux},%
                    pdfduplex=Simplex,%
                    pdfsubject={Mémoire de maîtrise, Université de Montréal},%
                    pdfkeywords={Weyl, nodal domains, spectral geometry, géométrie spectrale, ensembles nodaux, zeros of Bessel functions, zéros de fonctions de Bessel},
                    pdfdisplaydoctitle=true,%
                    pdfstartpage=1,%
                    pdftex]{hyperref}
        \hypersetup{colorlinks,
                    citecolor=green,%
                    filecolor=red,%
                    linkcolor=blue,%
                    urlcolor=red}

\else 
    \usepackage[dvips]{color,graphicx,hyperref}
    \usepackage{epsfig}
    \hypersetup{colorlinks,
                citecolor=red,%
                filecolor=red,%
                linkcolor=red,%
                urlcolor=red,%
                dvips}
\fi

\newtheorem{Def}{Définition}[section]

\newtheorem{Lem}{Lemme}[section]
\newtheorem{Prop}{Proposition}[section]
\newtheorem{Rem}{Remarque}[section]
\newtheorem{Thm}{Théorème}[section]

\theoremstyle{definition}\newtheorem{Prflemme}{Preuve du
Lemme}[section]
\newtheorem{Dem}{Démonstration}[section]

\numberwithin{equation}{section}
\numberwithin{table}{chapter}
\numberwithin{figure}{chapter}

\usepackage{attachfile}

\begin{document}

\title{Géométrie spectrale sur le disque: loi de Weyl et ensembles nodaux}

\author{Claude Gravel}

\copyrightyear{2007}


\directeur{Iosif Polterovich}

\membrejury{Marlène Frigon (Université de Montréal), Nilima Nigam
(McGill University)}


\sujet{mathématiques}
\orientation{mathématiques pures}

\pagenumbering{roman} \maketitle

\chapter*{Sommaire}

\noindent Soit $D$ le disque de rayon 1 ou soit $S(\alpha)$ le
secteur du disque de rayon 1 et d'angle $\alpha\in(0,2\pi)$. Soit le
Laplacien
$\triangle=\frac{\partial^{2}}{\partial{x}^{2}}+\frac{\partial^{2}}{\partial{y}^{2}}$
sur $A$ où $A=D$ ou $A=S(\alpha)$.

\noindent Dans ce mémoire, des questions relatives au problème des
valeurs propres de l'opérateur $\triangle$ sur $A$ c'est-à-dire
relative aux solutions de l'équation $\triangle(u)+\lambda{u}=0$
avec des conditions de Dirichlet aux frontières seront abordées.
D'abord nous rappelons que l'ensemble des valeurs admises pour
$\lambda$ c'est-à-dire l'ensemble des solutions pour $\lambda$ sont
les zéros des fonctions de Bessel de premier type. Cet ensemble de
valeurs propres s'appelle le spectre, c'est l'ensemble des valeurs
de $\lambda$ étant solutions de l'équation précédente. Si le spectre
est dénoté par $\mathrm{Sp}_{A}(-\triangle)$ et si $j_{\nu,k}$
dénote le k\ieme{} zéro de la fonction de Bessel d'ordre $\nu\geq0$
de premier type, alors
\begin{eqnarray*}
\mathrm{Sp}_{S(\alpha)}(-\triangle)&=&\big\{j^{2}_{\frac{n\pi}{\alpha},k}\phantom{1};\phantom{1}n=1,2,\ldots\phantom{1}\textrm{et}\phantom{1}k=1,2,\ldots\big\}\\
\mathrm{Sp}_{D}(-\triangle)&=&\big\{j^{2}_{n,k}\phantom{1};\phantom{1}n=0,1,2,\ldots\phantom{1}\textrm{et}\phantom{1}k=1,2,\ldots\big\}.
\end{eqnarray*}

\noindent Le premier problème abordé consistera à déterminer la
première ligne où la deuxième fonction propre s'annule dans le cas
d'un secteur $S(\alpha)$. L'ensemble nodal d'une fonction propre
$u_{j}$ est le sous-ensemble $\mathrm{Z}_{j}\subset A$ sur lequel la
fonction propre s'annule, c'est l'ensemble des zéros de $u_{j}$. Les
domaines nodaux de $u_{j}$ sont les composantes connexes du
complément de l'ensemble nodal, c'est l'ensemble des points du
domaines où $u_{j}$ ne s'annule pas. Il est utile de se rappeler que
la question précédente ne se pose pas pour la première fonction
propre puisque cette dernière n'a que la frontière comme seule ligne
nodale. Comme il y sera montré, si l'angle du secteur est
suffisamment petit, alors la ligne nodale est décrite par l'équation
en coordonnée polaire $r=\frac{j^{2}_{\nu,1}}{j^{2}_{\nu,2}}<1$ où
$\nu=\frac{2\pi}{\alpha}$ sinon elle est décrite par
$\theta=\frac{\alpha}{2}$. Nous verrons également qu'il existe un
angle critique $\alpha_{0}$ pour lequel la ligne nodale est
indéterminée et que cet angle critique est donné implicitement par
l'équation $j_{1}(2\nu)-j_{2}(\nu)=0$ où $\nu=\frac{\pi}{\alpha}$.
De même, lorsque l'angle $\alpha=\alpha_{0}$, nous verrons que la
multiplicité de 2\ieme{} valeur propre est double.

\noindent Le deuxième problème consiste à construire un algorithme
pour obtenir les $m$ premières valeurs propres, cet algorithme est
également efficace pour des valeurs larges de $m$. Cet algorithme ne
doit pas résoudre les systèmes d'équations différentielles
ordinaires relatives au problème $\triangle(u)+\lambda{u}=0$ puisque
la résolution de tels systèmes est très coûteuse en temps et en
mémoire voire presque irréalisable pour des valeurs propres très
larges. Le problème d'ordonner les valeurs propres pour lesquelles
une caractérisation de deux paramètres existe est équivalent à
connaître la structure des domaines nodaux des fonctions propres
$\{u_{j}\}_{j=1}^{m}$. La clef de voûte derrière ce problème est
d'utliser un théorème dû à Courant sur le nombre de domaines nodaux
d'une fonction propre.

\noindent Troisièmement, soit $N(\lambda)$ la fonction de compte des
valeurs propres c'est-à-dire
\begin{displaymath}
N(\lambda)=\sum_{\lambda_{j}\leq\lambda}{1}.
\end{displaymath}
Une explication approfondie d'un théorème dû à Kuznetsov et Fedosov
montrant que
\begin{displaymath}
N(\lambda)=\frac{\lambda}{4}-\frac{\sqrt{\lambda}}{2}+O\big(\lambda^{\frac{1}{3}}\big)
\end{displaymath}
sera donnée. Pour y arriver, un \og grand \fg{} détour sur
l'approximation uniforme des fonctions de Bessel de premier type
d'ordre réel strictement positif et d'argument positif.

\noindent Finalement, un algorithme permettant d'évaluer exactement
et efficacement $N(\lambda)$ sans avoir à trouver toutes les valeurs
propres $\lambda_{j}\leq\lambda$ sera donné. Cet algorithme
n'utilise que la monotonicité des zéros des fonctions de Bessel et
peut trouver par exemple $N(9\cdot{10^{8}})$ rapidement. Ce 2\ieme{}
algorithme est nettement plus efficace à utiliser pour évaluer
$N(\lambda)$ que si nous trouvions naïvement toutes les valeurs
propres $\lambda_{j}\leq \lambda$ en utilisant le 1\ier{}
algorithme.

\chapter*{Summary}

\noindent Let $D$ be the unit disc, $S(\alpha)$ be a sector of the
unit disc with angle $\alpha \in (0,2\pi)$, and let
$\triangle=\frac{\partial^{2}}{\partial{x}^{2}}+\frac{\partial^{2}}{\partial{y}^{2}}$
be the Laplacian.

\noindent In this master's thesis, we will consider several
questions arising from the study of the eigenvalue problem
$\triangle(u) + \lambda u = 0$ on $A$ with Dirichlet boundary
conditions, where $A = D$ or $S(\alpha)$. We proceed by first
recalling that the spectrum of the negative Laplacian (i.e. the
admissible set of eigenvalues $\lambda$ for the above eigenvalue
problem) is given by the zeros of the Bessel function of the first
kind. In other words, if we denote by $\mathrm{Sp}_{A}(-\triangle)$
the spectrum of the operator $-\triangle$ on $A$, and by $j_{\nu,k}$
the k$^{\mathrm{th}}$ zero of the Bessel function of the fisrt kind
of order $\nu\geq 0$, then
\begin{eqnarray*}
\mathrm{Sp}_{S(\alpha)}(-\triangle)&=&\big\{j^{2}_{\frac{n\pi}{\alpha},k}\phantom{1};\phantom{1}n=1,2,\ldots\phantom{1}\textrm{et}\phantom{1}k=1,2,\ldots\big\}\\
\mathrm{Sp}_{D}(-\triangle)&=&\big\{j^{2}_{n,k}\phantom{1};\phantom{1}n=0,1,2,\ldots\phantom{1}\textrm{et}\phantom{1}k=1,2,\ldots\big\}.
\end{eqnarray*}

\noindent The first problem consists of the determination of the
first nodal line of the second eigenfuntion of a sector $S(\alpha)$.
The nodal set of an eigenfunction $u_{j}$ is the subset
$\mathrm{Z}_{j} \subset A$ over which the eigenfunction vanishes.
Nodal domains of an eigenfunction $u_{j}$ are the connected
components of the complement of the nodal set i.e. of the set at
which the eigenfunction does not vanish. Note that this
corresponding problem for the first eigenfunction is uninteresting
because the nodal set would consist of the boundary $\partial{A}$
itself. We will see that when $\alpha$ is small enough, the nodal
curve is described by the equation
$r=\frac{j^{2}_{\nu,1}}{j^{2}_{\nu,2}}<1$, where
$\nu=\frac{2\pi}{\alpha}$, whereas for large $\alpha$ it is
described by $\theta=\frac{\alpha}{2}$. Thus, there exists a
critical angle $\alpha_0$ for which the second nodal curve is
undefined, and this angle is given implicitly by solving
$j_{2}(\nu)-j_{1}(2\nu)=0$ where $\nu=\frac{\pi}{\alpha}$. Also the
second eigenvalue has multiplicity 2 when the angle of the sector is
$\alpha_{0}$.

\noindent The second problem is to devise an algorithm for obtaining
the first $m$ eigenvalues that is also efficient for large $m$. Such
an algorithm must not rely on the explicit integration of the
underlying ODEs, since the solution of such ODEs would become
prohibitively expensive when large eigenvalues are involved. We will
see that the problem of properly ordering the eigenvalues (for which
a two-parameter characterization is known) is equivalent to knowing
the structure of the nodal sets of the first $m$ eigenfunctions
$\{u_{j}\}_{j=1}^{m}$. Thus, the key insight behind the algorithm is
to perform the ordering using nodal set information, which in turn
is provided by Courant's theorem on the number of nodal domains of a
given eigenfunction.

\noindent Thirdly, let $N(\lambda)$ be the counting function of the
eigenvalues i.e
\begin{displaymath}
N(\lambda)=\sum_{\lambda_{j}\leq\lambda}{1}.
\end{displaymath}
A detailed explanation of a theorem given by Kuznetsov and Fedosov
showing that
\begin{displaymath}
N(\lambda)=\frac{\lambda}{4}-\frac{\sqrt{\lambda}}{2}+O\big(\lambda^{\frac{1}{3}}\big)
\end{displaymath}
will be given. In order to get there, we will look at some results
concerning the uniform expansion Bessel functions of the first kind
of strictly positive real order and evaluated at $x \in R^+$.

\noindent Finally, we will present an algorithm that allows the
exact and efficient evaluation of $N(\lambda)$, without having to
compute all of the smaller eigenvalues $\lambda_{j} \leq \lambda$.
The proposed algorithm only relies on the monotonicity of the zeros
of the Bessel functions and can, for instance, find quickly
$N(9\cdot{10^{8}})$. The latter algorithm provides a much more
efficient alternative to the naive approach of finding, with the
help of the 1$^{\mathrm{st}}$ algorithm, all the eigenvalues
$\lambda_{j}$ such that $\lambda_{j}\leq \lambda$.

\tableofcontents \listoffigures \listoftables

\chapter*{Remerciements}

\noindent Je remercie mon directeur, Iosif Polterovich, pour ses
conseils et son apport financier durant la rédaction de ce mémoire.
Je remercie également les membres du jury, Marlène Frigon et Nilima
Nigam, pour leurs remarques qui ont amélioré mon mémoire. Je
remercie ma famille et mes amis pour leurs encouragements et le
temps qu'ils prennent avec moi pour me changer les idées lorsque
cela est nécessaire. Un merci spécial pour les discussions
éclaircissantes à Dominique Rabet et Igor Wigman pour la section 3.2
de ce mémoire. Je remercie également Felix Kwok pour la correction
de mon résumé en anglais.

\NoChapterPageNumber \pagenumbering{arabic}

\chapter*{Introduction}

\noindent Ce mémoire est un exemple de ce qu'il est possible
d'étudier en géométrie spectrale. Puisqu'il est généralement
difficile de définir et de délimiter exactement un domaine de la
connaissance, le lecteur trouvera dans ce mémoire quelques exemples
de ce qu'il est possible d'étudier en géométrie spectrale.
Néanmoins, personne ne se trompe en affirmant que la géométrie
spectrale est un domaine alliant la géométrie et l'analyse
spectrale. Traditionnellement, l'analyse spectrale étudie le spectre
des opérateurs comme les matrices ou encore des opérateurs
différentiels comme l'opérateur $\frac{d^{2}}{dx^{2}}$. Généralement
les opérateurs agissent sur des objets comme des vecteurs ou des
fonctions. Les domaines sur lesquels ces objets sont définis sont
plats généralement comme un carré en dimension 2 ou un cube en
dimension 3. En changeant les domaines sur lesquels ces objects sont
définis, par exemple en considérant des vecteurs sur la sphère en
identifiant chaque point antipodal comme équivalent, les opérateurs
agissant sur ces objets ont des \og comportements \fg{} différents.
Le spectre des opérateurs est entre autres différent d'une géométrie
à une autre. La géométrie spectrale étudie le spectre des opérateurs
agissant sur des objets définis sur des géométries courbes comme la
sphère ou le tore. Comme il n'est pas nécessaire de faire intervenir
l'ensemble des vérités admises concernant un énoncé, la phrase
précédente est certainement vraie mais non complète.

\noindent Plusieurs problèmes peuvent être étudiés en théorie
spectrale. Ces problèmes varient selon que les opérateurs possèdent
certaines proprités comme entre autres être compact, linéaire,
auto-adjoint, auto-adjoint non-borné, etc. Pour certains types
d'opérateurs, comme le Laplacien $\triangle$ qui sera d'intérêt dans
ce qui suit, il est possible de décomposer leur spectre au sens
exposé au chapitre 1. De par la forme du domaine, il peut être
possible de connaître explicitement la forme des fonctions propres.

\noindent De façon générale, lorsque les mathématiciens ont une
fonction sous leurs yeux, ils aiment connaître où elle s'annule.
C'est la même la chose avec les fonctions propres. L'ensemble des
zéros d'une fonction propre donnée est un ensemble fermé de la
variété sur laquelle lesdites fonctions propres sont définies. Cet
ensemble porte un nom particulier: il s'appelle le domaine nodal de
la fonction propre donnée. Le complément de l'ensemble des zéros est
ouvert et \og découpe \fg{} la variété en composantes connexes.

\noindent Également, lorsque l'ensemble des fonctions propres
forment une base séparable et complète pour l'espace
$\mathrm{L}^{2}$ sur le domaine en question, il est intéressant de
s'intéresser aux problèmes asymptotiques. Un de ces problèmes est
relié à la distribution empirique des valeurs propres. Les valeurs
propres des opérateurs définis positifs ont des valeurs propres
positives. Par conséquent, étant donné une valeur fixe positive,
combien de valeurs propres sont plus petites que cette valeur fixée
à priori en tenant compte de la multiplicité?

\noindent Dans ce présent manuscrit, l'analyse spectrale se fera sur
le disque $D$ de rayon 1 ou sur un secteur $S(\alpha)$ d'angle
$\alpha\in(0,2\pi)$ de rayon 1. $\partial D$ ou $\partial S(\alpha)$
dénotent la frontière de $D$ ou $S(\alpha)$. L'opérateur d'intérêt
sera évidemment le Laplacien dénoté $\triangle$ ou $\nabla^{2}$.
Connaître le spectre du Laplacien permet de connaître dans ce cas,
et dans bien d'autres cas, le spectre de l'opérateur des ondes, de
la chaleur et Schrödinger. L'espace sur lequel $\triangle$ est
défini se dénote $H_{0}^{1}$ et il s'appelle l'espace de Sobolev.
Les fonctions propres formant une base pour $H_{0}^{1}(D)$ avec des
conditions de Dirichlet aux frontières sont les fonctions de Bessel
de premier type d'ordre $0,1,2,\ldots$. Les fonctions propres
formant une base pour $H_{0}^{1}(S(\alpha))$ sont les fonctions de
Bessel de premier type $1,2,\ldots$.

\noindent Deux problèmes concernant l'opérateur $\triangle$ sur $D$
seront étudiés plus particulièrement. Le 1\ier{} problème consistera
à étudier la configuration des ensembles nodaux des $m$ premières
fonctions propres $u_{j}$ pour $j=1,\ldots,m$. Ce problème sera
équivalent à calculer les $m$ premières valeurs propres
$\lambda_{j}$ pour $j=1,\ldots,m$. Bref, c'est en voulant obtenir un
algorithme efficace pour obtenir les $m$ premières valeurs propres
qu'il sera évident comment obtenir la configuration des ensembles
nodaux. Un théorème dû à Courant donnant une borne sur le nombre de
composantes connexes de la m\ieme{} valeur propre sera l'outil
principal afin d'établir les configurations. Les ensembles nodaux
seront des lignes angulaires partant du centre du cercle et des
series de cercles concentriques. De même, un autre problème relié à
l'identification de la prèmière ligne nodale de la deuxième fonction
propre sur un secteur $S(\alpha)$ sera étudié. On mentionnera la
valeur critique de l'angle $\alpha$ dans le cas du secteur
$S(\alpha)$ où la ligne décrite par l'équation en coordonnée polaire
$r=\mathrm{const.}$ est préférée à la ligne décrite par l'équation
$\theta=\mathrm{const.}$. Nous verrons que si l'angle du secteur est
critique, alors la 1\iere{} ligne nodale n'est pas définie. De même
nous verrons que la multiplicité de la 2\ieme{} valeur propre
$\lambda_{2}$ est double pour ce secteur d'angle critique.

\noindent Un autre problème étudié est relié à la distribution
asymptotique des valeurs propres de $\triangle$. En d'autres termes,
soit $\textrm{Sp}_{D}(-\triangle)=\{\lambda_{j}\}_{j>0}$ le spectre
de l'opérateur $\triangle$ sur le disque $D$ de rayon 1 avec les
valeurs propres ordonnées naturellement c'est-à-dire
$0<\lambda_{1}\leq\lambda_{2}\leq\ldots$ Soit $N(\lambda)$, la
fonction de compte des valeurs propres $\lambda_{j}$ plus petites
que $\lambda$ c'est-à-dire soit
\begin{displaymath}
N(\lambda)=\sum_{\lambda_{j}\leq\lambda}{1}
\end{displaymath} Pour des valeurs de $\lambda$
très grandes,
\begin{displaymath}
N(\lambda)=\frac{\lambda}{4}-\frac{\sqrt{\lambda}}{2}+R(\lambda)\phantom{12}\textrm{(Loi
de Weyl)}
\end{displaymath}
où $R(\lambda)$ est le terme d'erreur. Il sera d'abord réexpliqué
clairement comment il a été démontré dans \cite{KuFe_1965} que
$R(\lambda)=O(\lambda^{\frac{1}{3}})$. Enfin, un deuxième algorithme
permettant d'évaluer $N(\lambda)$ exactement et efficacement est
donné. Cet algorithme n'utilise que la monotonicité des zéros des
fonctions de Bessel. Comme il y sera expliqué, l'algorithme fait une
\og marche \fg{} qui consiste à faire des retours en arrière et des
montées en alternance sur des paires d'entiers bien déterminées par
le problème en soit. Les coordonnées de ces paires sont l'ordre et
l'index des zéros des fonctions de Bessel de premier type pouvant
être les candidats possibles pour les valeurs propres.

\chapter[]{Préliminaires}

\noindent Voici le plan pour ce chapitre. Les références sont
données au début de chaque section.

\noindent Dans la \textbf{1\iere{}} section, certaines propriétés
d'intérêt pour ce mémoire de l'opérateur Laplacien dénoté
$\triangle$ sur un ouvert régulier borné $\Omega$ de
$\mathbb{R}^{n}$ sont mentionnées sans preuve. Plusieurs définitions
sont données avant d'y arriver permettant de bien définir le domaine
de $\triangle$ et d'arriver à quelques théorèmes intéressants. La
loi de Weyl (fonction de compte des valeurs propres) sera également
présentée avec quelques théorèmes pour l'approximer asymptotiquement
sur des domaines $\Omega$ arbitraires.

\noindent Dans la \textbf{2\ieme{}} section, il sera montré que
résoudre grâce à la séparation de variables le problème à valeurs
propres $\triangle{u}+\lambda{u}=0$ sur un domaine circulaire
c'est-à-dire lorsque $\Omega$ est un secteur ou un disque \og
engendre \fg{} l'équation de Bessel. Il sera montré quelles valeurs
$\lambda$ peut prendre engendrant ainsi la suite spectrale ou le
spectre. Il sera montré que les valeurs de $\lambda$ sont des carrés
des zéros de fonctions de Bessel.

\noindent Dans la \textbf{3\ieme{}} section, les zéros des fonctions
de Bessel seront appronfondis. Plus spécifiquement, le développement
asymptotique dû à Olver des zéros sera présenté sans preuve. Ce
développement est important pour quiconque veut calculer les valeurs
propres de $\triangle$ sans avoir à résoudre les équations
différentielles directement. Au chapitre 2, il sera montré comment
en effet il est possible d'obtenir le spectre, les $m$ premières
valeurs propres pour $m$ très grand, du Laplacien. Cette section
introduira succinctement les fonctions d'Airy comme des solutions
d'une équation différentielle du même nom. Les zéros de ces
dernières et leurs relations aux fonctions de Bessel seront abordés
brièvement pour nous permettre de comprendre le développement
asymptotique d'Olver.

\noindent Dans la \textbf{4\ieme{}} section, les ensembles des zéros
des fonctions propres appelés communément les ensembles nodaux des
fonctions propres seront étudiés. Un théorème dû à Courant sur le
nombre de domaines nodaux sera énoncé sans preuve. Cette section
avec la précédente seront utiles pour le chapitre 2.

\section{Rappels - définition du Laplacien, décomposition du spectre et fonction de compte}

\noindent Le lecteur retrouvera les définitions suivantes et
certains théorèmes sur les fonctions et les valeurs propres du
Laplacien dans l'ouvrage de références \cite{Ra_1999} et
\cite{Wi_2003}.

\noindent Soit $\Omega$, un ouvert borné régulier de
$\mathbb{R}^{n}$. Soit $\mathrm{C}^{\infty}(\Omega)$ l'espace des
fonctions définies sur $\Omega$ différentiables indéfiniment.

\begin{Def}[Le Laplacien $\triangle$]
\noindent  Soit $u\in\mathrm{C}^{\infty}(\Omega)$. Alors le
Laplacien est l'opérateur
$\triangle:\mathrm{C}^{\infty}(\Omega)\rightarrow\mathrm{C}^{\infty}(\Omega)$
tel que
\begin{displaymath}
\triangle{u}=\frac{\partial^{2}u}{\partial{x_{2}}^{2}}+\ldots+\frac{\partial^{2}u}{\partial{x_{n}}^{2}}.
\end{displaymath}
\end{Def}

\begin{Def}[Fonction propre et valeur propre]
Une \textbf{fonction propre} relative à la \textbf{valeur propre}
$\lambda$ est une solution non nulle du problème
\begin{displaymath}
\left\{ \begin{array}{ll}
    -\triangle{u} = \lambda{u} & \textrm{dans}\phantom{1}\Omega\\
    \phantom{-T}{u} = 0 & \textrm{sur}\phantom{1}\partial{\Omega}
\end{array}\right.
\end{displaymath}
\end{Def}

\noindent Le théorème suivant (page 115 de \cite{Wi_2003}) montre
que les fonctions propres du Laplacien constituent une base de
$\mathrm{L}^{2}(\Omega)$.
\begin{Thm}[Décomposition du spectre]
Il existe une suite non bornée de valeurs propres de
\begin{displaymath}
0<\lambda_{1}\leq\lambda_{2}\leq\ldots
\end{displaymath}
et une suite de fonctions propres relatives
$\{\varphi_{j}\}_{j}^{\infty}$ qui est une base hilbertienne de
$\mathrm{L}^{2}(\Omega)$.
\end{Thm}

\noindent Les valeurs propres peuvent être exprimées comme dans la
proposition suivante due Rayleigh (page 115 de \cite{Wi_2003}).
\begin{Prop}[Quotient de Rayleigh]
Pour tout $n\geq{1}$,
\begin{displaymath}
\lambda_{n}=\min\Big\{\int_{\Omega}{|\nabla
u|^{2}dx}\phantom{1}\big|\phantom{1}\int_{\Omega}{u^{2}dx}=1,\phantom{1}\int_{\Omega}{u\varphi_{j}dx}=0\phantom{1}\textrm{pour}\phantom{1}j=1,\ldots,n-1\Big\}
\end{displaymath}
\end{Prop}

\noindent Pour les résultats qui suivent concernant la fonction de
compte des valeurs propres $N_{\Omega}(\lambda)$, ils se retrouvent
dans \cite{Be_2003}.

\noindent Soit maintenant $N_{\Omega}(\lambda)$ la fonction de
compte des valeurs propres c'est-à-dire
\begin{displaymath}
N_{\Omega}(\lambda)=\sum_{\lambda_{j}\leq\lambda}{1}.
\end{displaymath}
La fonction de compte est une fonction à valeurs entières définie
sur les réels c'est-à-dire
\begin{displaymath}
N_{\Omega}:\mathbb{R}\to\mathbb{N}
\end{displaymath}
qui est également continue à droite c'est-à-dire pour $h>0$
\begin{displaymath}
N_{\Omega}(\lambda)=\lim_{h\to 0^{+}}{N_{\Omega}(\lambda+h)}
\end{displaymath}
Si $\lambda\in\mathbb{R}$ est une valeur propre, alors sa
multiplicité est donnée par
\begin{displaymath}
\lim_{h\to
0^{+}}{\big(N_{\Omega}(\lambda+h)-N_{\Omega}(\lambda-h)\big)}.
\end{displaymath}
Si $\lambda\in\mathbb{R}$ n'est pas une valeur propre, alors
\begin{displaymath}
\lim_{h\to
0^{+}}{\big(N_{\Omega}(\lambda+h)-N_{\Omega}(\lambda-h)\big)}=0.
\end{displaymath}
Par conséquent, si $\lambda$ est une valeur propre, alors il existe
$j\in\mathbb{N}$ tel que $\lambda=\lambda_{j}$ et tel que
\begin{displaymath}
\lim_{h\to 0^{+}}{N_{\Omega}(\lambda-h)}=\lim_{h\to
0^{+}}{N_{\Omega}(\lambda_{j}-h)}<j\leq\lim_{h\to
0^{+}}{N_{\Omega}(\lambda+h)}=\lim_{h\to
0^{+}}{N_{\Omega}(\lambda_{j}+h)}.
\end{displaymath}
La connaissance de $N_{\Omega}(\lambda)$ permet donc de retrouver
les valeurs propres $\lambda_{j}$ et vice-versa.

\noindent Plusieurs résultats pour approximer $N(\lambda)$ ont été
obtenus. Voici certains d'entre eux. Le premier (voir p. 89 de
\cite{Be_2003}) est dû Weyl en 1911 qui a montré que
\begin{Thm}[Weyl]
\begin{displaymath}
N_{\Omega}(\lambda)\sim\frac{V_{n}}{(2\pi)^{n}}\mu(\Omega)\lambda^{\frac{n}{2}}\phantom{12}\textrm{lorsque}\phantom{1}\lambda\to\infty
\end{displaymath}
\end{Thm}
\noindent où $V_{n}$ est le volume de la sphère de rayon 1 dans
$\mathbb{R}^{n}$ et $\mu(\Omega)$ la mesure de Lebesgue du domaine
borné régulier $\Omega$.

\noindent En 1952, 1953 et 1955 Levitan écrivait une série
d'articles où il démontrait le théorème suivant donnant une
estimation du terme d'erreur. Les articles de Levitan sont
respectivement \cite{Le_1952}, \cite{Le_1953} et \cite{Le_1955}. En
1956, Avakumovi\'c montrait également le théorème suivant dans son
article \cite{Av_1956}. Hörmander en donna également une autre
preuve en 1968 dans son article \cite{Ho_1968}.
\begin{Thm}[Levitan]
\begin{displaymath}
N_{\Omega}(\lambda)=\frac{V_{n}}{(2\pi)^{n}}\mu(\Omega)\lambda^{\frac{n}{2}}+O\big(\lambda^{\frac{(n-1)}{2}}\big)
\end{displaymath}
\end{Thm}

\noindent Lorsque Weyl a fait la conjecture suivante qui est vraie
pour tout domaine $\Omega$ convexe et qui a été prouvée par Ivri\u i
dans \cite{Iv_1980} et, indépendamment par Melrose, dans
\cite{Melr_1980} pour d'autres variétés riemanniennes satisfaisant
des hypothèses plus générales.
\begin{Thm}
\begin{displaymath}
N_{\Omega}(\lambda)=\frac{V_{n}}{(2\pi)^{n}}\mu(\Omega)\lambda^{\frac{n}{2}}-\frac{V_{n-1}}{(2\pi)^{n-1}}\mu(\partial\Omega)\lambda^{\frac{n-1}{2}}+o\big(\lambda^{\frac{n-1}{2}}\big)
\end{displaymath}
\end{Thm}

\begin{Rem}[$\triangle{u}+\lambda{u}=0$ ou $\triangle{u}+\lambda^{2}{u}=0$]
Si au lieu d'étudier l'équation $\triangle{u}+\lambda{u}=0$,
l'équation $\triangle{u}+\lambda^{2}{u}=0$ est étudiée, alors il
faut substituer $\lambda$ par $\lambda^{2}$ dans les résultats qui
précèdent.
\end{Rem}

\noindent Comme le disque sera l'objet géométrique d'intérêt dans ce
mémoire, si $\Omega=D$ où $D$ le disque de rayon $r$, alors, dans
\cite{KuFe_1965}, il est montré le théorème suivant spécialisé au
cas du disque. \textbf{Notons que \cite{KuFe_1965} étudie le
problème $\triangle{u}+\lambda^{2}{u}=0$ plutôt que le problème
$\triangle{u}+\lambda{u}=0$. Le théorème est cité presque tel qu'il
est dans l'article afin de faciliter son explication au chapitre 3.2
qui lui est consacré entièrement. Nous comprenons bien ici que
$k=\lambda$ dans le théorème suivant.}
\begin{Thm}[Kuznetsov et Fedosov, cas du disque $D$ de rayon $r$]
Soit le problème à valeur propre $-\triangle{u}=k^{2}{u}$ avec
$u|_{\partial{D}}=0$. Soit la fonction de compte $N_{D}(k)$ des
valeurs $\{k_{n}\}_{n=1}^{\infty}$ ordonnées en ordre croissant
c'est-à-dire soit
\begin{displaymath}
N_{D}(k)=\sum_{k_{n}<k}{1}
\end{displaymath}
alors
\begin{displaymath}
N_{D}(k)=\frac{S}{4\pi}k^{2}-\frac{L}{4\pi}k+O(k^{\frac{2}{3}})
\end{displaymath}
où $S=\pi{r^{2}}$ et où $L=2\pi{r}$.
\end{Thm}

\section{Séparation de variables, équation différentielle de Bessel, fonctions et valeurs propres du Laplacien sur un domaine circulaire}

\noindent Soit $S(\alpha)$, le secteur du disque de rayon 1 d'angle
$\alpha\in(0,2\pi)$ ou soit $D$ le disque de rayon 1. Cette section
montre comment la séparation de variables lors de la résolution de
l'équation $\triangle{u}+\lambda{u}=0$ sur $S(\alpha)$ ou sur $D$
\og engendre \fg{} l'équation différentielle de Bessel. Les
fonctions de Bessel qui sont les solutions de l'équation en question
représentent la partie radiale des fonctions propres du Laplacien
sur un domaine circulaire. Il sera montré que dans le cas du disque
$D$ les fonctions de Bessel impliquées dans la résolution du
problème à valeurs propres $\triangle{u}+\lambda{u}=0$ ont des
ordres entièrs. De même pour le secteur $S(\alpha)$ il sera montré
que les ordres des fonctions de Bessel impliqués ont des ordres
réels de la forme $\frac{n\pi}{\alpha}$ pour $n=1,2,\ldots$

\noindent Pour un exposé plus détaillé de ce qui va suivre dans
cette partie, le lecteur peut se référer entres autres à
~\cite{CoHi1_1953}. En ce qui concerne les fonctions de Bessel plus
spécifiquement, le lecteur peut lire \cite{Wa_1944}.

\noindent \emph{Le cas du secteur $S(\alpha)$ de rayon 1 et d'angle
$\alpha\in(0,2\pi)$ est étudié plus en profondeur dans les lignes
qui suivent. Le cas du disque est presque similaire à celui du
secteur.} Soit donc un secteur $S(\alpha)$ du disque de rayon 1 et
d'angle $\alpha$. Soit l'équation $(\triangle+\lambda)u=0$ avec la
condition de Dirichlet aux frontières c'est-à-dire $u(\partial
S(\alpha))=0$. En coordonnées polaires, l'équation est
\begin{displaymath}
(\triangle+\lambda)u(r,\theta) = \frac{\partial^{2} u}{\partial
r^{2}} + \frac{1}{r}\frac{\partial u}{\partial r} +
\frac{1}{r^{2}}\frac{\partial^{2} u}{\partial \theta^{2}}+\lambda u
=0
\end{displaymath}
En utilisant la séparation de variables avec $u(r,\theta) =
R(r)M(\theta)$ et en substituant,
\begin{displaymath}
\frac{r^{2}\bigg(\frac{d^{2}R}{dr^{2}}+\frac{1}{r}\frac{dR}{dr}+\lambda
R\bigg)}{R}=-\frac{\frac{d^{2}M}{d\theta^{2}}}{M}=\textnormal{constante}=c
\end{displaymath}
Il faut donc résoudre
\begin{displaymath}
M''(\theta)+cM(\theta)=0
\end{displaymath}
où la solution de l'équation précédente est
\begin{displaymath}
M(\theta)=A_{c}\cos(\sqrt{c}\theta)+B_{c}\sin(\sqrt{c}\theta)
\end{displaymath}

\noindent Comme il faut que la fonction $M$ soit périodique
c'est-à-dire $M(0)=M(\alpha)$, alors
\begin{displaymath}
\sqrt{c}=\frac{\pi n}{\alpha}:=:\sqrt{c_{n}}
\end{displaymath}
Puisque $c$ ne dépend que de $n$, il est commode de dénoter par
$M_{c_{n}}$ la solution de l'équation correspondante. Généralement,
\begin{displaymath}
M_{c_{n}}(\theta)=A_{c_{n}}\cos(\sqrt{c_{n}}\theta)+B_{c_{n}}\sin(\sqrt{c_{n}}\theta)\phantom{1},
\end{displaymath}
mais
\begin{eqnarray*}
M(0) & = & M(\alpha)\\
A_{c_{n}}1+B_{c_{n}}0 & = & A_{c_{n}}\cos(n\pi)+B_{c_{n}}0\\
& \Downarrow &\\
A_{c_{n}} & = & 0
\end{eqnarray*}
et, par conséquent,
\begin{displaymath}
M_{c_{n}}(\theta) = B_{c_{n}}\sin(\sqrt{c_{n}}\theta)
\end{displaymath}

\noindent Maintenant, en fixant la valeur $n$ et en dénotant par
$R(r)=R_{n}(r)$ la solution de l'équation
\begin{displaymath}
\frac{r^{2}\bigg(\frac{d^{2}R_{n}}{dr^{2}}+\frac{1}{r}\frac{dR_{n}}{dr}+\lambda_{n}
R_{n}\bigg)}{R_{n}} = c_{n}
\end{displaymath}
Il est démontré que l'équation précédente possède des solutions que
pour un nombre dénombrable de valeurs de $\lambda_{n}$ (voir
~\cite{CoHi1_1953}). En dénotant ces valeurs par
$\{\lambda_{n,k}\}_{k=1}^{\infty}$ et leurs solutions
correspondantes par $R_{n,k}(r)$, alors
\begin{displaymath}
\frac{r^{2}\bigg(\frac{d^{2}R_{n,k}}{dr^{2}}+\frac{1}{r}\frac{dR_{n,k}}{dr}+\lambda_{n,k}
R_{n,k}\bigg)}{R_{n,k}} = c_{n}\phantom{1}.
\end{displaymath}

\begin{Rem}[Fonction de Bessel et zéros des fonctions de Bessel]
La solution $R_{n,k}$ porte le nom de fonction de Bessel de premier
type d'ordre $\frac{n\pi}{\alpha}$. L'ordre sera évident par ce qui
suit. Les valeurs $\{\lambda_{n,k}\}_{k=1}^{\infty}$ seront les
valeurs propres de l'équation originale $\triangle{u}+\lambda{u} =
0$ avec les conditions de Dirichlet aux frontières. Les valeurs
$\{\lambda_{n,k}\}_{k=1}^{\infty}$ correspondent aux \textbf{carrés}
\textbf{des} \textbf{zéros} des fonctions de Bessel qui sont en
quantité dénombrable. Les fonctions de Bessel sont analytiques, ce
qui a du sens puisque toutes les fonctions propres du Laplacien
doivent être au moins lisses.
\end{Rem}

\noindent Pour plus d'information sur les fonctions de Bessel, nous
pouvons lire ~\cite{Wa_1944}. De façon générale, la fonction de
Bessel de premier type est analytique. C'est la solution en série de
l'équation différentielle portant son nom et elle possède le
développement que voici
\begin{eqnarray*}
J_{n}(z)&=&\Big(\frac{1}{2}z\Big)^{n}\sum_{s=0}^{\infty}{\frac{(-1)^{s}\big(\frac{z}{2}\big)^{s}}{s!(n+s)!}}\\
&=&\frac{1}{\pi}\int_{0}^{\pi}{\cos\big(n\theta-z\sin(\theta)\big)d\theta}\phantom{12}\textrm{pour}\phantom{1}n\in\mathbb{Z}\phantom{1}\textrm{et}\phantom{1}z\in\mathbb{C}
\end{eqnarray*}
Lorsque $n$ est remplacé $\nu\in\mathbb{C}$, alors $J_{\nu}(z)$ est
définie par la série suivante
\begin{displaymath}
J_{\nu}(z)=\Big(\frac{1}{2}z\Big)^{\nu}\sum_{s=0}^{\infty}{\frac{(-1)^{s}\big(\frac{z}{2}\big)^{s}}{s!(\nu+s+1)!}}\phantom{12}\textrm{pour}\phantom{1}\nu\in\mathbb{C}\phantom{1}\textrm{et}\phantom{1}z\in\mathbb{C}.
\end{displaymath}

\noindent Soit donc $J_{\nu}(z)$ la fonction de Bessel de premier
type d'ordre réel $\nu$ (l'ordre comme l'argument peuvent être
complexes de façon générale) et d'argument réel $z$. Soit également
le k\ieme{} zéros de $J_{\nu}(z)$ noté par $j_{k}(\nu)$. Ainsi,
\begin{eqnarray*}
R_{n,k}(r)&=&J_{\frac{n\pi}{\alpha}}(\sqrt{\lambda_{n,k}}r)\\
&=&J_{\frac{n\pi}{\alpha}}\Big(j_{k}\Big(\frac{n\pi}{\alpha}\Big)r\Big)\\
&\Downarrow&\\
u_{n,k}(r,\theta)&=&J_{\frac{n\pi}{\alpha}}(\sqrt{\lambda_{n,k}}r)\sin\Big(\frac{n\pi\theta}{\alpha}\Big)
\end{eqnarray*}
Pour vérifier que $u_{n,k}(r,\theta)$ est bien la bonne solution, il
est commode pour l'instant seulement de réécrire
\begin{eqnarray*}
\lambda_{n,k}&:=:&\lambda\\
\frac{n\pi}{\alpha}&:=:\nu\phantom{1}.
\end{eqnarray*}
Voici deux faits concernant les fonctions de Bessel de premier type:
\begin{eqnarray}
J'_{\nu}(z) & = & \frac{\nu}{z}J_{\nu}(z) - J_{\nu+1}(z)\\
J_{\nu+2}(z) & = &
\bigg(\frac{2(\nu+1)}{z}J_{\nu+1}(z)-J_{\nu}(z)\bigg).
\end{eqnarray}
Par conséquent,
\begin{eqnarray*}
(J_{\nu}(\sqrt{\lambda} r))' & = & J'_{\nu}(\sqrt{\lambda} r)(\sqrt{\lambda} r)'\\
& = & \sqrt{\lambda} J'_{\nu}(\sqrt{\lambda} r)\\
& = & \sqrt{\lambda} \bigg(\frac{\nu}{\sqrt{\lambda} r}J_{\nu}(\sqrt{\lambda} r)-J_{\nu+1}(\sqrt{\lambda} r)\bigg)\\
& = & \frac{\nu}{r}J_{\nu}(\sqrt{\lambda} r)-\sqrt{\lambda}
J_{\nu+1}(\sqrt{\lambda} r).
\end{eqnarray*}
Aussi,
\begin{eqnarray*}
(J_{\nu}(\sqrt{\lambda} r))'' & = & (J'_{\nu}(\sqrt{\lambda} r))'\\
& = & \bigg(\frac{\nu}{r}J_{\nu}(\sqrt{\lambda} r)-\sqrt{\lambda} J_{\nu+1}(\sqrt{\lambda} r)\bigg)'\\
& = & \nu\bigg(\frac{rJ'_{\nu}(\sqrt{\lambda}
r)-J_{\nu}(\sqrt{\lambda}
r)}{r^{2}}\bigg)-\sqrt{\lambda}\big(J_{\nu+1}(\sqrt{\lambda} r)\big)'\\
& \ldots & \\
& =
&\bigg(\frac{\nu^{2}}{r^{2}}-\frac{\nu}{r^{2}}\bigg)J_{\nu}(\sqrt{\lambda}
r)\ldots\\
&&\phantom{1234567890}-\sqrt{\lambda}\bigg(\frac{\nu}{r}+\frac{\nu+1}{r}\bigg)J_{\nu+1}(\sqrt{\lambda}
r) + \lambda J_{\nu+2}(\sqrt{\lambda} r)\phantom{1},
\end{eqnarray*}
ce qui donne
\begin{eqnarray*}
(\triangle + \lambda)J_{\nu}(\sqrt{\lambda} r)\sin(\nu\theta) & = &
\sin(\nu\theta)\bigg(\frac{\nu^{2}}{r^{2}}-\frac{\nu}{r^{2}}+\frac{\nu}{r^{2}}\bigg)J_{\nu}(\sqrt{\lambda}r)\\
& &-\sqrt{\lambda}\sin(\nu\theta)\bigg(\frac{\nu}{r}+\frac{\nu+1}{r}+\frac{1}{r}\bigg)J_{\nu+1}(\sqrt{\lambda} r)\\
& &+\sin(\nu\theta)\lambda J_{\nu+2}(\sqrt{\lambda} r)\\
& &-\frac{\nu^{2}}{r^{2}}\sin(\nu\theta)J_{\nu}(\sqrt{\lambda} r)\\
& &+\lambda\sin(\nu\theta)J_{\nu}(\sqrt{\lambda} r)\\
\end{eqnarray*}
qui implique que
\begin{eqnarray*}
(\triangle + \lambda)J_{\nu}(\sqrt{\lambda} r)\sin(\nu\theta) & = &
\bigg(\frac{\nu^{2}}{r^{2}}-\frac{\nu}{r^{2}}+\frac{\nu}{r^{2}}-\frac{\nu^{2}}{r^{2}}+\lambda\bigg)J_{\nu}(\sqrt{\lambda}
r)\sin(\nu\theta)\\
& &
-\sqrt{\lambda}\bigg(\frac{\nu}{r}+\frac{\nu+1}{r}+\frac{1}{r}\bigg)J_{\nu+1}(\sqrt{\lambda}r)\sin(\nu\theta)\\
& &+\lambda J_{\nu+2}(\sqrt{\lambda} r)\sin(\nu\theta)\\
& = &
\bigg(\frac{\nu^{2}}{r^{2}}-\frac{\nu}{r^{2}}+\frac{\nu}{r^{2}}-\frac{\nu^{2}}{r^{2}}+\lambda\bigg)J_{\nu}(\sqrt{\lambda}
r)\sin(\nu\theta)\\
& &
-\sqrt{\lambda}\bigg(\frac{\nu}{r}+\frac{\nu+1}{r}+\frac{1}{r}\bigg)J_{\nu+1}(\sqrt{\lambda}
r)\sin(\nu\theta)\\
& &+\lambda\bigg(\frac{2(\nu+1)}{\sqrt{\lambda}
r}J_{\nu+1}(\sqrt{\lambda} r) - J_{\nu}(\sqrt{\lambda} r)\bigg)\sin(\nu\theta)\\
\end{eqnarray*}
qui implique que
\begin{eqnarray*}
(\triangle + \lambda)J_{\nu}(\sqrt{\lambda} r)\sin(\nu\theta) & = &
\bigg(\frac{\nu^{2}}{r^{2}}-\frac{\nu}{r^{2}}+\frac{\nu}{r^{2}}-\frac{\nu^{2}}{r^{2}}+\lambda-\lambda\bigg)J_{\nu}(\sqrt{\lambda}
r)\sin(\nu\theta)\\
& & +
\bigg(\frac{2\nu\sqrt{\lambda}}{r}+\frac{2\sqrt{\lambda}}{r}-\frac{\sqrt{\lambda}\nu}{r}\ldots\\
&&\phantom{1234567890}-\frac{\sqrt{\lambda}\nu}{r}-\frac{2\sqrt{\lambda}}{r}\bigg)J_{\nu+1}(\sqrt{\lambda}
r)\sin(\nu\theta)\\
\end{eqnarray*}
et ainsi
\begin{displaymath}
(\triangle + \lambda)J_{\nu}(\sqrt{\lambda} r)\sin(\nu\theta) = 0
\end{displaymath}

\begin{Rem}[Multiplicité des valeurs propres pour un secteur]
Dans ce cas, toutes les valeurs propres $\lambda_{n,k}$ pour $n\geq
1$ et $k\geq 1$ du Laplacien sur un secteur du disque sont simples
en général sauf pour des angles critiques (voir entre autres la
remarque \ref{RemAngleCrit}). Les fonctions de Bessel d'ordre 0 ne
sont jamais des fonctions propres du Laplacien sur un secteur.
\end{Rem}

\begin{Rem}[Multiplicité des valeurs propres pour le disque]
Soit $D$, le disque de rayon 1 dans $\mathbb{R}^{2}$. Les fonctions
propres sont encore des fonctions de Bessel. Cependant les fonctions
de Bessel d'ordre 0 sont admises. De façon générale,
\begin{eqnarray*}
u_{0,k}(r,\theta)&=&J_{0}(\sqrt{\lambda_{0,k}}r)\phantom{12}k\geq 1\\
u_{n,k}(r,\theta)&=&J_{n}(\sqrt{\lambda_{n,k}}r)\sin(n\theta)\phantom{12}\textrm{ou}\\
u_{n,k}(r,\theta)&=&J_{n}(\sqrt{\lambda_{n,k}}r)\cos(n\theta)\phantom{12}\textrm{pour
}n\geq1\textrm{ et }k\geq1
\end{eqnarray*}
La multiplicité du cas correspondant à l'ordre $n=0$ est ainsi
simple et celles correspondant aux ordres $n\geq 1$ sont doubles.
\end{Rem}

\section{Développement asymptotique d'Olver des zéros des fonctions de Bessel}


\noindent Les zéros de fonctions de Bessel sont les valeurs propres
du Laplacien $\triangle$ sur un domaine circulaire, le disque ou un
secteur du disque, comme il a été discuté précédemment.

\noindent Plusieurs développements asymtotiques des zéros des
fonctions de Bessel et des zéros des dérivées des fonctions de
Bessel existent. Par exemples, les développements d'Olver, de
McMahon ou de Meissel (pour $J_{n}(n)$) pour ne citer que ces
développements les plus connus. Ces développements sont utiles selon
que l'ordre ou l'argument sont dans une certaine région du plan
complexe. Ils existent par exemple des développement de Taylor pour
les zéros de fonctions de Bessel d'ordre $\frac{2k+1}{2}$ pour
$k=0,1,2,\ldots$ ou simplement pour les ordres entiers. Le lecteur
peut lire \cite{Rog_2005} en ce qui concerne les développements de
Taylor reliés aux zéros des fonctions de Bessel où une liste de
problèmes ouverts est également présentée. Le développement qui est
sans doute le plus utile est celui fourni par F.W.J. Olver durant
les années 50. Le développement asymptotique d'Olver a permis en
autre une implémentation numérique stable dans plusieurs logiciels.
Tout lecteur intéressé peut lire les nombreux articles d'Olver cités
dans les références comme \cite{Ol_1951}, \cite{Ol_1952},
~\cite{Ol_1960}, ~\cite{Ol_1974} et ~\cite{FLO_2004}.

\noindent Avant de donner le développement asymptotique d'Olver, une
brève introduction aux fonctions d'Airy s'impose. Le lecteur peut se
référer à \cite{Ol_1974}, ce dernier ouvrage traite des fonctions
\og spéciales \fg{} comme les fonctions d'Airy et de Bessel. Une
fonction d'Airy est une fonction qui est une solution de l'équation
différentielle
\begin{displaymath}
\frac{d^{2}w}{dt^{2}}-tw=0.
\end{displaymath}
\noindent L'équation précédente a deux solutions dénotées par $Ai$
et $Bi$ qui sont nommées respectivement fonction d'Airy de premier
type et fonction d'Airy de deuxième type.

\noindent Les fonctions d'Airy $Ai$ et $Bi$ sont reliées aux
fonctions de Bessel de premier type $J_{\nu}$ avec l'ordre
$\nu=\frac{1}{3}$ comme suit
\begin{eqnarray*}
Ai(t)&=&\frac{1}{3}\sqrt{x}\bigg(J_{-\frac{1}{3}}\Big(\frac{2}{3}t^{\frac{3}{2}}\Big)-J_{\frac{1}{3}}\Big(\frac{2}{3}t^{\frac{3}{2}}\Big)\bigg)\phantom{12}t>0\\
Bi(t)&=&\sqrt{\frac{x}{3}}\bigg(J_{-\frac{1}{3}}\Big(\frac{2}{3}t^{\frac{3}{2}}\Big)+J_{\frac{1}{3}}\Big(\frac{2}{3}t^{\frac{3}{2}}\Big)\bigg)\phantom{12}t>0
\end{eqnarray*}
\noindent Les relations précédentes sont connues sous le nom de \og
connection formulas \fg{}. Les zéros de $Ai$ seront dénotés par
$a_{k}$ avec $k=1,2,\ldots$ et $0>a_{1}>a_{2}>\ldots$

\noindent Seule la fonction $Ai$ sera utile pour développer les
zéros de la fonction de Bessel $J_{\nu}$ d'ordre $\nu>0$. Il est à
noter que les zéros de $Ai$ et de $Bi$ sont tous négatifs. $a_{k}$
possède également un développement asymptotique d'Olver (voir
\cite{Ol_1974}). Il est possible de démonter en utilisant le
développement asymptotique que
\begin{displaymath}
\lim_{k\to\infty}{\frac{-a_{k}}{\big(\frac{3}{8}\pi(4k-1)\big)^{\frac{2}{3}}}}
= 1.
\end{displaymath}

\noindent Pour la suite, nous nous concentrerons sur les zéros de la
fonction de Bessel de premier type $J_{\nu}$ pour $\nu>0$. En
dénotant par $j_{\nu}(k)$ le k\ieme{} zéro de $J_{\nu}(z)$ pour
$z>0$ et en dénotant par $\gamma_{k}=-a_{k}2^{-\frac{1}{3}}$ alors
\begin{eqnarray*}
j_{k}(\nu)&=&\nu+\gamma_{k}\nu^{\frac{1}{3}}+\frac{3}{10}\gamma_{k}^{2}\nu^{-\frac{1}{3}}+\frac{5-\gamma_{k}^{3}}{350}\nu^{-1}\\
&
&-\frac{479\gamma_{k}^{4}+20\gamma_{k}}{63000}\nu^{-\frac{5}{3}}+\frac{20231\gamma_{k}^{5}-27550\gamma_{k}^{2}}{8085000}\nu^{-\frac{7}{3}}+\mathcal{O}(\nu^{-3}).
\end{eqnarray*}



\noindent La preuve de ce résultat est longue et complexe. Les idées
de base se retrouvent principalement dans les deux longs articles
~\cite{Ol_1951} et ~\cite{Ol_1952}.




\noindent Il est possible de constater grâce entres autres au
développement ci-dessus de $j_{\nu}(k)$ que
\begin{eqnarray}\label{capoute1}
\lim_{\nu\to\infty}\frac{j_{\nu}(k)}{\nu}&=&1\phantom{12}\forall{k}\in\mathbb{N}.
\end{eqnarray}
\noindent Il est même prouvé pour $\nu>0$ dans \cite{QuWo_1999} que
\begin{eqnarray}\label{capoute2}
\nu-\frac{a_{k}}{2^{\frac{1}{3}}}\nu^{\frac{1}{3}}<j_{k}(\nu)<\nu-\frac{a_{k}}{2^{\frac{1}{3}}}\nu^{\frac{1}{3}}+\frac{3}{20}a_{k}^{2}\frac{2^{\frac{1}{3}}}{\nu^{\frac{1}{3}}}.
\end{eqnarray}
\noindent Grâce à (\ref{capoute1}) et à (\ref{capoute2}), il est
possible d'évaluer l'ordre de grandeur des zéros. Ceci sera utile
ultérieurement (voir les remarques \ref{remcapoute1} et
\ref{remprecalgocompte}) afin de connaître le degré de précision
requis pour évaluer numériquement les zéros des fonctions de Bessel.

\section{Les ensembles nodaux des fonctions propres et leurs
structures, cas du disque et du secteur en exemples}

\noindent Le théorème principal de cette section concerne le nombre
de domaines nodaux qu'une fonction propre peut avoir. Il est énoncé
entre autres dans \cite{Sh_abcd}. Avant d'arriver au théorème, un
rappel de ce qu'est l'ensemble nodal d'une fonction propre s'impose.
De même, les exemples pour le secteur du disque $S(\alpha)$ de rayon
1 et d'angle $\alpha\in(0,2\pi)$ et le disque $D$ de rayon 1 sont
donnés.

\noindent L'ensemble nodal d'une fonction propre est son ensemble de
zéro. Si la dimension du domaine est $d$, alors la dimension des
ensemles nodaux est $d-1$. En d'autres termes, pour la j\ieme{}
fonction propre $\varphi_{j}$ défini sur un domaine ouvert assez
régulier $\Omega\subset\mathbb{R}^{d}$, si $\mathrm{Z}_{j}$ est
l'ensemble des zéros, alors
\begin{displaymath}
\mathrm{Z}_{j}=\Big\{\varphi_{j}(x)=0,\phantom{12}x\in\Omega\Big\}.
\end{displaymath}
$Z_{j}$ est un ensemble fermé.

\noindent Dans le cas d'un secteur $S(\alpha)$ de \textbf{rayon} 1
et d'angle $\alpha\in(0,2\pi)$, une fonction propre typique est
donnée par
\begin{displaymath}
u_{n,k}(r,\theta)=J_{\frac{n\pi}{\alpha}}(\sqrt{\lambda_{n,k}}r)\sin\Big(\frac{n\pi\theta}{\alpha}\Big)
\end{displaymath}
Par conséquent,
\begin{displaymath}
u_{n,k}(r,\theta)=0\Leftrightarrow
J_{\frac{n\pi}{\alpha}}(\sqrt{\lambda_{n,k}}r)=0\phantom{12}\textrm{ou}\phantom{12}
\sin\Big(\frac{n\pi\theta}{\alpha}\Big)=0
\end{displaymath}
Ainsi
\begin{eqnarray*}
\sin\Big(\frac{n\pi\theta}{\alpha}\Big)&=&0\\
&\Updownarrow&\\
\frac{n\pi\theta}{\alpha}&=&d\pi\phantom{12}\textrm{pour
$d=0,1,\ldots,n$}\\
&\Updownarrow&\\
\theta&=&\frac{d\alpha}{n}
\end{eqnarray*}
lequel cas donne une suite de lignes angulaires décrites par les
multiples entiers $d=0,1,\ldots, n$ de l'angle $\frac{\alpha}{n}$.
Le cas où
\begin{eqnarray*}
J_{\frac{n\pi}{\alpha}}(\sqrt{\lambda_{n,k}}r)&=&0\\
&\Updownarrow&\\
r&=&\frac{\sqrt{\lambda_{n,j}}}{\sqrt{\lambda_{n,k}}}\phantom{12}\textrm{pour }j=1,\ldots k\\
\end{eqnarray*}
donne une suite de cercles concentriques avec le cas correspondant à
la frontière lorsque $j=k$, c'est-à-dire lorsque $r=1$.

\noindent L'ensemble des zéros d'une fonction propre découpe le
domaine de définition des fonctions propres en composantes connexes
où la fonction propre alterne de signe. Dans le cas du disque ou
d'un secteur circulaire, les composantes connexes ont la forme de
secteurs circulaires comme il est possible de le voir sur les
figures de la section \ref{Pics}. Une image vaut milles mots.

\begin{Rem}[Nombre de composantes connexes]$\phantom{1}$\\
(\textbf{cas d'un secteur})\\
Pour une fonction propre donnée $u_{n,k}(r,\theta)$, le nombre de
composantes connexes sur lesquelles la fonction propre alterne de
signe est donnée par $nk$.
\\(\textbf{cas du disque})\\
Pour les fonctions propres sur le disque de multiplicité 1
c'est-à-dire pour celles de la forme $u_{0,k}(r,\theta)$, le nombre
de composantes connexes est $k$. Pour les fonctions propres de
multiplicité 2 c'est-à-dire celles de la forme $u_{n,k}(r,\theta)$
avec $n=1,2,\ldots$, le nombre de composantes connexes est $2nk$.
\end{Rem}

\noindent Avant de conclure cette section, voici un théorème célèbre
dû à Courant qui sera très utile pour la prochaine section qui peut
être retrouvé entre autres dans \cite{Sh_abcd}.

\begin{Thm}[Théorème de Courant sur les ensembles nodaux]
Soit $\Omega\subset\mathbb{R}^{n}$, un domaine tel que
$\partial\Omega$ est suffisamment régulier par morceaux. Soit
$\lambda_{m}$, la m\ieme{} valeur propre. Alors pour une fonction
propre $u_{m}$ correspondant à $\lambda_{m}$, le nombre de domaines
nodaux est au plus $m$.
\end{Thm}

\chapter[]{Les ensembles nodaux des fonctions propres}

\noindent Soit $D$ ou $S(\alpha)$, le disque de rayon 1 ou un
secteur de rayon 1 et d'angle $\alpha\in (0,2\pi)$. Soit $\triangle$
l'opérateur différentiel de Laplace. Soit $\lambda$, une valeur
propre $\triangle$ telle que $(\triangle+\lambda)u=0$ avec $u\in
C^{\infty}$ étant la fonction propre associée à $\lambda$.

\noindent Dans un premier temps, une conjecture concernant la
1\iere{} ligne nodale du problème de Dirichlet sur un domaine
bidimensionnel est étudiée. La conjecture en question a été soulevée
par Payne et affirme que pour tout domaine bidimensionnel, la
1\iere{} ligne nodale touche toujours la frontière. Deux cas
extrêmes s'offrent à la conjecture, le rectangle et le secteur
circulaire. Le cas du rectangle étant très simple puisqu'il est
facile d'identifier où se trouve exactement la première ligne nodale
comme il sera montré plus loin. Le cas du secteur étant un peu plus
intéressant, cela motivant en soi la prochaine partie. À la fin de
la 1\iere{} section, le lecteur trouvera une proposition donnant
l'angle critique $\alpha_{0}$ pour lequel tout secteur d'angle
$\alpha<\alpha_{0}$ a une ligne nodale en coordonnées polaires de la
forme $r=\mathrm{const}$ et pour lequel tout secteur d'angle
$\alpha>\alpha_{0}$ a une ligne nodale en coordonnées polaires de la
forme $\theta=\mathrm{const}$. De même, lorsque $\alpha=\alpha_{0}$,
nous verrons que la ligne nodale est indéterminée.

\noindent Dans un deuxième temps, soit un entier $m>1$ généralement
très grand comme $m=1\cdot{10}^{6}$ par exemple. Comment construire
la suite des valeurs propres $\{\lambda_{j}\}_{j=1}^{m}$ en ordre
croissant de l'opérateur $\triangle$ sur $D$ ou $S(\alpha)$? En
d'autres termes, il faut établir un algorithme efficace pour
déterminer le spectre jusqu'au rang $m$. Par \textbf{efficace}, cela
veut dire ne pas résoudre l'équation différentielle avec des
méthodes comme les éléments finis qui ne permettraient pas de
construire la suite spectrale pour des valeurs très grandes comme
$m=1\cdot{10}^{6}$ et ni même d'identifier les ensembles nodaux de
façon exacte. Par \textbf{efficace}, il faut se servir de toute la
théorie décrite jusqu'à maintenant et du fait que dans certains
logiciels, en l'occurence Maple, il existe des routines très stables
permettant de trouver les zéros des fonctions de Bessel. La
construction de la suite spectrale jusqu'à $m$ est équivalente à
dénombrer le nombre de composantes connexes disjointes obtenues pour
chaque $u_{j}$ en enlevant du domaine original les ensembles nodaux
des fonctions propres $u_{j}$ pour $j=1,2,\ldots,m$. Les valeurs
propres sont données en fonction de deux paramètres $n$ et $k$
entiers non-négatifs et connaître le spectre jusqu'à l'entier $m$,
consistera à établir une bijection entre $\mathbb{N}^{2}$ et
$\mathbb{N}$.

\section{La première ligne nodale de la deuxième fonction propre}

\noindent Dans cette partie, la première ligne nodale de la deuxième
fonction propre de $\triangle$ sur un domaine convexe du plan est
étudiée. Cette ligne est l'ensemble des zéros de la deuxième
fonction propre $u_{2}$ du Laplacien d'un domaine convexe du plan,
avec des conditions de Dirichlet sur la frontière. Pour la première
fonction propre, $u_{1}$, elle ne possède qu'une seule ligne nodale
qui s'annule sur la frontière du domaine en question et qui, par
conséquent, est triviale. Quant à $u_{2}$, par le théorème de
Courant sur les ensembles nodaux, elle possède au plus 2 lignes; la
première étant toujours triviale et s'annulant sur la frontière et
quant à la deuxième ligne, non triviale celle-ci, elle peut toucher
à la frontière ou être entièrement contenue à l'intérieur du
domaine.

\noindent Mathématiquement parlant, soit $\Omega$ un domaine convexe
du plan. Soit $u_{2}$, la deuxième fonction propre de $\triangle$,
c'est-à-dire $\triangle{u_{2}}+\lambda_{2}u_{2}=0$ et
$u_{2}(\partial\Omega)=0$. Soit $\Gamma=\{z\in\Omega,u_{2}(z)=0\}$,
alors deux possibilités s'imposent:
$\bar{\Gamma}\cap\partial\Omega=\emptyset$ ou
$\bar{\Gamma}\cap\partial\Omega\neq\emptyset$. Dans le cas où
$\bar{\Gamma}\cap\partial\Omega=\emptyset$, alors
$\bar{\Gamma}\subsetneq\Omega$.

\noindent L.E. Payne a conjecturé que
$\bar{\Gamma}\cap\partial\Omega\neq\emptyset$ pour tout domaine du
plan. S.-T. Yau a conjecturé que
$\bar{\Gamma}\cap\partial\Omega\neq\emptyset$ pour tout domaine
convexe du plan. Pour un domaine borné convexe
$\Omega\subset\mathbb{R}^{2}$, Melas en 1992 dans \cite{Mela_1992} a
montré le théorème suivant.
\begin{Thm}[Melas]
Soit $\Omega\subset\mathbb{R}^{2}$ un domaine borné convexe avec des
frontières C$^{\infty}$, alors la ligne nodale de la deuxième
fonction propre $u_{2}$ doit intersecter la frontière
$\partial\Omega$ exactement en deux points.
\end{Thm}

\noindent Pour un secteur $S(\alpha)$ de rayon 1 qui est un domaine
borné convexe du plan, la ligne nodale de la deuxième fonction
propre peut donc toucher $\partial{S(\alpha)}$ de deux façons. La
ligne nodale intersecte la partie de $\partial{S(\alpha)}$ décrite
par $r=1$ ou les deux parties décrites par $\theta=0$ et
$\theta=\alpha$.

\noindent \textbf{Question.} Quelle doit être la valeur de $\alpha$
faisant en sorte que nous avons une situation plutôt qu'une autre?

\begin{Rem}[Deux possibilités pour $u_{2}$]\label{RemPfProp}
$u_{2}$ a donc deux sous domaines simplement connexes disjoints
séparés par une ligne nodale, radiale ou angulaire. Les seuls zéros
admis des fonctions de Bessel sont
$j_{1}\big({\frac{\pi}{\alpha}}\big)$,
$j_{1}\big({2\frac{\pi}{\alpha}}\big)$ et
$j_{2}\big({\frac{\pi}{\alpha}}\big)$. Puisque
$\lambda_{1}=j_{1}\big({\frac{\pi}{\alpha}}\big)$ nécessairement,
alors les deux seules autres possibiltés sont
$j_{1}\big({2\frac{\pi}{\alpha}}\big)$ et
$j_{2}\big({\frac{\pi}{\alpha}}\big)$ chacun donnant respectivement
une ligne nodale angulaire et une ligne nodale radiale et, dans
chaque cas, deux sous domaines simplement connexes où $u_{2}$
alterne de signe.
\end{Rem}

\noindent La fonction $j_{k}(\nu)$ avec $\nu=\frac{\pi}{\alpha}$
étant lisse et croissante avec l'ordre $\nu$, alors la ligne nodale
radiale se produit si et seulement si
\begin{displaymath}
j_{2}^{2}(\nu) < j_{1}^{2}(2\nu) \Leftrightarrow j_{2}(\nu) <
j_{1}(2\nu)
\end{displaymath}
et l'équation de la ligne radiale est donnée par
\begin{displaymath}
r = \frac{j_{1}^{2}(\nu)}{j_{2}^{2}(\nu)} =
\frac{\lambda_{1}(\alpha)}{\lambda_{2}(\alpha)} \phantom{1}
\textrm{avec}\phantom{1} \nu = \frac{\pi}{\alpha}.
\end{displaymath}
Quant à la ligne nodale angulaire, elle se produit si et seulement
si
\begin{displaymath}
j_{2}^{2}(\nu) > j_{1}^{2}(2\nu) \Leftrightarrow j_{2}(\nu) >
j_{1}(2\nu)
\end{displaymath}
et l'équation de la ligne angulaire est donnée par
\begin{displaymath}
\theta=\frac{\alpha}{2}.
\end{displaymath}

\noindent Lorsque $\nu=\frac{\pi}{\alpha}$ satisfait l'équation
$j_{1}(2\nu)=j_{2}(\nu)$, alors l'angle $\alpha$ est critique et la
deuxième valeur propre a une multiplicité double. Si nous résolvons
l'équation
\begin{displaymath}
j_{2}(\nu) - j_{1}(2\nu) = 0,
\end{displaymath}
alors la solution est
\begin{displaymath}
\nu \approx 2.823823 \Leftrightarrow \alpha \approx 1.112531 \approx
63.743334\phantom{1}\textrm{degrées}
\end{displaymath}
et la 1\iere{} ligne nodale n'est donc pas définie. En effet, soit
\begin{eqnarray}
u_{1,2}&=&J_{\frac{\pi}{\alpha_{0}}}\Big(j_{2}\Big(\frac{\pi}{\alpha_{0}}\Big)r\Big)\sin\Big(\frac{\pi\theta}{\alpha_{0}}\Big)\label{2ndNoLiEF1}\\
u_{2,1}&=&J_{\frac{2\pi}{\alpha_{0}}}\Big(j_{1}\Big(2\frac{\pi}{\alpha_{0}}\Big)r\Big)\sin\Big(\frac{2\pi\theta}{\alpha_{0}}\Big)\label{2ndNoLiEF2}
\end{eqnarray}
les deux fonctions propres associées respectivement aux zéros
$j_{2}\big(\frac{\pi}{\alpha_{0}}\big)$ et
$j_{1}\big(2\frac{\pi}{\alpha_{0}}\big)$. Comme
\begin{displaymath}
\lambda_{2}=j_{2}^{2}\Big(\frac{\pi}{\alpha_{0}}\Big)=j_{1}^{2}\Big(2\frac{\pi}{\alpha_{0}}\Big)
\end{displaymath}
alors n'importe quelles combinaisons linéaires des fonctions propres
(\ref{2ndNoLiEF1}) et (\ref{2ndNoLiEF2}) est une solution de
l'équation
$\triangle(u)+\lambda{u}=0$ sur $S(\alpha_{0})$.\\

\begin{Rem}[1\iere{} ligne nodale et secteur
critique]\label{RemAngleCrit} Étant donné un angle $\alpha$ et le
secteur du disque correspondant d'angle $\alpha$, les valeurs
propres sont généralement de multiplicité 1. Soit l'angle critique
$\alpha_{0}$ et l'ordre associé $\nu_{0}=\frac{\pi}{\alpha_{0}}$ tel
que $j_{2}(\nu_{0}) - j_{1}(2\nu_{0}) = 0$. Lorsque l'angle
$\alpha=\alpha_{0}$, la deuxième valeur propre a une multiplicité de
2. La première ligne nodale n'est donc pas définie pour cette valeur
critique $\alpha$ puisque n'importe quelles combinaisons linéaires
des fonctions propres (\ref{2ndNoLiEF1}) et (\ref{2ndNoLiEF2}) est
une solution du problème $\triangle(u)+\lambda{u}=0$ sur
$S(\alpha_{0})$.
\end{Rem}

\noindent Il est intéressant de visualiser les ensembles nodaux de
la 2\ieme{} valeur propre pour quelques valeurs $\alpha$ ou, en
d'autres termes, la 1\iere{} ligne nodale. Nous voyons en effet que
lorsque $\alpha\to\alpha_{0}$, la ligne nodale est indéterminée, en
effet, Matlab prend au \og hasard \fg{} une combinaison linéaire des
deux fonctions propres. \textbf{Les angles mentionnés en-dessous de
chaque graphique sont exactes et ne sont pas approximatives. Le but
étant de donner à Matlab des angles se rapprochant de plus en plus
près de $\alpha_{0}$.}

\newpage

\begin{figure}[ht]
\centering
\includegraphics[width = 8.5cm, height = 8.5cm]{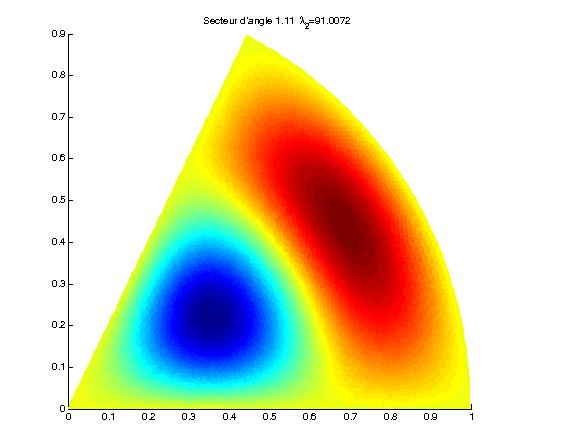}
\caption{1\iere{} ligne nodale, $\alpha=1.11$,
$\lambda_{2}\approx91.0072$}
\end{figure}

\begin{figure}[ht]
\centering
\includegraphics[width = 8.5cm, height = 8.5cm]{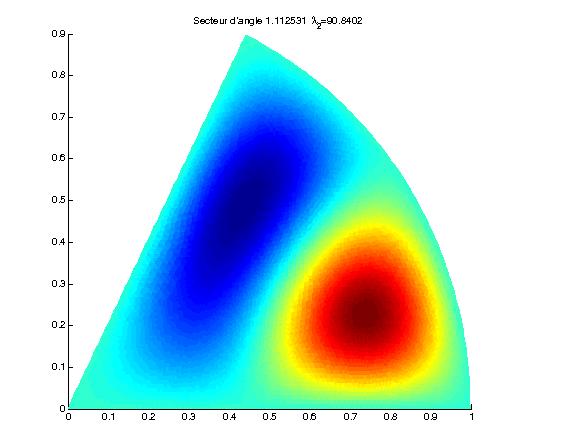}
\caption{1\iere{} ligne nodale, $\alpha=1.112531$,
$\lambda_{2}\approx90.8402$}
\end{figure}

\begin{figure}[ht]
\centering
\includegraphics[width = 8.5cm, height = 8.5cm]{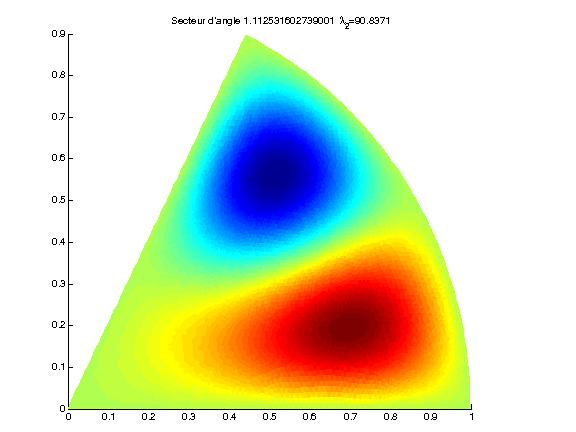}
\caption{1\iere{} ligne nodale, $\alpha=1.112531602739001$,
$\lambda_{2}\approx90.8371$}
\end{figure}

\begin{figure}[ht]
\centering
\includegraphics[width = 8.5cm, height = 8.5cm]{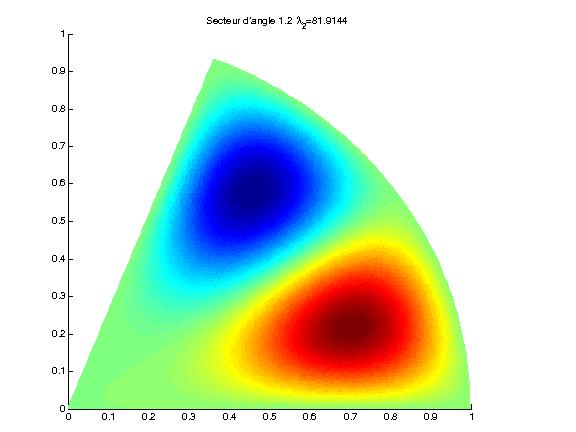}
\caption{1\iere{} ligne nodale, $\alpha=1.2$,
$\lambda_{2}\approx81.9144$}
\end{figure}

\newpage

\noindent Voici donc à la lumière de ce qui précède la petite
proposition suivante concernant l'\og emplacement\fg{} de la
première ligne nodale de la deuxième fonction propre.
\begin{Prop}[\og Emplacement\fg{} de la première ligne nodale]
Soit $S(\alpha)$ un secteur du disque de rayon 1 et d'angle
$\alpha$. Soit $\alpha_{0}$ la solution de l'équation
\begin{displaymath}
j_{2}\Big(\frac{\pi}{\alpha_{0}}\Big)=j_{1}\Big(\frac{2\pi}{\alpha_{0}}\Big)
\end{displaymath}
où $j_{k}(\nu)$ est le k\ieme{} zéro de la fonction de Bessel de
premier type d'ordre $\nu>0$. Alors si l'angle $\alpha$ du secteur
$S(\alpha)$ est supérieur à $\alpha_{0}$, alors la première ligne
nodale de la deuxième fonction propre est de la forme
$\theta=\frac{\alpha}{2}$ sinon elle de la forme
$r=\frac{j_{1}^{2}(\nu)}{j_{2}^{2}(\nu)}$ avec
$\nu=\frac{\pi}{\alpha}$.
\end{Prop}
\begin{Dem}
Lire la remarque \ref{RemPfProp}.
\end{Dem}

\section{Algorithme pour ordonner les valeurs propres et déterminer la structure des ensembles nodaux}

\noindent \textbf{Question} Supposons donné un entier $m>1$. Est-il
possible de déterminer efficacement le spectre du Laplacien sur le
disque $D$ ou un secteur $S(\alpha)$ de rayon 1 jusqu'à l'entier m?
En d'autres termes, comment obtenir toutes les valeurs propres
$\lambda_{j}$ pour $j=1,\ldots,m$ telles que
\begin{displaymath}
0<\lambda_{1}\leq\lambda_{2}\leq\ldots\leq\lambda_{m}
\end{displaymath}

\noindent \textbf{Réponse} Répondre à la question précédente, c'est
connaître la bijection entre $j_{k}^{2}(\nu)$ et $\lambda_{m}$ où
$\nu=n$ pour le disque $D$ et où $\nu=\frac{n\pi}{\alpha}$ pour le
secteur $S(\alpha)$. Nous avons donc un problème de classification
des valeurs $j_{\nu}^{2}(k)$. Étant donné $m$, quelles valeurs de
$n$ et $k$ sont admissibles en fonctions de $m$? Dans ce qui suit,
nous écrirons $\nu=\nu(n)$.

\noindent \textbf{Idée de l'algorithme} Par le théorème de Courant,
il faut trouver toutes les paires $(n,k)\in\mathbb{N}^{2}$ telles
que $nk\leq m$, ce sont les paires admissibles. De même, la
connaissance de toutes ces paires $(n,k)$ implique la connaissance
des paires $(n',k')$ telles que $n'k'\leq m'$ pour $m'\leq m$. Par
conséquent, évaluer $j_{\nu(n)}^{2}(k)$ pour toutes les paires
$(n,k)$ telles que $nk\leq m$ assure d'obtenir
$\{\lambda_{j}\}_{j=1}^{m}\subsetneqq \{j_{\nu(n)}^{2}(k)\}_{nk\leq
m}$. Il ne reste donc qu'à ordonner et à prendre en compte la
multiplicité dans le cas du disque. Dans le cas de $S(\alpha)$, il
ne faut qu'ordonner linéairement $\{j_{\nu(n)}^{2}(k)\}_{nk\leq m}$
et prendre les $m$ premières valeurs, ce qui donne
$0<\lambda_{1}\leq\ldots\leq\lambda_{m}$. Dans le cas du disque, il
faut regarder les valeurs où $n=0$ et celles $n\neq 0$ et extraire
les valeurs de $\{j_{\nu(n)}^{2}(k)\}_{nk\leq m}$ en répétant deux
fois celles où $n\neq0$ pour obtenir
$0<\lambda_{1}\leq\ldots\leq\lambda_{m}$.


\noindent L'algorithme précédent a été implémenté en Maple et est
donné en annexe.

%
%

\begin{Rem}[Améliorations possibles]
Comme à n'importe quelle programme informatique, plusieurs
améliorations peuveut être apportées.\\
1) La première, pour ceux dont les ressouces le permettent, seraient
de
paralléliser le code avec MPI.\\
2) Afin de réduire la taille en mémoire utilisée et également de
permettre de sauvegarder en mode binaire (car Maple ne fait pas la
distinction entre le mode binaire et ascii sur les stations Unix),
reprogrammer en C++ avec par exemple la librairie NTL (Number Theory
Library) compilée avec GMP (Gnu Multiple Precision arithmetic)
procureait sans aucun doute des économies de temps et d'espace. Cela
implique cependant d'implémenter une fonction pour le calcul des
zéros de fonctions Bessel de premier type qui serait stable, ce qui
n'est pas nécessairement évident à faire efficacement.
\end{Rem}

\begin{Rem}[Avantages de l'algorithme]
À la connaissance de l'auteur, l'algorithme donné ici est un des
plus efficace pour déterminer les valeurs propres et les fonctions
propres de l'opérateur $\triangle$ sur le disque. Il ne requiert
aucunement de résoudre l'équation différentelle
$\triangle{u}+\lambda{u}=0$ sur le disque. Calculer les zéros des
fonctions de Bessel étant beaucoup plus rapide et stable que de
résoudre l'équation différentielle, l'algorithme est d'autant plus
rapide et efficace. Également plusieurs des méthodes numériques pour
résoudre l'équation différentielle requiert l'utilisation de
matrices donc le nombre d'entrées est proportionnel au carré du
nombre de valeurs propres désirées et, par conséquent, il serait
inutile d'essayer de calculer par exemple les $1\cdot{10}^{6}$
premières valeurs propres en résolvant l'équation différentielle
avec ces méthodes numériques, ce qui par contre a été fait
raisonnablement avec l'algorithme.
\end{Rem}

\subsection{Résultats comparatifs} \label{Pics}

\noindent Le spectre du disque pour les valeurs propres
$\{\lambda_{j}\}_{j=1}^{10^{6}}$ a été calculé environ en
114576.4710 secondes de calcul soit environ 31 heures 50 minutes sur
un simple PC à deux processeurs avec 2GB de RAM. La liste d'une
taille de 30 \textbf{MB} est disponible auprès de l'auteur. Dans ce
ficher texte, la 1\iere{} colonne indique le rang de la
fonction/valeur propre, la 2\ieme{} colonne indique la valeur
propre, la 3\ieme{} indique l'ordre de la fonction de Bessel et la
4\ieme{} indique l'index du zéro de fonction Bessel.

\noindent Dans ce qui suit, pour chacune des valeurs propres listées
dans les tableaux, la représentation en courbes de niveaux de la
fonction propre correspondante est montrée plus loin sur les
graphiques. Grâce aux courbes de niveaux, il est facile d'y compter
le nombre de domaine nodaux et de \textbf{comparer} avec les valeurs
de $k$ (index) et $n$ (ordre) dans les tableaux. Les courbes de
niveaux des fonctions propres ont été obtenues à l'aide de Matlab en
solutionnant les équations différentielles à l'aide de la méthode
des éléments finis. Par exemple, à des fins comparatives, les 100
premières valeurs propres et leurs fonctions propres associées ont
été obtenues en près de 3 jours de calcul en Matlab. En plus de
prendre un temps excessivement long, la visulation des courbes à
niveau devient très difficile lorsque le rang des fonctions propres
est élevé.

\noindent Des programmes en Matlab sont également donnés en annexe.
Ces programmes permettent de résoudre les problèmes de valeurs
propres de même que de visualiser les ensembles nodaux des fonctions
propres.

\begin{table}[hb]
\begin{center}
\begin{tabular}{|c|c|c|c|c|}
\hline m $=$ rang & $\lambda_{m}$ (Matlab) & $j_{k}^{2}(n)$ (Maple) & n $=$ ordre & k $=$ index\\
\hline 3 & 149.4956 & 149.4529 & 2 & 1\\
\hline 6 & 278.9782 & 278.8316 & 3 & 1\\
\hline 7 & 310.5226 & 310.3223 & 1 & 4\\
\hline 16 & 646.8263 & 646.0310 & 5 & 1\\
\hline 25 & 991.8437 & 989.7291 & 3 & 5\\
\hline 27 & 1088.4330 & 1085.9440 & 4 & 4\\
\hline 30 & 1157.8175 & 1155.2319 & 5 & 3\\
\hline
\end{tabular}
\vspace{5mm}\caption{Résultats pour le secteur $S(\alpha)$ avec
$\alpha = \frac{\pi}{4}$ de rayon 1} \vspace{5mm}
\begin{tabular}{|c|c|c|c|c|}
\hline m $=$ rang & $\lambda_{m}$ (Matlab) & $j_{k}^{2}(n)$ (Maple) & n $=$ ordre & k $=$ index\\
\hline 3 & 23.0006 & 22.9968 & 3 & 1\\
\hline 9 & 58.4266 & 58.4019 & 7 & 1\\
\hline 11 & 69.3865 & 69.3521 & 8 & 1\\
\hline 19 & 108.3025 & 108.2183 & 2 & 3\\
\hline 26 & 152.1173 & 151.9596 & 14 & 1\\
\hline 29 & 165.6505 & 165.4521 & 5 & 3\\
\hline
\end{tabular}
\vspace{5mm}\caption{Résultats pour le secteur $S(\alpha)$ avec
$\alpha = 2\exp(1) \approx 5.4366$ de rayon 1}\vspace{-3mm}
\noindent Les figures 2.5 et 2.6 représentent les domaines nodaux
des fonctions propres correspondantes aux valeurs propres du
tableau, le nombre de composantes connexes est donné par $nk$.
\end{center}
\end{table}

\newpage
\begin{figure}[ht]
\begin{center}
\includegraphics[width = 7.0cm, height = 7.0cm]{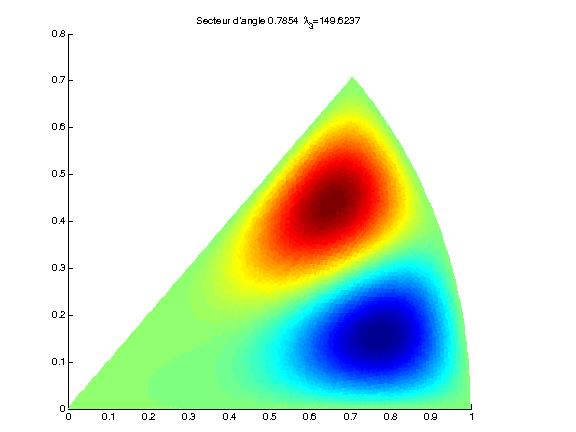}
\includegraphics[width = 7.0cm, height = 7.0cm]{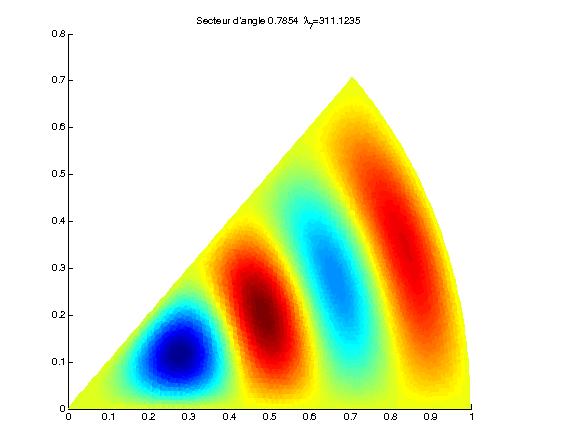}
\includegraphics[width = 7.0cm, height = 7.0cm]{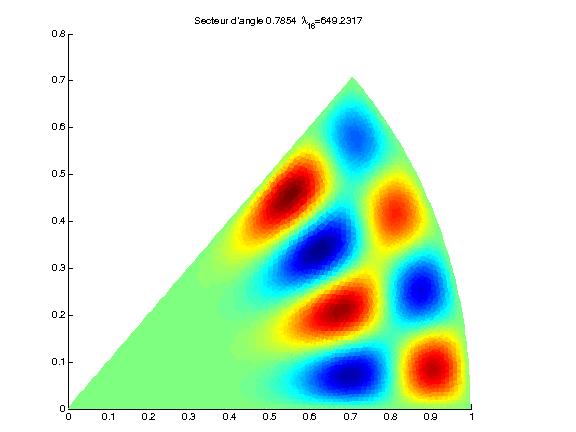}
\includegraphics[width = 7.0cm, height = 7.0cm]{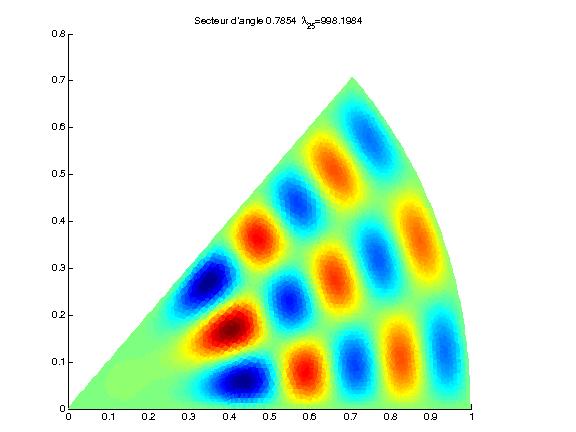}
\includegraphics[width = 7.0cm, height = 7.0cm]{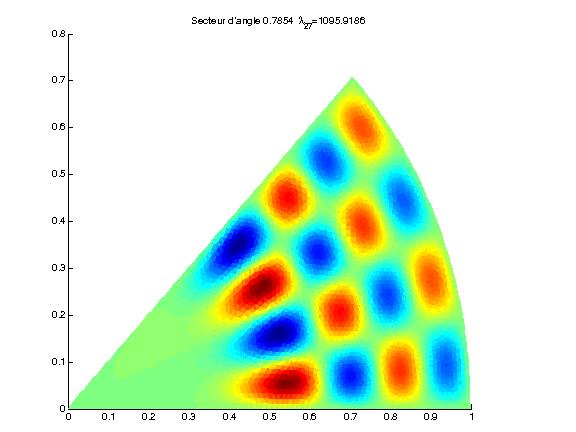}
\includegraphics[width = 7.0cm, height = 7.0cm]{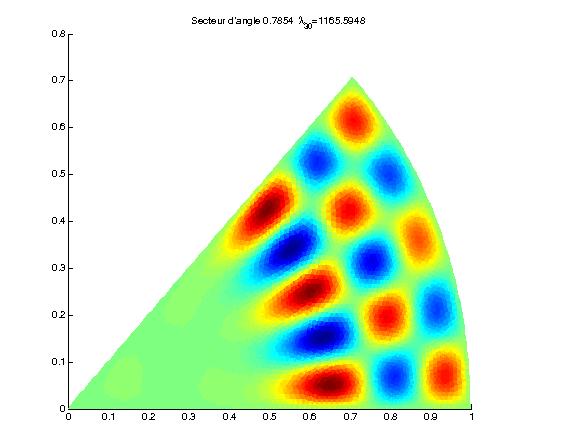}
\caption{Domaines nodaux de la 3\ieme{}, 7\ieme{}, 16\ieme{},
25\ieme{}, 27\ieme{} et 30\ieme{} fonction propre pour le secteur
$S(\alpha)$ avec $\alpha = \frac{\pi}{4}$ de rayon 1}
\end{center}
\end{figure}


\newpage
\begin{figure}[ht]
\begin{center}
\includegraphics[width = 7.0cm, height = 7.0cm]{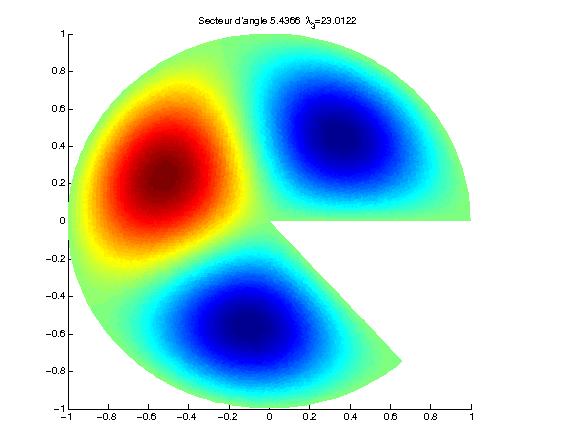}
\includegraphics[width = 7.0cm, height = 7.0cm]{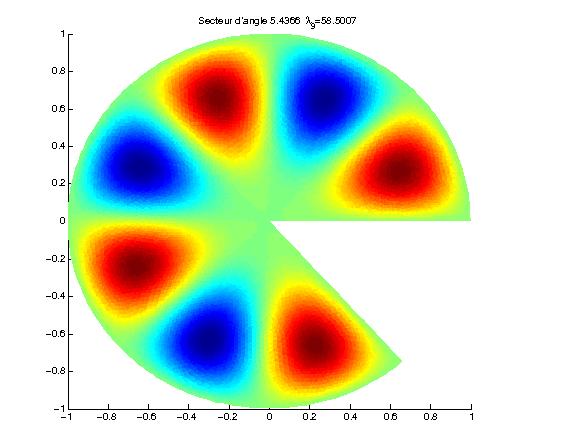}
\includegraphics[width = 7.0cm, height = 7.0cm]{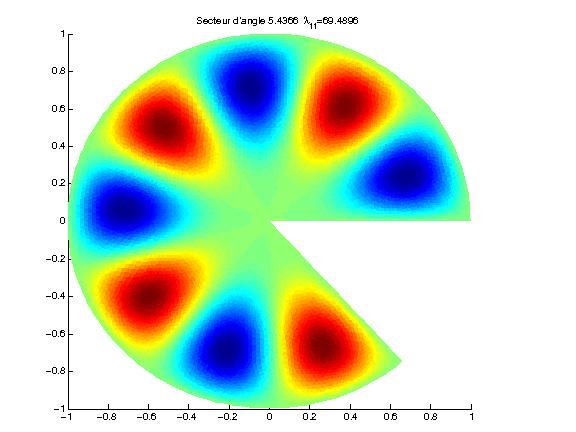}
\includegraphics[width = 7.0cm, height = 7.0cm]{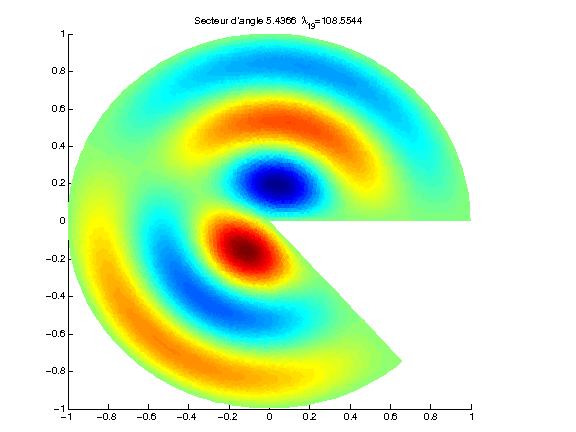}
\includegraphics[width = 7.0cm, height = 7.0cm]{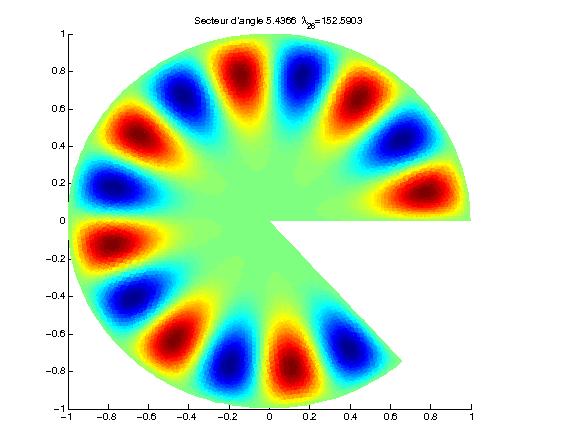}
\includegraphics[width = 7.0cm, height = 7.0cm]{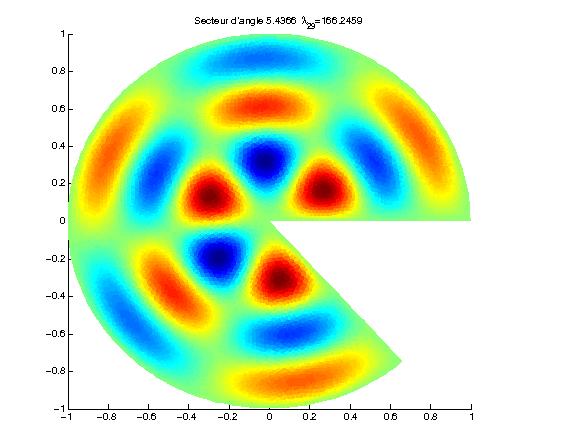}
\caption{Domaines nodaux de la 5\ieme{}, 9\ieme{}, 11\ieme{},
19\ieme{}, 26\ieme{} et 29\ieme{} fonction propre pour le secteur
$S(\alpha)$ avec $\alpha \approx 5.4366$ de rayon 1}
\end{center}
\end{figure}

\newpage
\begin{table}[ht]
\begin{center}
\begin{tabular}{|c|c|c|c|}
\hline m $=$ rang & $j_{k}^{2}(n)$ (Maple) & n $=$ ordre & k $=$ index \\
\hline 1 & 5.783186 & 0 & 1\\
\hline 4 & 26.374616 & 2 & 1\\
\hline 5 & 26.374616 & 2 & 1\\
\hline 6 & 30.471262 & 0 & 2\\
\hline 7 & 40.706466 & 3 & 1\\
\hline 8 & 40.706466 & 3 & 1\\
\hline 28 & 135.020709 & 2 & 3\\
\hline 29 & 135.020709 & 2 & 3\\
\hline 30 & 139.040284 & 0 & 4\\
\hline 31 & 149.452881 & 8 & 1\\
\hline 32 & 149.452881 & 8 & 1\\
\hline 33 & 152.241154 & 5 & 2\\
\hline 34 & 152.241154 & 5 & 2\\
\hline \ldots & & &\\
\hline 1000000 & 4004017.840283 & 1533 & 69\\
\hline 2000000 & 8005695.299643 & 391 & 714\\
\hline 3000000 & 12006894.927313 & 305 & 955\\
\hline 4000000 & 16008031.629743 & 1820 & 498\\
\hline 5000000 & 20008771.752617 & 1624 & 707\\
\hline 6000000 & 24009692.370612 & 713 & 1220\\
\hline 7000000 & 28010694.266301 & 3663 & 276\\
\hline 7912680 & 31685063.008767 & 0 & 1792\\
\hline
\end{tabular}
\vspace{5mm}\caption{Résultats pour le disque}\vspace{-3mm}
\noindent Les figures 2.7 et 2.8 représentent les domaines nodaux
des fonctions propres correspondantes aux valeurs propres du
tableau, le nombre de composantes connexes est donné par $2nk$ si
$n\geq 1$ (multiplicité 2) et par $k$ si $n=0$ (multiplicité 1).
\end{center}
\end{table}

\newpage

\begin{figure}[ht]
\begin{center}

\includegraphics[width = 7.0cm, height = 7.0cm]{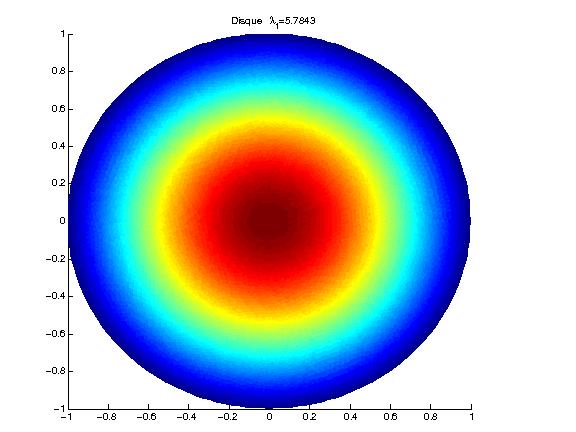}
\includegraphics[width = 7.0cm, height = 7.0cm]{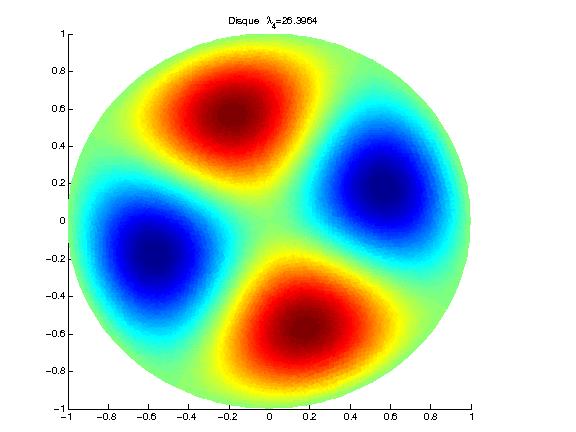}
\includegraphics[width = 7.0cm, height = 7.0cm]{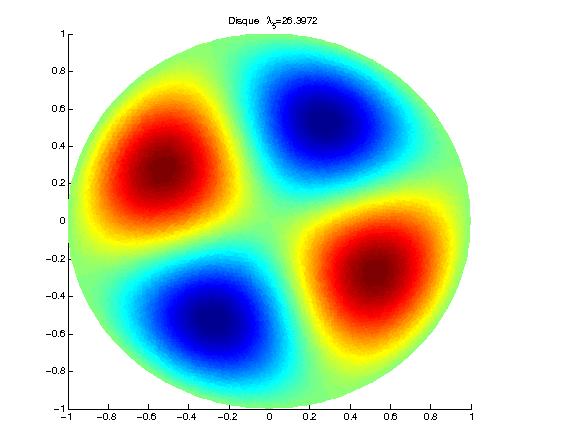}
\includegraphics[width = 7.0cm, height = 7.0cm]{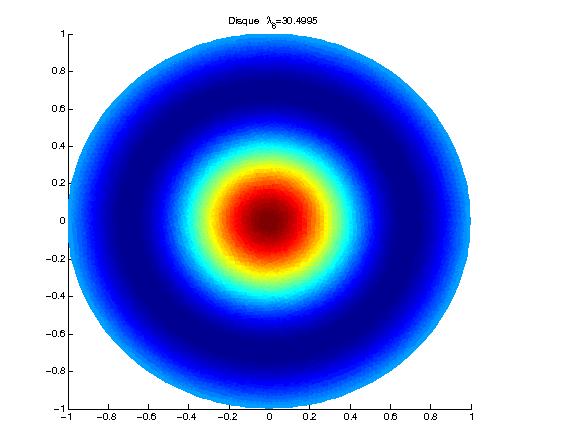}
\includegraphics[width = 7.0cm, height = 7.0cm]{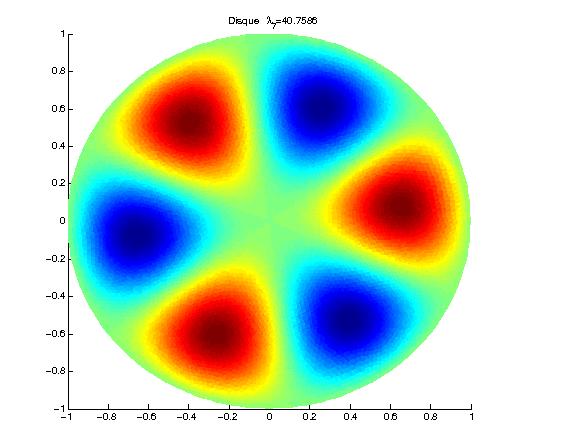}
\includegraphics[width = 7.0cm, height = 7.0cm]{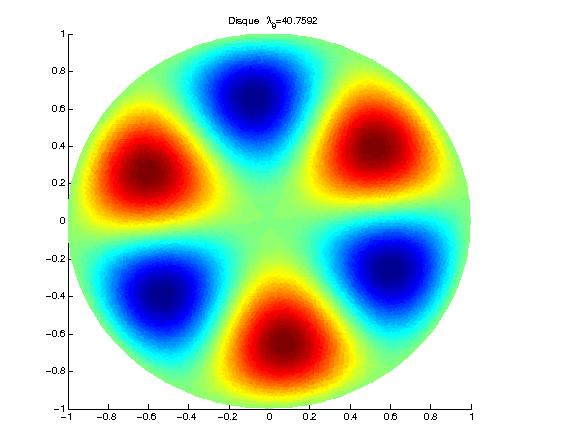}
\caption{Domaines nodaux de la 1\iere{}, 4\ieme{}, 5\ieme{},
6\ieme{}, 7\ieme{} et 8\ieme{} fonction propre pour le disque $D$ de
rayon 1}
\end{center}
\end{figure}

\newpage
\begin{figure}[ht]
\begin{center}

\includegraphics[width = 7.0cm, height = 7.0cm]{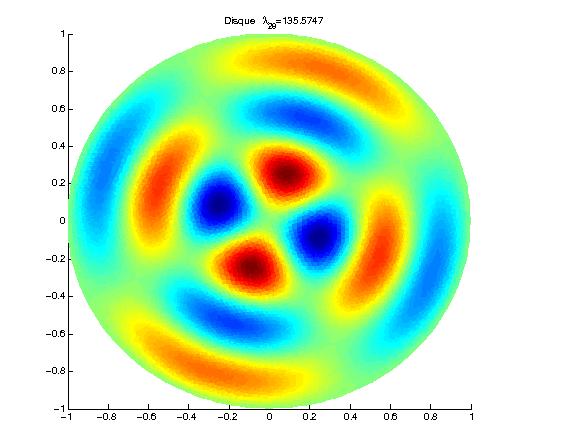}
\includegraphics[width = 7.0cm, height = 7.0cm]{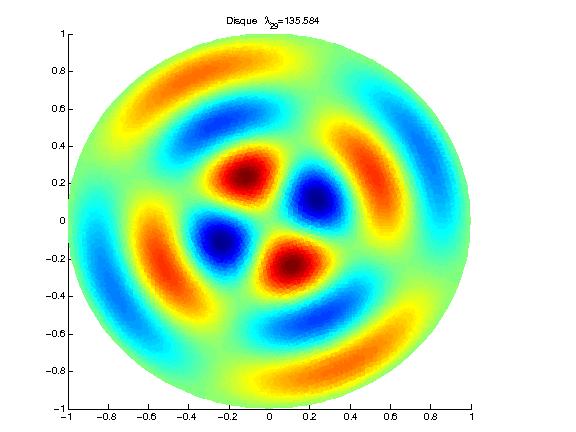}
\includegraphics[width = 7.0cm, height = 7.0cm]{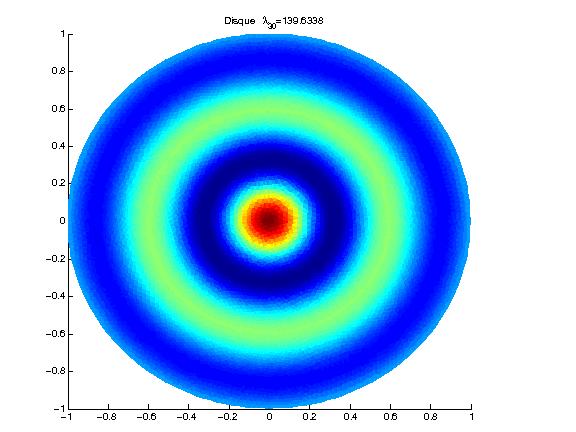}
\includegraphics[width = 7.0cm, height = 7.0cm]{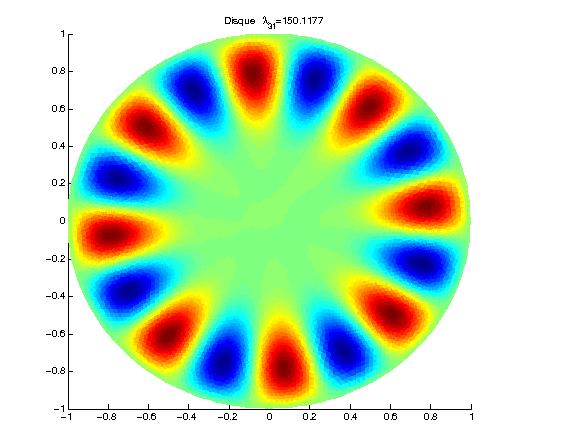}
\includegraphics[width = 7.0cm, height = 7.0cm]{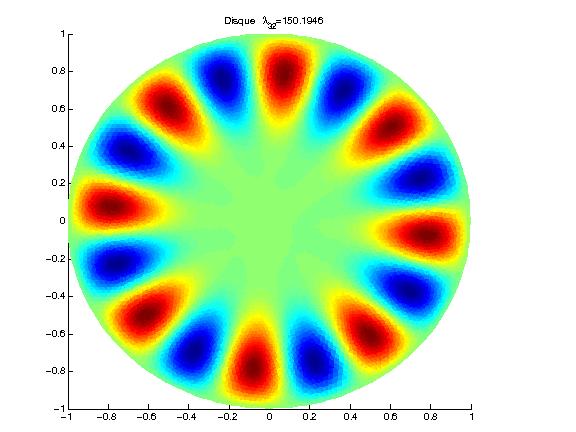}
\includegraphics[width = 7.0cm, height = 7.0cm]{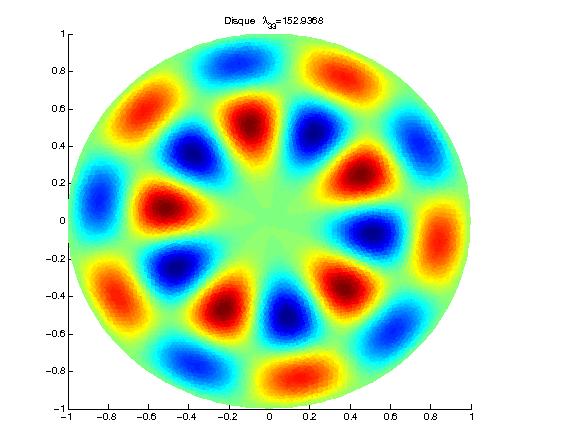}
\caption{Domaines nodaux de la 28\ieme{}, 29\ieme{}, 30\ieme{},
31\ieme{}, 32\ieme{} et 33\ieme{} fonction propre pour le disque $D$
de rayon 1}
\end{center}
\end{figure}

\newpage

\chapter[]{Loi de Weyl pour le disque}

\noindent Dans cette section, une explication approfondie de
l'article ~\cite{KuFe_1965} est donnée. C'est dans cet article qu'il
a été montré comment la loi de Weyl pour le disque $D$ de rayon r
est dans l'ordre $2/3$. Je tiens à remercier Igor Wigman et
Dominique Rabet pour les discussions très enrichissantes que nous
avons tenues afin d'éclaicir les multiples points très obscures de
l'article.

\noindent L'article montre que la loi de Weyl, en débutant avec le
problème $\triangle{u}+k^{2}u=0$, pour le disque $D$ de rayon r est
dans l'ordre $2/3$. En d'autres termes, soit le problème à valeur
propre $-\triangle{u}=k^{2}u$ avec $u|_{\partial{D}}=0$. Soit la
fonction de compte, $N_{D}(k)$, des valeurs propres
$\{k_{n}\}_{n=1}^{\infty}$ ordonnées en ordre croissant c'est-à-dire
soit
\begin{displaymath}
N_{D}(k)=\sum_{k_{n}<k}{1}
\end{displaymath}
alors il est montré que
\begin{displaymath}
N_{D}(k)=\frac{S}{4\pi}k^{2}-\frac{L}{4\pi}k+O(k^{\frac{2}{3}})
\end{displaymath}
où $S=\pi{r^{2}}$ et où $L=2\pi{r}$.

\noindent \textbf{Afin d'alléger la notation, nous écrirons
simplement $N(\lambda)$ au lieu de $N_{D}(\lambda)$. Seulement pour
l'explication du théorème de Kuznetsov et Fedosov, nous utiliserons
$k$ au lieu de $\lambda$ afin de conserver la notation identique à
celle de l'article facilitant ainsi les références.}

\noindent Avant de commencer à expliquer l'article, je tiens à faire
un \og détour \fg{} théorique sur l'approximation uniforme de la
fonction de Bessel de premier type, $J_{n}(x)$ pour $x>n>0$. Les
fonctions de Bessel de premier type d'ordre $n=0,1,2,\ldots$ sont au
coeur de l'étude du spectre du Laplacien sur le disque. Ce \og
détour \fg{} se justifie d'une part parce que les auteurs donne une
formule asymptotique pour $J_{n}(x)$ pour laquelle il faut
travailler quelque peu avant d'y arriver et, d'autre part, parce que
la formule donnée n'est pas effective tandis que je montre comment
obtenir une expansion uniforme asymptotique effective c'est-à-dire
sans terme O. La théorie sur les expansions uniformes des fonctions
a été principalement développée par F.W.J. Olver dans les années
1950.

\noindent À la fin de ce chapitre, je donne un algorithme permettant
de calculer efficacement et exactement $N(\lambda)$. En utilisant la
monotonicité des zéros des fonctions de Bessel ainsi qu'un principe
de \og marche \fg{} qui consiste à faire des retours en arrière et
des montées en alternance sur des paires d'entiers bien déterminés
par le problème, nous verrons qu'il est possible d'évaluer
$N(\lambda)$ efficacement sans avoir à calculer naïvement toutes les
valeurs propres telles que $\lambda_{j}\leq\lambda$. Par méthode
naïve, je veux signifier l'utilisation de l'algorithme de la section
2.2 afin de déterminer la valeur du rang maximale de la valeur
propre correspondant à la valeur de $N(\lambda)$.

\section{Approximation uniforme de la fonction de Bessel de premier type}

\noindent Le lemme suivant découle de deux théorèmes qui seront
exposés bientôt. La partie la plus difficile a été d'évaluer les
variations totales de certaines fonctions qui apparaîtront bientôt.
Il est à noter que les auteurs de \cite{LaWo_1995} affirment avoir
analytiquement trouvé les points stationnaires des fonctions en
question que vous pouvez consulter à la section 3.1.1 sans le
montrer dans leur article, ils ne se contentent que de donner les
valeurs numériques de variations totales qui sont évidemment les
mêmes que les miennes sans même donner explicitement les fonctions
comme je l'ai fait. Certaines des fonctions impliquées permettant de
borner l'approximation d'Olver avaient déjà été obtenues par Olver,
je me suis contenté de les recalculer numériquement pour obtenir
plus de décimales qu'Olver a pu le faire avec les ordinateurs des
années 50. Les résultats utiles pour cette section se retrouvent
dans \cite{Ol_1954_1}, \cite{Ol_1954_2} et \cite{Ol_1962}.

\begin{Lem}\label{lemcaca}
Soit $x>n>0$, $J_{n}(x)$ la fonction de Bessel de premier type et
\begin{eqnarray}
f_{n}(x)&=&(1+\delta_{3})\sqrt{\frac{\pi}{2}}(x^{2}-n^{2})^{\frac{1}{4}}J_{n}(x)\label{lemeqfn}
\end{eqnarray}
alors l'inégalité suivante est satisfaite
\begin{displaymath}
\Big|f_{n}(x)-\cos\Big(\eta-\frac{\pi}{4}\Big)\Big|\leq\frac{u_{1}}{\eta}+\frac{s_{2}}{n}\phantom{11111111111111111111111111111111111}
\end{displaymath}
\begin{displaymath}
+\bigg(|R_{2}|+\frac{s_{1}}{n^{2}}+\frac{u_{1}s_{1}}{n^{2}\eta}+\frac{|R_{2}|s_{1}}{n^{2}}+\frac{|v_{1}|s_{2}}{n\eta}+\frac{|R'_{2}|s_{2}}{n}+\sqrt{\pi}|\epsilon_{3}|\bigg)
\end{displaymath}
où $s_{1},s_{2},u_{1},v_{1},R_{2},R'_{2},\delta_{3},\epsilon_{3}$
sont exposés plus loin et où $\eta$ est comme suit
\begin{displaymath}
\eta = \sqrt{x^{2}-n^{2}}-n\arccos\Big(\frac{n}{x}\Big)
\end{displaymath}
De même, pour une valeur fixe arbitraire de $x$, il existe $n_{x}<x$
telle que
\begin{displaymath}
f_{n}(x)=\cos\Big(\eta-\frac{\pi}{4}\Big)+O\Big(\frac{1}{n}\Big)+O\Big(\frac{1}{\eta}\Big).
\end{displaymath}
\end{Lem}

\begin{Rem}[$n_{x}$]\label{remIgor}
La valeur $n_{x}$ peut être bien approximé en utilisant l'excellente
approximation de $\tilde{\eta}(y)=\sqrt{1-y^{2}}-y\arccos(y)$ pour
$y\in(0,1)$ qu'est
\begin{displaymath}
1<\frac{(1-y)^{\frac{3}{2}}}{\tilde{\eta}(y)}<\frac{3}{4}\sqrt{2}.
\end{displaymath}
\end{Rem}

\noindent En effet, de par l'inégalité de la remarque précédente,
qui m'a été fournie par Igor Wigman, avec $0<y=\frac{n}{x}<1$,
\begin{displaymath}
(1-y)^{\frac{3}{2}}<\frac{3\sqrt{2}}{4}\tilde{\eta}(y)\Rightarrow\frac{4x}{3\sqrt{2}}\Big(1-\frac{n}{x}\Big)^{\frac{3}{2}}<x\tilde{\eta}\Big(\frac{n}{x}\Big)=\eta(n,x)
\end{displaymath}
et, par conséquent,
\begin{displaymath}
1<\frac{4x}{3\sqrt{2}}\Big(1-\frac{n}{x}\Big)^{\frac{3}{2}}<\eta(x,n)\Leftrightarrow\Big(\frac{3\sqrt{2}}{4x}\Big)^{\frac{2}{3}}<1-\frac{n}{x}\Leftrightarrow
n<x-\Big(\frac{3\sqrt{2}}{4}\Big)^{\frac{2}{3}}x^{\frac{1}{3}}
\end{displaymath}

\noindent Il est à noter que
$x-x^{\frac{4}{9}}<x-\big(\frac{3\sqrt{2}}{4}\big)^{\frac{2}{3}}x^{\frac{1}{3}}$
pour $x>\big(\frac{3\sqrt{2}}{4}\big)^{6}\approx1.423828$.

\noindent De façon équivalente \textbf{comme il en ressortira du
théorème \ref{pourlaconne}},
\begin{eqnarray}
f_{n}(x)&=&\theta_{1}(\eta)+\theta_{1}(\eta)\frac{A_{1}(\zeta)}{n^{2}}+\frac{1}{n}\theta_{2}(\eta)(-\zeta)^{\frac{1}{2}}B_{0}(\zeta)+\sqrt{\pi}\epsilon_{3}\label{eqnfn}
\end{eqnarray}
où
\begin{eqnarray*}
\theta_{1}(\eta)&=&\cos\Big(\eta-\frac{1}{4}\pi\Big)+\sin\Big(\eta-\frac{1}{4}\pi\Big)\frac{u_{1}}{\eta^{1}}+R_{2}\\
\theta_{2}(\eta)&=&\sin\Big(\eta-\frac{1}{4}\pi\Big)-\cos\Big(\eta-\frac{1}{4}\pi\Big)\frac{v_{1}}{\eta^{1}}+R'_{2}
\end{eqnarray*}
alors
\begin{displaymath}
J_{n}(x)=\frac{1}{1+\delta_{3}}\sqrt{\frac{2}{\pi}}(x^{2}-n^{2})^{-\frac{1}{4}}f_{n}(x)
\end{displaymath}

\begin{Rem}[Dénombrement des zéros de $J_{n}(x)$]
Seule la connaissance de $f_{n}(x)$ est nécessaire au dénombrement
des zéros de $J_{n}(x)$.
\end{Rem}

\begin{center}\line(1,0){300}\end{center}

\noindent Avant de débuter la preuve, voici les deux théorèmes qui
seront utilisés pour prouver le lemme et se retrouvant
respectivement dans ~\cite{Ol_1964} et ~\cite{FLO_2004}.

\noindent Avant d'énoncer les théorèmes, voici les définitions des
fonctions qui apparaîtront dans les théorèmes. Dans ce qui suit,
$Ai(x)$ et $Bi(x)$ dénotent respectivement la fonction d'Airy de
premier type et deuxième type. Ainsi,
\begin{eqnarray*}
E(x)&=&\bigg(\frac{Bi(x)}{Ai(x)}\bigg)^{\frac{1}{2}}\phantom{123}\textrm{pour $c\leq x \leq\infty$}\\
E(x)&=&1\phantom{123}\textrm{pour $-\infty\leq x \leq c$}\\
M(x)&=&(E^{2}(x)Ai^{2}(x)+E^{-2}(x)Bi^{2}(x))^{\frac{1}{2}}\phantom{123}\textrm{pour
$x\in\mathbb{R}$}\\
\lambda&=&\sup_{x\in\mathbb{R}}(\pi|x|^{\frac{1}{2}}M^{2}(x))\\
&\approx&1.039522542988\phantom{12}\textrm{à}\phantom{12}x\approx1.321915092767\\
\end{eqnarray*}
où $c$ est la plus grande valeur négative telle que $Ai(x)=Bi(x)$.
\textbf{Nous nous entendons bien que $\lambda$ n'a aucun lien avec
les valeurs propres. C'est la lettre utilisée par Olver pour
désigner le suprémum ci-dessus.}

\noindent Quant au symbole $\mathcal{V}_{(a,b)}(g)$, il dénote la
variation totale d'une fonction $g$ sur un intervalle
$(a,b)\subseteq\mathrm{dom}(g)$. Pour plus d'information, voir le
1\ier{} chapitre de ~\cite{Ol_1974}.

\begin{Thm}[Approximation uniforme de
$J_{n}(nz)$]\label{pourlaconne}
Soit $J_{n}(nz)$ la fonction de Bessel de premier type d'ordre
$n\in\mathbb{R}$ et d'argument $nz$ avec $z\in\mathbb{C}$, alors
l'expansion uniforme de $J_{n}(nz)$ est donné par
\begin{eqnarray}\label{eqbesapp}
J_{n}(nz)&=&\frac{1}{1+\delta_{2l+1}}\frac{1}{n^{\frac{1}{3}}}\bigg(\frac{4\zeta}{1-z^{2}}\bigg)^{\frac{1}{4}}\bigg(Ai(n^{\frac{2}{3}}\zeta)\sum_{s=0}^{l}{\frac{A_{s}(\zeta)}{n^{2s}}}\\
&&+\frac{Ai'(n^{\frac{2}{3}}\zeta)}{n^{\frac{4}{3}}}\sum_{s=0}^{l-1}{\frac{B_{s}(\zeta)}{n^{2s}}}+\epsilon_{2l+1}\bigg)
\end{eqnarray}
où les termes d'erreur $\delta_{2l+1}$ et $\epsilon_{2l+1}$ sont
donnés comme suit
\begin{eqnarray}
|\delta_{2l+1}|&\leq&\frac{2}{n^{2l+1}}\exp\Big(2\lambda\frac{1}{n}\mathcal{V}_{-\infty,\infty}(|\zeta|^{\frac{1}{2}}B_{0}(\zeta))\Big)\mathcal{V}_{-\infty,\infty}(|\zeta|^{\frac{1}{2}}B_{l}(\zeta))\label{boulette1}\\
|\epsilon_{2l+1}|&\leq&\frac{2}{n^{2l+1}}\exp\Big(2\lambda\frac{1}{n}\mathcal{V}_{\zeta,\infty}(|\zeta|^{\frac{1}{2}}B_{0}(\zeta))\Big)\mathcal{V}_{\zeta,\infty}(|\zeta|^{\frac{1}{2}}B_{l}(\zeta))\frac{M(n^{\frac{2}{3}}\zeta)}{E(n^{\frac{2}{3}}\zeta)}\label{boulette2}
\end{eqnarray}
où
$\mathcal{V}_{\zeta,\infty}(|\zeta|^{\frac{1}{2}}B_{l}(\zeta))\leq\mathcal{V}_{-\infty,\infty}(|\zeta|^{\frac{1}{2}}B_{l}(\zeta))$.
\end{Thm}

\noindent \textbf{Le lecteur intéressé à connaître les récurrences
permettant de calculer $A_{s}$ et $B_{s}$ dans le théorème qui suit
peuvent consulter ~\cite{Ol_1954_2}. J'ai trouvé $A_{1}(\zeta)$ et
$B_{1}(\zeta)$ en résolvant les récurrences, mais cela ne requiert
qu'une bonne dose de patience et je me suis contenté de les donner à
la section 3.1.1. $\zeta$ est donné plus loin à l'équation
(\ref{eqnzeta}).}

\begin{Rem}[$n$ et $z$]
Dans ce qui suit, seuls les cas où $n\in\mathbb{N}$ et
$z\in\mathbb{R}^{+}$ seront utiles.
\end{Rem}

\begin{Thm}[Approximations uniformes de $Ai$ et
$Ai'$.]\label{thmappaibi} Soit $Ai$ la fonction de Airy de premier
type sa dérivé $Ai'$. Soit $x>0$ et
$\frac{2}{3}x^{\frac{3}{2}}=\xi$. Alors
\begin{eqnarray}\label{eqaiapp}
Ai(-x)&=&\frac{1}{\sqrt{\pi}x^{\frac{1}{4}}}\bigg(\cos\Big(\xi-\frac{1}{4}\pi\Big)\sum_{s=0}^{\lfloor\frac{1}{2}S-\frac{1}{2}\rfloor}{(-1)^{s}\frac{u_{2s}}{\xi^{2s}}}\\
&&+\sin\Big(\xi-\frac{1}{4}\pi\Big)\sum_{s=0}^{\lfloor\frac{1}{2}S-1\rfloor}{(-1)^{s}\frac{u_{2s+1}}{\xi^{2s+1}}}+R_{S}\bigg)
\end{eqnarray}
où
\begin{displaymath}
|R_{S}|\leq 2\exp\Big(\frac{5}{36\xi}\Big)\frac{u_{S}}{\xi^{S}}
\end{displaymath}
et
\begin{eqnarray}\label{eqaipapp}
Ai'(-x)&=&\frac{x^{\frac{1}{4}}}{\sqrt{\pi}}\bigg(\sin\Big(\xi-\frac{1}{4}\pi\Big)\sum_{s=0}^{\lfloor\frac{1}{2}S-\frac{1}{2}\rfloor}{(-1)^{s}\frac{v_{2s}}{\xi^{2s}}}\\
&&-\cos\Big(\xi-\frac{1}{4}\pi\Big)\sum_{s=0}^{\lfloor\frac{1}{2}S-1\rfloor}{(-1)^{s}\frac{v_{2s+1}}{\xi^{2s+1}}}+R'_{S}\bigg)
\end{eqnarray}
où
\begin{displaymath}
|R'_{S}|\leq 2\exp\Big(\frac{7}{36\xi}\Big)\frac{|v_{S}|}{\xi^{S}}
\end{displaymath}
avec
\begin{eqnarray}
u_{s}&=&\frac{2^{s}\Gamma(3s+\frac{1}{2})}{\sqrt{\pi}3^{3s}\Gamma(2s+1)}\phantom{123}s>0\label{equs}\\
v_{s}&=&-\frac{(6s+1)}{(6s-1)}u_{s}\phantom{123}s>0\label{eqvs}\\
u_{0}&=&1\label{equ0}
\end{eqnarray}
\end{Thm}

\noindent Voici maintenant la preuve du lemme \ref{lemcaca}
\begin{Prflemme}
\noindent À partir des équations (\ref{eqbesapp}), (\ref{boulette1})
et (\ref{boulette2}) du théorème \ref{pourlaconne} avec $l=1$ dans
l'approximation de $J_{n}(nz)$, nous avons que
\begin{displaymath}
J_{n}(nz)=\frac{1}{1+\delta_{3}}\frac{1}{n^{\frac{1}{3}}}\bigg(\frac{4\zeta}{1-z^{2}}\bigg)^{\frac{1}{4}}\bigg(Ai(n^{\frac{2}{3}}\zeta)\Big(1+\frac{A_{1}(\zeta)}{n^{2}}\Big)+\frac{Ai'(n^{\frac{2}{3}}\zeta)}{n^{\frac{4}{3}}}B_{0}(\zeta)+\epsilon_{3}\bigg)\\
\end{displaymath}
avec
\begin{eqnarray}
|\delta_{3}|&\leq&\frac{0.013227}{n^{3}}\exp\Big(\frac{0.218423}{n}\Big)\label{eqd3}\\
|\epsilon_{3}|&\leq&\frac{0.008352}{n^{3}}\exp\Big(\frac{0.218423}{n}\Big)\label{eqe3}.
\end{eqnarray}
\noindent Nous reviendrons à la section \ref{pipicaca} sur les
évaluations numériques de $\delta_{3}$ et $\epsilon_{3}$.

\noindent Maintenant, un changement de variable qui \og défait \fg{}
l'homothétie $nz$ est donné par $x=nz$. En se concentrant seulement
sur $z>1$ c'est-à-dire sur $z\in(1,\infty)$, alors (voir
~\cite{Ol_1962})
\begin{eqnarray}\label{eqnzeta}
\zeta&=&-\bigg(\frac{3}{2}\int_{1}^{z}{\frac{\sqrt{z^{2}-1}}{z}dz}\bigg)^{\frac{2}{3}}\\
&=&-\bigg(\frac{3}{2}\bigg)^{\frac{2}{3}}\eta^{\frac{2}{3}}\frac{1}{n^{\frac{2}{3}}}\phantom{123}x>n>0
\end{eqnarray}
où
\begin{displaymath}
\eta
=\sqrt{x^{2}-n^{2}}-n\arccos\Big(\frac{n}{x}\Big)\phantom{123}x>n>0
\end{displaymath}

\noindent Ainsi
\begin{eqnarray*}
Ai(n^{\frac{2}{3}}\zeta)&=&Ai(-(n^{\frac{2}{3}}(-\zeta)))\\
Ai'(n^{\frac{2}{3}}\zeta)&=&Ai'(-(n^{\frac{2}{3}}(-\zeta)))\\
\xi&=&\eta\\
\end{eqnarray*}

\noindent En substituant $S=2$ dans les approximations de $Ai$ et
$Ai'$ données respectivement par (\ref{eqaiapp}) et
(\ref{eqaipapp}), et en dénotant pour des raisons de commodité
$\theta_{1}(\eta)$ et $\theta_{2}(\eta)$ par
\begin{eqnarray}
\theta_{1}(\eta)&=&\cos\Big(\eta-\frac{1}{4}\pi\Big)+\sin\Big(\eta-\frac{1}{4}\pi\Big)\frac{u_{1}}{\eta^{1}}+R_{2}\label{eqthe1}\\
\theta_{2}(\eta)&=&\sin\Big(\eta-\frac{1}{4}\pi\Big)-\cos\Big(\eta-\frac{1}{4}\pi\Big)\frac{v_{1}}{\eta^{1}}+R'_{2}\label{eqthe2}\\
|R_{2}|&\leq&\exp\Big(\frac{5}{36\eta}\Big)\frac{45}{74}\frac{1}{\eta^{2}}\label{eqR2}\\
|R'_{2}|&\leq&\exp\Big(\frac{7}{36\eta}\Big)\frac{65}{74}\frac{1}{\eta^{2}}\label{eqRP2}
\end{eqnarray}
alors
\begin{eqnarray*}
Ai(n^{\frac{2}{3}}\zeta)&=&\frac{1}{\sqrt{\pi}}\frac{1}{n^{\frac{1}{6}}}\frac{1}{(-\zeta)^{\frac{1}{4}}}\theta_{1}(\eta)\\
Ai'(n^{\frac{2}{3}}\zeta)&=&\frac{1}{\sqrt{\pi}}n^{\frac{1}{6}}(-\zeta)^{\frac{1}{4}}\theta_{2}(\eta)
\end{eqnarray*}

\noindent Par conséquent, pour $z>1$ ou bien $\zeta<0$,
\begin{eqnarray*}
J_{n}(nz)&=&\frac{1}{1+\delta_{3}}\frac{1}{n^{\frac{1}{3}}}\bigg(\frac{4}{z^{2}-1}\bigg)^{\frac{1}{4}}\bigg(\frac{1}{\sqrt{\pi}}\frac{1}{n^{\frac{1}{6}}}\theta_{1}(\eta)\Big(1+\frac{A_{1}(\zeta)}{n^{2}}\Big)\phantom{1234567890}\\
&&\phantom{12345678901234567890}+\frac{1}{\sqrt{\pi}}\frac{n^{\frac{1}{6}}}{n^{\frac{4}{3}}}\theta_{2}(\eta)(-\zeta)^{\frac{1}{2}}B_{0}(\zeta)+\epsilon_{3}\bigg)
\end{eqnarray*}
et donc pour $x>n>0$ avec $\zeta<0$,
\begin{eqnarray*}
J_{n}(x)&=&\frac{1}{1+\delta_{3}}\sqrt{\frac{2}{\pi}}(x^{2}-n^{2})^{-\frac{1}{4}}\bigg(\theta_{1}(\eta)+\theta_{1}(\eta)\frac{A_{1}(\zeta)}{n^{2}}\phantom{1234567890}\\
&&\phantom{12345678901234567890}+\frac{1}{n}\theta_{2}(\eta)(-\zeta)^{\frac{1}{2}}B_{0}(\zeta)+\sqrt{\pi}\epsilon_{3}\bigg)
\end{eqnarray*}

\noindent Maintenant soit (voir ~\cite{Ol_1962})
\begin{eqnarray}
s_{1}&=&\sup_{\zeta<0}{|A_{1}(\zeta)|}=\frac{1}{225}\label{eqs1}\\
s_{2}&=&|\zeta(z_{1})|^{\frac{1}{2}}|B_{0}(z_{1})|\label{eqs2}\\
&\approx&0.010862854400
\end{eqnarray}
et où $z_{1}$ est le seul point stationnaire de
$|\zeta(z)|^{\frac{1}{2}}B_{0}(\zeta(z))$ pour
$z>1\Leftrightarrow\zeta<0$ tel que
\begin{displaymath}
z_{1}=1.979495483061\phantom{12}\textrm{avec}\phantom{12}\zeta(z_{1})=-1.000459796360\phantom{1}.
\end{displaymath}
Par conséquent,
\begin{eqnarray*}
|\theta_{1}(\eta)A_{1}(\zeta)|&\leq&|\theta_{1}(\eta)|s_{1}\\
&=&\frac{s_{1}}{n^{2}}+\frac{u_{1}s_{1}}{n^{2}\eta}+\frac{|R_{2}|s_{1}}{n^{2}}\\
\frac{1}{n}|\theta_{2}(\eta)(-\zeta)^{\frac{1}{2}}B_{0}(\zeta)|&\leq&\frac{1}{n}|\theta_{2}(\eta)|s_{2}\\
&\leq&\frac{s_{2}}{n}+\frac{|v_{1}|s_{2}}{n\eta}+\frac{|R'_{2}|s_{2}}{n}
\end{eqnarray*}

\noindent En dénotant par
\begin{eqnarray*}
f_{n}(x)&=&(1+\delta_{3})\sqrt{\frac{\pi}{2}}(x^{2}-n^{2})^{\frac{1}{4}}J_{n}(x)\\
&=&\theta_{1}(\eta)+\theta_{1}(\eta)\frac{A_{1}(\zeta)}{n^{2}}+\frac{1}{n}\theta_{2}(\eta)(-\zeta)^{\frac{1}{2}}B_{0}(\zeta)+\sqrt{\pi}\epsilon_{3}
\end{eqnarray*}
alors pour $x>n>0$, l'inégalité suivante est satisfaite
\begin{displaymath}
\Big|f_{n}(x)-\cos\Big(\eta-\frac{\pi}{4}\Big)\Big|\leq\frac{u_{1}}{\eta}+\frac{s_{2}}{n}\phantom{11111111111111111111111111111111111}
\end{displaymath}
\begin{displaymath}
+\bigg(|R_{2}|+\frac{s_{1}}{n^{2}}+\frac{u_{1}s_{1}}{n^{2}\eta}+\frac{|R_{2}|s_{1}}{n^{2}}+\frac{|v_{1}|s_{2}}{n\eta}+\frac{|R'_{2}|s_{2}}{n}+\sqrt{\pi}|\epsilon_{3}|\bigg)
\end{displaymath}

\noindent \textbf{Ceci termine la preuve de la première inégalité.
Pour la deuxième partie du lemme, procédons comme ce qui suit.}

\noindent Soit la fonction $f_{n}(x)$ donnée par l'équation
(\ref{eqnfn}), il est possible de regarder à $\eta(x,n)$ pour
$x>n>0$ de deux façons. La première façon est lorsque $n$ est fixe
et $x$ varie, lequel cas la fonction $g_{1}(x):=:\eta(x,n)$ est
croissante en $x>n$. La deuxième façon, ce qui sera utile, est
lorsque $x$ est fixe et $n$ varie, dans ce cas la fonction
$g_{2}(n):=:\eta(x,n)$ est décroissante pour tout $n\in(0,x)$. Ainsi
pour une valeur fixe arbitraire $x>0$, il existe $n_{x}<x$ telle que
pout tout $n$ avec $n\leq n_{x}$, la fonction $\eta^{-1}(x,n)<1$.
Par conséquent,
\begin{displaymath}
\frac{1}{\eta^{2}}\leq\frac{1}{\eta}\leq
1\phantom{12}\textrm{pour}\phantom{1}n\leq n_{x}
\end{displaymath}
De même, pour toute constante $k$,
\begin{displaymath}
\exp\Big(\frac{k}{\eta}\Big)\leq
\exp(k)\phantom{12}\textrm{pour}\phantom{1}n\leq n_{x}
\end{displaymath}

%

\noindent De par la remarque \ref{remIgor} et les équations
(\ref{equs}), (\ref{eqvs}), (\ref{eqs1}) et (\ref{eqs2}), alors
$u_{1}=\frac{5}{72}$, $v_{1}=-\frac{7}{72}$, $s_{1}=\frac{1}{225}$
et $s_{2}\approx0.010862854400$. Si, en plus, $n\leq n_{x}$,
et, en posant les constantes
\begin{eqnarray*}
A&\approx&s_{2}+s_{1}+u_{1}s_{1}+\exp\Big(\frac{5}{36}\Big)\frac{45}{74}s_{1}+|v_{1}|s_{2}+\exp\Big(\frac{7}{36}\Big)\frac{65}{74}s_{2}\phantom{1234567890}\\
&&\phantom{123456789012345678901234567890}+0.008352\sqrt{\pi}\exp(0.218423)\\
&\approx&0.049784723505\\
B&=&u_{1}+\exp\Big(\frac{5}{36}\Big)\frac{45}{74}\\
&\approx&0.768158487672
\end{eqnarray*}
alors
\begin{displaymath}
\Big|f_{n}(x)-\cos\Big(\eta-\frac{\pi}{4}\Big)\Big|\leq\frac{A}{n}+\frac{B}{\eta}
\end{displaymath}
et
\begin{displaymath}
f_{n}(x)=\cos\Big(\eta-\frac{\pi}{4}\Big)+O\Big(\frac{1}{n}\Big)+O\Big(\frac{1}{\eta}\Big)
\end{displaymath}
\begin{flushright}$\blacksquare$\end{flushright}
\end{Prflemme}

\subsection{Calculs numériques des variations
totales}\label{pipicaca}

\noindent Pour plus de détails concernant ce qui suit, se référer
principalement à ~\cite{Ol_1962}.

\noindent La valeur de
$\mathcal{V}_{\mathbb{R}}(|\zeta|^{\frac{1}{2}}B_{0}(\zeta))$ est
calculée en trouvant les points stationnaires de
$|\zeta|^{\frac{1}{2}}B_{0}(\zeta)$ et en utilisant la formule de la
variation totale. L'unique point stationnaire autre que $\zeta=0$ se
situe à $\zeta=\zeta(z_{1})\approx-1.000459796360$ où
$z_{1}\approx1.979495483061$. Par conséquent, en dénotant par
\begin{displaymath}
g(z)=|\zeta(z)|^{\frac{1}{2}}B_{0}(\zeta(z))
\end{displaymath}
alors
\begin{eqnarray*}
\mathcal{V}_{\mathbb{R}}(|\zeta|^{\frac{1}{2}}B_{0}(\zeta))&=&
\lim_{z\to 0}{g(z)} - 2g(z=1) + 2g(z=z_{1}) + \lim_{z\to \infty}{g(z)}\\
&\approx&\frac{1}{12} - 0 + 2(0.010862854400) + 0\\
&\approx&0.105059042134
\end{eqnarray*}
Également
\begin{eqnarray*}
\lim_{\zeta\to\infty}{|\zeta|^{\frac{1}{2}}B_{0}(\zeta)}&=&\lim_{\zeta\to\infty}\bigg(-\frac{3}{24}+\frac{5}{24}-\frac{5}{48\zeta^{\frac{3}{2}}}\bigg)\\
&=&\lim_{z\to0}{\bigg(-\frac{3}{24}+\frac{5}{24}-\frac{5}{48\zeta(z)^{\frac{3}{2}}}\bigg)}\\
&=&\bigg(-\frac{3}{24}+\frac{5}{24}-0\bigg)\\
&=&\frac{1}{12}\\
\lim_{\zeta\to-\infty}{|\zeta|^{\frac{1}{2}}B_{0}(\zeta)}&=&\lim_{z\to\infty}{g(z)}\\
&=&0\\
\end{eqnarray*}

\noindent Pour trouver la valeur numérique de
$\mathcal{V}_{\mathbb{R}}(|\zeta|^{\frac{1}{2}}B_{1}(\zeta))$, il
faut les points stationnaires de
$|\zeta|^{\frac{1}{2}}B_{1}(\zeta)$.

\noindent En dénotant par
\begin{displaymath}
g(z)=|\zeta(z)|^{\frac{1}{2}}B_{1}(\zeta(z))
\end{displaymath}
alors ses points stationnaires autres que $z=1$ sont
\begin{eqnarray*}
z_{1}&\approx&0.138560281581\\
z_{2}&\approx&1.418538099456\\
\end{eqnarray*}
\noindent Par conséquent, la variation totale est
\begin{eqnarray*}
\mathcal{V}_{\mathbb{R}}(|\zeta|^{\frac{1}{2}}B_{1}(\zeta))&=&
\lim_{z\to 0}{g(z)} - 2g(z_{1}) + 2g(z=1) - 2g(z=z_{2}) + \lim_{z\to\infty}{g(z)}\\
&\approx&-0.002681327160-2(-0.004008186698
)\\
&&\phantom{1}+2(0)-2(-0.000639161111)+0\\
&\approx&0.006613368457\\
\end{eqnarray*}
\noindent avec les valeurs des limites suivantes
\begin{eqnarray*}
\lim_{\zeta\to\infty}{|\zeta|^{\frac{1}{2}}B_{1}(\zeta)}&=&\lim_{z\to 0}{g(z)}\\
&=&-\frac{139}{51840}\\
\lim_{\zeta\to-\infty}{|\zeta|^{\frac{1}{2}}B_{1}(\zeta)}&=&\lim_{z\to\infty}{g(z)}\\
&=&0\\
\end{eqnarray*}

\noindent Maintenant, je donne les coefficients $A_{s}(\zeta)$ et
$B_{s}(\zeta)$ qui n'ont pas été faciles à trouver pour $s=1$.

\noindent Les coefficients $A_{s}(\zeta)$ et $B_{s}(\zeta)$ pour
$s=0,1$ suivent et, dans ce qui suit, $t=(1-z^{2})^{-\frac{1}{2}}$,
$t=-it_{1}$ avec $t_{1}=(z^{2}-1)^{-\frac{1}{2}}$ et
$\zeta_{1}=-\zeta$.
\begin{eqnarray*}
B_{0}(\zeta)&=&-\frac{3t-5t^{3}}{24\zeta^{\frac{1}{2}}}-\frac{5}{48\zeta^{2}}\phantom{12}\zeta>0\phantom{12}\textrm{ou}\phantom{12}0\leq
z\leq1\\
B_{0}(\zeta)&=&\frac{3t_{1}+5t_{1}^{3}}{24\zeta_{1}^{\frac{1}{2}}}-\frac{5}{48\zeta_{1}^{2}}\phantom{12}\zeta<0\phantom{12}\textrm{ou}\phantom{12}z\geq1\\
B_{1}(\zeta)&=&\frac{-(30375t^{3}-369603t^{5}+765765t^{7}-425425t^{9})}{414720\zeta^{\frac{1}{2}}}\\
&&\phantom{1}-\frac{(405t^{2}-2310t^{4}+1925t^{6})}{55296\zeta^{2}}\\
&&\phantom{1}-\frac{(1155t-1925t^{3})}{110592\zeta^{\frac{7}{2}}}-\frac{85085}{663552\zeta^{5}}\phantom{12}\zeta>0\phantom{12}\textrm{ou}\phantom{12}0\leq z\leq 1\\
B_{1}(\zeta)&=&\frac{(30375t_{1}^{3}+369603t_{1}^{5}+765765t_{1}^{7}+425425t_{1}^{9})}{414720\zeta_{1}^{\frac{1}{2}}}\\
&&\phantom{1}+\frac{(405t_{1}^{2}+2310t_{1}^{4}+1925t_{1}^{6})}{55296\zeta_{1}^{2}}\\
&&\phantom{1}+\frac{(1155t_{1}+1925t_{1}^{3})}{110592\zeta_{1}^{\frac{7}{2}}}+\frac{85085}{663552\zeta_{1}^{5}}\phantom{12}\zeta<0\phantom{12}\textrm{ou}\phantom{12}z\geq 1\\
\end{eqnarray*}
\begin{eqnarray*}
A_{1}(\zeta)&=&\frac{81t^{2}-462t^{4}+385t^{8}}{1152}-\frac{7(3t-5t^{3})}{1152\zeta^{\frac{3}{2}}}-\frac{455}{4608\zeta^{3}}\phantom{12}\zeta>0\phantom{12}\textrm{ou}\phantom{12}0\leq
z\leq1\\
A_{1}(\zeta)&=&-\frac{81t_{1}^{2}+462t_{1}^{4}+385t^{8}}{1152}-\frac{7(3t_{1}+5t_{1}^{3})}{1152\zeta_{1}^{\frac{3}{2}}}+\frac{455}{4608\zeta_{1}^{3}}\phantom{12}\zeta<0\phantom{12}\textrm{ou}\phantom{12}z\geq1\\
\end{eqnarray*}

\noindent Les valeurs de $A_{s}(0)$ et de $B_{s}(0)$ peuvent être
calculées en utilisant les développements de Meissel de $J_{n}(n)$
que nous pouvons retrouver dans ~\cite{Wa_1944}. Ces valeurs se
retrouvent également dans ~\cite{Ol_1962} et, par conséquent,
\begin{eqnarray*}
A_{1}(0)&=&-\frac{1}{225}\\
B_{0}(0)&=&\frac{2^{\frac{4}{3}}}{140}\\
B_{1}(0)&=&-\frac{(1213.2)^{\frac{4}{3}}}{2047500}
\end{eqnarray*}

\begin{Rem}[comparaisons des calculs]
Les coefficients mentionnés précédemment sont ceux qui été codés
avec Maple pour pouvoir trouver les points stationnaires et, par
suite, les variations totales. Pour la variation totale de
$\mathcal{V}_{\mathbb{R}}(|\zeta|^{\frac{1}{2}}B_{0}(\zeta))$, il
est possible de la retrouver avec moins de décimales dans de
nombreux documents dont ~\cite{Ol_1962}. Quant à
$\mathcal{V}_{\mathbb{R}}(|\zeta|^{\frac{1}{2}}B_{1}(\zeta))$, elle
peut être retrouvée à la page 315 dans ~\cite{LaWo_1995} où les
auteurs ne donnent pas cependant la fonction $B_{1}(\zeta)$.
\end{Rem}

\noindent Finalement, les routines en Maple pour évaluer les points
stationnaires numériquement sont dont données en annexe.

%

\section{L'article de Kuznetsov et Fedosov, explications et
corrections}

\noindent Ayant maintenant une meilleure idée concernant l'expansion
uniforme de $J_{n}(x)$, alors voici, avec un peu plus détails que
l'article original ~\cite{KuFe_1965}, où il est montré que la loi de
Weyl pour le disque $D$ de rayon $r$ est dans l'ordre $2/3$.

\noindent Sans perte de généralité, le rayon du disque est fixé à
$r=1$. Il est bien connu que les valeurs propres du problème de
Dirichlet sur le disque de rayon 1 sont les zéros positifs des
fonctions de Bessel $J_{n}(x)$ d'ordre $n=0,1,2,\ldots$. Les valeurs
propres correspondant aux ordres $n=1,2,\ldots$ ont une multiplicité
double et celle correspondant à l'ordre $n=0$ sont simples.

\noindent Voici une expansion uniforme de la fonction de Bessel de
premier type obtenue précédemment qui sera utile pour la suite. Soit
$x>n>0$, $J_{n}(x)$ la fonction de Bessel de premier type et soit
\begin{eqnarray}
f_{n}(x)&=&(1+\delta_{3})\sqrt{\frac{\pi}{2}}(x^{2}-n^{2})^{\frac{1}{4}}J_{n}(x)
\end{eqnarray}
avec
\begin{displaymath}
\eta(x,n) = \sqrt{x^{2}-n^{2}}-n\arccos\Big(\frac{n}{x}\Big)
\end{displaymath}
et $\delta_{3}>0$, $\delta_{3}$ ne dépend que de $n$,
$|\delta_{3}|<0.016456$ et $\lim_{n\to\infty}{\delta_{3}}=0$. Alors
pour une valeur fixe de $x$, il existe $n_{x}<x$ telle que
l'inégalité suivante est vérifiée
\begin{eqnarray}
f_{n}(x)&=&\cos\Big(\eta-\frac{\pi}{4}\Big)+O\Big(\frac{1}{n}\Big)+\Big(\frac{1}{\eta}\Big)\label{approx}
\end{eqnarray}
comme il a été démontré précédemment.

\begin{Rem}
Les nombres de zéros de la fonction $f_{n}$ et $J_{n}$ dans un
intervalle donné sont les mêmes, car $f_{n}$ n'est que $J_{n}$
multiplié par une fonction strictement positive.
\end{Rem}

\begin{Rem}
L'approximation donnée de $J_{n}$ donnée dans l'article est comme
suit
\begin{displaymath}
J_{n}(k)=\sqrt{\frac{2}{\pi}}(k^{2}-n^{2})^{-\frac{1}{4}}\bigg(\cos\Big(\eta(k,n)-\frac{\pi}{4}\Big)+O\Big(\frac{1}{\eta}\Big)+O\Big(\frac{1}{n}\Big)\bigg)
\end{displaymath}
\end{Rem}

\noindent Ainsi, le problème consistant au dénombrement du nombre de
zéros positifs des fonctions $J_{n}(x)$ n'excédant pas $k$ est
transformé en un problème de dénombrement de points entiers dans un
domaine bien précis d'où l'utilisation prochaine du théorème de Van
der Corput.

\noindent Soit $N_{n}(k)$ le nombre de zéros de la fonction
$J_{n}(x)$ dans l'intervalle $n<x<k$. Alors,
\begin{displaymath}
N(k)=N_{0}(k)+\sum_{n=1}^{[k]}{N_{n}(k)}
\end{displaymath}
où $[k]$ dénote la partie entière de $k$. En fixant les nombres
$\nu_{0}=[k^{\frac{1}{3}}]+\frac{1}{2}$ et
$\nu_{1}=[k-k^{\frac{4}{9}}]+\frac{1}{2}$, la somme $N(k)-N_{0}(k)$
est divisée en trois parties comme suit
\begin{displaymath}
\sum_{n=1}^{[k]}{N_{n}(k)}=\sum_{n=1}^{\nu_{0}-\frac{1}{2}}{1}+\sum_{n=\nu_{0}+\frac{1}{2}}^{\nu_{1}-\frac{1}{2}}{1}+\sum_{n=\nu_{1}+\frac{1}{2}}^{[k]}{1}=\Sigma_{1}+\Sigma_{2}+\Sigma_{3}.
\end{displaymath}
Une estimation de $\Sigma_{2}$ sera donnée en utilisant le théorème
de Van der Corput. Quant à $\Sigma_{1}$ et $\Sigma_{3}$, leurs
estimations utiliseront la monotonicité de $N_{n}(k)$. $N_{0}(k)$
sera estimé à la toute fin.

\noindent En utilisant (\ref{approx}) pour $n\in(\nu_{0},\nu_{1})$
qui est équivalent à affirmer qu'il existe $A$ et $B$ tels que
\begin{displaymath}
\Big|f_{n}(k)-\cos\Big(\eta(k,n)-\frac{\pi}{4}\Big)\Big|\leq\frac{A}{n}+\frac{B}{\eta(k,n)}
\end{displaymath}
et il s'ensuit que le nombre de zéros de $J_{n}(x)$ pour $x\in(n,k)$
est le même que le nombre de zéros de la fonction
$\cos\big(x-\frac{\pi}{4}\big)$ dans l'intervalle
$0<x\leq\eta(k,n)+\frac{A}{n}+\frac{B}{\eta(k,n)}$. Par conséquent,
\begin{displaymath}
N_{n}(k)=\bigg[\frac{1}{\pi}\eta(k,n)+\frac{1}{4}+\frac{A}{n}+\frac{B}{\eta(k,n)}\bigg]
\end{displaymath}

\begin{Rem}[Justification du nombre de zéros]
Une justification de la formule pour $N_{n}(k)$ a été donnée
rigourousement par Dominique Rabet. Ses résultats paraîtront
ultérieurement dans ses publications.
\end{Rem}

\noindent Dans cette partie, une estimation de $\Sigma_{2}$ est
donnée en détails. Il faut utiliser le théorème de Van der Corput
mentionné dans \cite{KuFe_1965} que voici
\begin{Thm}[Van der Corput]
Soit les nombres $\nu_{0}-\frac{1}{2}$, $\nu_{1}-\frac{1}{2}$ et
$x_{0}-\frac{1}{2}$ entiers tel que $\nu_{0}<\nu_{1}$. Soit la
fonction à valeurs réelles $f(\nu)$ telle que
$f\in\mathrm{C}^{2}(\nu_{0},\nu_{1})$ et telle que
\begin{displaymath}
0<\sigma\leq f'(\nu) \leq
\tau,\phantom{12}f''(\nu)>\frac{1}{\mu},\phantom{12}\mu>1,\mu>\sigma^{-3}.
\end{displaymath}
Soit $N$ le nombre de points entiers à l'intérieur du domaine
$\nu_{0}<\nu<\nu_{1}$ et $x_{0}\leq x \leq f(\nu)$. Soit $A$ l'aire
u domaine. Alors
\begin{displaymath}
N=A+O(\mu^{\frac{2}{3}}\tau)
\end{displaymath}
\end{Thm}

\noindent Par ce qu'il a été accompli auparavant, $\Sigma_{2}$ est
donné par le nombre de points entiers dans le domaine que voici
\begin{displaymath}
\nu_{0}<\nu<\nu_{1},\phantom{12}\frac{1}{2}<x<f(\nu)=\frac{1}{\pi}\eta(k,\nu)+\frac{1}{4}+\frac{A}{\nu}+\frac{B}{\eta(k,\nu)}.
\end{displaymath}

\noindent Soit $\delta=\max{\big|\frac{A}{\nu}\big|}$ et $a$ telle
que $\big|\frac{B}{\eta(k,\nu)}\big|<\frac{a}{\pi\eta}$ pour
$\nu\in(\nu_{0},\nu_{1})$. Soit
\begin{eqnarray*}
f_{1}(\nu)&=&\frac{1}{\pi}\eta+\frac{1}{4}-\frac{a}{\pi\eta}-\delta\\
f_{2}(\nu)&=&\frac{1}{\pi}\eta+\frac{1}{4}+\frac{a}{\pi\eta}+\delta
\end{eqnarray*}
alors $f_{1}(\nu)\leq f(\nu) \leq f_{2}(\nu)$. Par conséquent, en
dénotant $N_{1}$ et $N_{2}$ les nombres de points entiers donnés en
remplaçant $f$ respectivement par $f_{1}$ et $f_{2}$, alors
$N_{1}\leq \Sigma_{2}\leq N_{2}$.

\noindent Le théorème de Van der Corput sera utilisé sur $f_{1}$ et
$f_{2}$. Il faut montrer en détails que $f_{1}$ et $f_{2}$ satisfont
aux hypothèses du théorème. Il est suffisant de le faire pour
$f_{1}$, le cas de $f_{2}$ étant identique.

\noindent Pour montrer que $f_{1}$ satisfait aux hypothèses du
théorème, il faut d'abord évaluer les dérivées suivantes.
\begin{eqnarray*}
\frac{\partial\eta}{\partial\nu}&=&\frac{\partial}{\partial\nu}\Big(\sqrt{k^{2}-\nu^{2}}-\nu\arccos\Big(\frac{\nu}{k}\Big)\Big)\\
&=&\frac{1}{2}\frac{1}{\sqrt{k^{2}-\nu^{2}}}(-2\nu)-\arccos\Big(\frac{\nu}{k}\Big)+\nu\frac{1}{\sqrt{1-\frac{\nu^{2}}{k^{2}}}}\frac{1}{k}\\
&=&\frac{-1}{\sqrt{k^{2}-\nu^{2}}}-\arccos\Big(\frac{\nu}{k}\Big)+\frac{\nu}{k}\frac{k}{\sqrt{k^{2}-\nu^{2}}}\\
&=&-\arccos\Big(\frac{\nu}{k}\Big)\\
\frac{\partial\eta^{2}}{\partial\nu^{2}}&=&\frac{k}{\sqrt{k^{2}-\nu^{2}}}
\end{eqnarray*}

\noindent De même,
\begin{eqnarray*}
\frac{\partial{f_{1}}}{\partial\nu}&=&\frac{1}{\pi}\eta_{\nu}(k,\nu)+\frac{a}{\pi}\frac{\eta_{\nu}(k,\nu)}{\eta^{2}(k,\nu)}\\
&=&\frac{1}{\pi}\eta_{\nu}(k,\nu)\bigg(1+\frac{a}{\eta^{2}(k,\nu)}\bigg)\\
&=&-\frac{\arccos\big(\frac{\nu}{k}\big)}{\pi}\bigg(1+\frac{a}{\eta^{2}(k,\nu)}\bigg)
\end{eqnarray*}
et
\begin{eqnarray*}
\frac{\partial^{2}f_{1}}{\partial\nu^{2}}&=&-\frac{1}{\pi}\bigg(\frac{-k}{\sqrt{k^{2}-\nu^{2}}}\bigg)\bigg(1+\frac{a}{\eta^{2}(k,\nu)}\bigg)-\frac{\arccos\big(\frac{\nu}{k}\big)}{k}\bigg(0-\frac{2a\eta_{\nu}(k,\nu)}{\eta^{3}(k,\nu)}\bigg)\\
&=&\frac{k}{\pi}\frac{1}{\sqrt{k^{2}-\nu^{2}}}\bigg(1+\frac{a}{\eta^{2}(k,\nu)}\bigg)-\frac{\Big(\arccos\big(\frac{\nu}{k}\big)\Big)^{2}2a}{\pi\eta^{3}(k,\nu)}
\end{eqnarray*}

\noindent Comme $\eta(k,\nu)\searrow$ lorsque $\nu\nearrow$ avec $k$
fixe et que $\eta>0$, alors
$\min\eta(k,\nu)=\eta(k,\nu_{1})=\eta(k,k-k^{\frac{4}{9}})$ pour
$\nu\leq\nu_{1}$. L'inégalité $\arccos(1-\alpha)\geq\sqrt{\alpha}$
pour $\alpha\in(0,1)$ permet donc de déduire l'inégalité utile
suivante non démontrée dans l'article \cite{KuFe_1965}.
\begin{eqnarray*}
\eta(k,k-k^{\frac{4}{9}})&=&\sqrt{k^{2}-(k^{2}-2kk^{\frac{4}{9}}+k^{\frac{8}{9}})}-(k-k^{\frac{4}{9}})\arccos\Big(\frac{k-k^{\frac{4}{9}}}{k}\Big)\\
&\geq&\sqrt{2k^{\frac{13}{9}}-k^{\frac{8}{9}}}-(k-k^{\frac{4}{9}})\sqrt{k^{-\frac{5}{9}}}\\
&=&\sqrt{2k^{\frac{26}{18}}-k^{\frac{10}{18}}}-k^{\frac{3}{18}}(k^{\frac{10}{18}}-1)\\
&=&\sqrt{2k^{\frac{6}{18}}\Big(k^{\frac{20}{18}}-\frac{1}{2}k^{\frac{10}{18}}\Big)}-\sqrt{2}k^{\frac{3}{18}}(k^{\frac{10}{18}}-1)\frac{1}{\sqrt{2}}\\
&=&\sqrt{2k^{\frac{3}{18}}}\bigg(\sqrt{k^{\frac{20}{18}}-\frac{1}{2}k^{\frac{10}{18}}}-\frac{1}{\sqrt{2}}(k^{\frac{10}{18}}-1)\bigg)\\
&\geq&\sqrt{2}k^{\frac{1}{6}}\\
&\geq&k^{\frac{1}{6}}
\end{eqnarray*}
De la même façon,
\begin{displaymath}
\eta(k,k-k^{\frac{3}{5}})\geq\sqrt{2}k^{\frac{2}{5}}\geq{k^{\frac{2}{5}}}
\end{displaymath}

\noindent Également, l'approximation suivante pour
$\tilde{\eta}(z)=\sqrt{1-z^{2}}-z\arccos(z)$ est utile
\begin{displaymath}
1\leq\frac{(1-z)^{\frac{3}{2}}}{\sqrt{1-z^{2}}-z\arccos(z)}<\frac{3}{4}\sqrt{2}\phantom{12}z\in(0,1)
\end{displaymath}
Par example, avec $0<z=\frac{\nu}{k}<1$, l'inégalité suivante est
vérifiée
\begin{displaymath}
\frac{1}{\eta(k,\nu)}=\frac{1}{k\big(\sqrt{1-z^{2}}-z\arccos(z)\big)}\leq\frac{3}{4}\sqrt{2}\frac{1}{k}(1-z)^{-\frac{3}{2}}
\end{displaymath}

\noindent Puisque
$\arccos\big(\frac{\nu_{1}}{k}\big)\leq\arccos\big(\frac{\nu}{k}\big)$,
alors
\begin{displaymath}
k^{-\frac{5}{18}}=\sqrt{k^{-\frac{5}{9}}}\leq\arccos(1-k^{-\frac{5}{9}})=\arccos\Big(\frac{\nu_{1}}{k}\Big)\leq\arccos\Big(\frac{\nu}{k}\Big).
\end{displaymath}
Par conséquent,
\begin{displaymath}
0<\frac{k^{-\frac{5}{18}}}{\pi}\leq\frac{k^{-\frac{5}{18}}}{\pi}\bigg(1+\frac{a}{\eta^{2}}\bigg)\leq-\frac{\partial{f_{1}}}{\partial\nu}
\end{displaymath}
De même, comme
$0\leq\arccos\big(\frac{\nu}{k}\big)\leq\frac{\pi}{2}$, alors
\begin{displaymath}
-\frac{\partial{f_{1}}}{\partial\nu}\leq\frac{1}{2}\bigg(1+\frac{a}{\eta^{2}}\bigg)\leq\frac{1}{2}+\frac{a}{\eta^{2}(k,k-k^{\frac{4}{9}})}\leq\frac{1}{2}+\frac{a}{(k^{\frac{1}{6}})^{2}}=\frac{1}{2}+ak^{-\frac{1}{3}}
\end{displaymath}
Il s'ensuit donc qu'avec $C=\frac{1}{\pi}$
\begin{displaymath}
0<Ck^{-\frac{5}{18}}<-\frac{\partial{f_{1}}}{\partial\nu}\leq\frac{1}{2}+O(k^{-\frac{1}{3}})
\end{displaymath}
tel qu'écrit dans l'article \cite{KuFe_1965}.

\noindent Pour estimer $\frac{\partial^{2}f_{1}}{\partial\nu^{2}}$,
il faut d'abord montrer que
\begin{displaymath}
\frac{\Big(\arccos\big(\frac{\nu}{k}\big)\Big)^{2}}{\pi\eta^{3}(k,\nu)}=o\bigg(\frac{1}{\sqrt{k^{2}-\nu^{2}}}\bigg)
\end{displaymath}
uniformément pour $\nu\in(\nu_{0},\nu_{1})$ lorsque $k\to\infty$.
Pour ce faire, les auteurs séparent l'intervalle
$(\nu_{0},\nu_{1})=(\nu_{0},\nu')\cup(\nu',\nu_{1})$ où
$\nu'=k-k^{\frac{3}{5}}$. Dans l'intervalle $(\nu_{0},\nu')$, comme
$\eta(k,\nu)$ est toujours décroissante, alors
\begin{displaymath}
\frac{\Big(\arccos\big(\frac{\nu}{k}\big)\Big)^{2}}{\eta^{3}}\leq\frac{\Big(\arccos\big(\frac{\nu_{0}}{k}\big)\Big)^{2}}{\eta^{3}}\leq\frac{\pi^{2}}{4}\frac{1}{\eta^{3}(k,\nu')}\leq\frac{\pi^{2}}{4}\bigg(\frac{1}{k^{\frac{2}{5}}}\bigg)^{3}=\frac{\pi^{2}}{4}\frac{1}{k^{\frac{6}{5}}}\ldots
\end{displaymath}
\begin{displaymath}
\ldots\leq\frac{\pi^{2}}{4}\frac{1}{k}\leq\frac{\pi^{2}}{4}\frac{1}{\sqrt{k^{2}-\nu^{2}}}
\end{displaymath}
et, dans l'intervalle $(\nu',\nu_{1})$,
\begin{displaymath}
\frac{\Big(\arccos\big(\frac{\nu}{k}\big)\Big)^{2}}{\eta^{3}}\leq\frac{\Big(\arccos\big(\frac{\nu'}{k}\big)\Big)^{2}}{\eta^{3}(k,\nu_{1})}\leq\frac{\pi^{2}}{4}\bigg(\frac{1}{k\tilde{\eta}\big(\frac{\nu}{k}\big)}\bigg)^{3}\leq\frac{\pi^{2}}{4}\Bigg(\frac{3}{4}\sqrt{2}\frac{1}{k}\Big(1-\frac{\nu}{k}\Big)^{-\frac{3}{2}}\Bigg)^{3}\ldots
\end{displaymath}
\begin{displaymath}
\phantom{12345}\ldots\leq\frac{\pi^{2}}{4}\Bigg(\frac{3}{4}\sqrt{2}\frac{1}{k}\Big(1-\frac{\nu'}{k}\Big)^{-\frac{3}{2}}\Bigg)^{3}=\frac{\pi^{2}}{4}\Bigg(\frac{3}{4}\sqrt{2}\frac{1}{k}\Big(k^{-\frac{2}{5}}\Big)^{-\frac{3}{2}}\Bigg)^{3}=(\mathrm{const.})k^{-\frac{6}{5}}
\end{displaymath}
Par conséquent,
\begin{displaymath}
\frac{\partial^{2}f_{1}}{\partial\nu^{2}}=\frac{1}{\pi\sqrt{k^{2}-\nu^{2}}}+o\bigg(\frac{1}{\sqrt{k^{2}-\nu^{2}}}\bigg)>\frac{D}{k}
\end{displaymath}

\begin{Rem}[Satisfaction des hypothèses du théorème de Van der Corput]
Comme
\begin{displaymath}
0<Ck^{-\frac{5}{18}}<-\frac{\partial{f_{1}}}{\partial\nu}\leq\frac{1}{2}+O(k^{-\frac{1}{3}})
\end{displaymath}
et que
\begin{displaymath}
\frac{\partial^{2}f_{1}}{\partial\nu^{2}}>\frac{D}{k}
\end{displaymath}
alors $f_{1}$ satisfait aux hypothèses du théorème de Van der Corput
avec
\begin{eqnarray*}
\sigma&=&Ck^{-\frac{5}{18}}\\
\tau&=&\frac{1}{2}+O(k^{-\frac{1}{3}})\\
\mu&=&\frac{k}{D}.
\end{eqnarray*}
\end{Rem}

\noindent En dénotant ainsi par $N_{1}$ et $A_{1}$ respectivement le
nombre de points entiers et l'aire du domaine délimité par
$\nu_{0}<\nu<\nu_{1}$ et $x_{0}\leq x\leq f_{1}(\nu)$ et en
applicant le théorème de Van der Corput,
\begin{displaymath}
N_{1}=A_{1}+O(\mu^{\frac{2}{3}}\tau)=\int_{\nu_{0}}^{\nu_{1}}{\Big(f_{1}(\nu)-\frac{1}{2}\Big)d\nu}+O(k^{\frac{2}{3}}).
\end{displaymath}
Comme $\delta=\max\big(\frac{A}{\nu}\big)$, alors
\begin{displaymath}
\int_{\nu_{0}}^{\nu_{1}}{\delta{d\nu}}\leq(\nu_{1}-\nu_{0})\frac{A}{\nu_{0}}\leq\frac{\nu_{1}}{\nu_{0}}=O(k^{\frac{2}{3}}).
\end{displaymath}
Aussi, puisque $\eta_{\nu}=-\arccos\big(\frac{\nu}{k}\big)$, que
$k^{-\frac{5}{18}}\leq\arccos\big(\frac{\nu}{k}\big)$ pour
$\nu\in(\nu_{0},\nu_{1})$, alors
\begin{displaymath}
\int_{\nu_{0}}^{\nu_{1}}{\frac{1}{\eta}d\nu}=\int_{\nu_{0}}^{\nu_{1}}{\frac{1}{\eta}\frac{\eta_{\nu}}{\eta_{\nu}}d\nu}=\int_{\eta(k,\nu_{0})}^{\eta(k,\nu_{1})}{\frac{d\eta}{\eta\eta_{\nu}}}=\int_{\eta(k,\nu_{0})}^{\eta(k,\nu_{1})}{\frac{d\eta}{-\arccos\big(\frac{\nu}{k}\big)\eta}}=\ldots
\end{displaymath}
\begin{displaymath}
\ldots=\int_{\eta(k,\nu_{1})}^{\eta(k,\nu_{0})}{\frac{d\eta}{\arccos\big(\frac{\nu}{k}\big)\eta}}\leq\frac{1}{k^{-\frac{5}{18}}}\log\bigg(\frac{\eta(k,\nu_{0})}{\eta(k,\nu_{1})}\bigg)=O(k^{\frac{5}{18}}\log{k}).
\end{displaymath}

\begin{Rem}[correction]
Dans l'article, une erreur typographique s'est produite. En effet,
il est écrit que $N_{1}=O(k^{-\frac{5}{18}}\log{k})$.
\end{Rem}

\noindent Ainsi,
\begin{eqnarray*}
N_{1}&=&\int_{\nu_{0}}^{\nu_{1}}{\Big(f_{1}(\nu)-\frac{1}{2}\Big)d\nu}+O(k^{\frac{2}{3}})\\
&=&\int_{\nu_{0}}^{\nu_{1}}{\Big(\frac{\eta}{\pi}+\frac{1}{4}-\frac{a}{\pi\eta}-\delta-\frac{1}{2}\Big)d\nu}+O(k^{\frac{2}{3}})\\
&=&\int_{\nu_{0}}^{\nu_{1}}{\Big(\frac{\eta}{\pi}-\frac{1}{4}\Big)d\nu}+O(k^{\frac{2}{3}})+O(k^{\frac{5}{18}}\log{k})+O(k^{\frac{2}{3}})\\
&=&\int_{\nu_{0}}^{\nu_{1}}{\Big(\frac{\eta}{\pi}-\frac{1}{4}\Big)d\nu}+O(k^{\frac{2}{3}})
\end{eqnarray*}
Il est de même pour $N_{2}$ en montrant de la même façon que $f_{2}$
satisfait aux hypothèses du théorème de Van der Corput. Par
conséquent, $N_{1}\leq\Sigma_{2}\leq N_{2}$ et
\begin{displaymath}
\Sigma_{2}=\int_{\nu_{0}}^{\nu_{1}}{\Big(\frac{\eta}{\pi}-\frac{1}{4}\Big)d\nu}+O(k^{\frac{2}{3}})
\end{displaymath}

\noindent Maintenant pour estimer $\Sigma_{1}$ et $\Sigma_{2}$, il
faut utiliser la monotonicité de $N_{n}(k)$. De façon générale,
\begin{displaymath}
N_{n}(k)\leq{N_{n-1}(k)}\phantom{12}\textrm{pour tout}\phantom{1}n,k
\end{displaymath}

\noindent Pour $\Sigma_{3}$, comme
\begin{displaymath}
\Sigma_{3}=\sum_{n=\nu_{1}+\frac{1}{2}}^{[k]}{N_{n}(k)}=\sum_{[k-k^{\frac{4}{9}}]+1}^{[k]}{N_{n}(k)}\leq\sum_{[k-k^{\frac{4}{9}}]+1}^{[k]}{N_{\nu_{1}-\frac{1}{2}}(k)}\leq
k^{\frac{4}{9}}N_{\nu_{1}-\frac{1}{2}}(k)
\end{displaymath}
que
\begin{displaymath}
N_{\nu_{1}-\frac{1}{2}}(k)=\bigg[\frac{\eta(k,k-k^{\frac{4}{9}})}{\pi}+\frac{1}{4}+\frac{A}{n}+\frac{B}{\eta}\bigg]=O(\eta(k,k-k^{\frac{4}{9}}))=O(k^{\frac{1}{6}})
\end{displaymath}
que $\eta(k,k-k^{\frac{4}{9}})=k\tilde{\eta}(1-k^{-\frac{5}{9}})$ et
que
\begin{eqnarray*}
\tilde{\eta}(1-k^{-\frac{5}{9}})&<&\Big(1-\big(1-k^{-\frac{5}{9}}\big)\Big)^{\frac{3}{2}}\\
&=&k^{-\frac{5}{6}}\\
&\Downarrow&\\
\eta(k,k-k^{\frac{4}{9}})&=&k\tilde{\eta}(1-k^{-\frac{5}{9}})\\
&<&k^{\frac{1}{6}}
\end{eqnarray*}
alors
\begin{displaymath}
\Sigma_{3}=O(k^{\frac{4}{9}}k^{\frac{1}{6}})=O(k^{\frac{11}{18}})
\end{displaymath}

\noindent Pour $\Sigma_{1}$, il faut téléscoper la somme d'abord
comme suit
\begin{eqnarray*}
\Sigma_{1}&=&\sum_{n=1}^{\nu_{0}-\frac{1}{2}}{N_{n}(k)}\\
&=&\Big(\nu_{0}-\frac{1}{2}\Big)N_{\nu_{0}-\frac{1}{2}}(k)+\Big(N_{1}(k)-N_{\nu_{0}-\frac{1}{2}}(k)\Big)\\
&&\phantom{\big(\nu_{0}-\frac{1}{2}\big)N_{\nu_{0}-\frac{1}{2}}(k)}+\Big(N_{2}(k)-N_{\nu_{0}-\frac{1}{2}}(k)\Big)\\
&&\phantom{\big(\nu_{0}-\frac{1}{2}\big)N_{\nu_{0}-\frac{1}{2}}(k)}\vdots\\
&&\phantom{\big(\nu_{0}-\frac{1}{2}\big)N_{\nu_{0}-\frac{1}{2}}(k)}+\Big(N_{\nu_{0}-\frac{1}{2}}(k)-N_{\nu_{0}-\frac{1}{2}}(k)\Big)\\
&=&\Big(\nu_{0}-\frac{1}{2}\Big)N_{\nu_{0}-\frac{1}{2}}(k)+\sum_{n=1}^{\nu_{0}-\frac{1}{2}}{\Big(N_{n}(k)-N_{\nu_{0}-\frac{1}{2}}(k)\Big)}.
\end{eqnarray*}
Maintenant il faut estimer $N_{0}(k)$ et
$N_{\nu_{0}-\frac{1}{2}}(k)$. Dans le cas de $N_{0}(k)$, comme
$J_{0}(k)\approx\sqrt{\frac{2}{\pi{z}}}\cos\big(z-\frac{\pi}{4}\big)$,
alors $N_{0}(k)=\frac{k}{\pi}+O(1)$ et
$N_{\nu_{0}-\frac{1}{2}}(k)=\frac{k}{\pi}+O(k^{\frac{1}{3}})$. Par
la monotonicité de $N_{n}(k)$,
\begin{eqnarray*}
N_{n}(k)&<&N_{0}(k)\\
&\Downarrow&\\
N_{n}(k)-N_{\nu_{0}-\frac{1}{2}}(k)&<&N_{0}(k)-N_{\nu_{0}-\frac{1}{2}}(k)\\
&=&\frac{k}{\pi}+O(1)-\frac{k}{\pi}+O(k^{\frac{1}{3}})\\
&=&O(k^{\frac{1}{3}})
\end{eqnarray*}
et, par conséquent,
\begin{eqnarray*}
\Sigma_{1}&=&\Big(\nu_{0}-\frac{1}{2}\Big)N_{\nu_{0}-\frac{1}{2}}(k)+\sum_{n=1}^{\nu_{0}-\frac{1}{2}}{\Big(N_{n}(k)-N_{\nu_{0}-\frac{1}{2}}(k)\Big)}\\
&=&\Big(\nu_{0}-\frac{1}{2}\Big)\Big(\frac{k}{\pi}+O(k^{\frac{1}{3}})\Big)+\Big(\nu_{0}-\frac{1}{2}\Big)O(k^{\frac{1}{3}})\\
&=&\frac{k}{\pi}\Big(\nu_{0}-\frac{1}{2}\Big)+O(k^{\frac{2}{3}})
\end{eqnarray*}

\begin{Rem}[correction]
Des erreurs typographiques pour l'expression de $\Sigma_{1}$ se sont
produites dans l'article. En effet, il est écrit que
\begin{displaymath}
\sum_{n=1}^{\nu_{0}-\frac{1}{2}}{N_{n}(k)}=N_{\nu_{0}-\frac{1}{2}}\Big(\nu_{0}-\frac{1}{2}\Big)+\sum_{n=1}^{\nu_{0}-\frac{1}{2}}{\Big(N_{n}(k)-N_{\nu_{0}-\frac{1}{2}}(k)\Big)}
\end{displaymath}
qui devrait s'écrire
\begin{displaymath}
\sum_{n=1}^{\nu_{0}-\frac{1}{2}}{N_{n}(k)}=\Big(\nu_{0}-\frac{1}{2}\Big)N_{\nu_{0}-\frac{1}{2}}(k)+\sum_{n=1}^{\nu_{0}-\frac{1}{2}}{\Big(N_{n}(k)-N_{\nu_{0}-\frac{1}{2}}(k)\Big)}
\end{displaymath}
\end{Rem}

\noindent Ainsi,
\begin{eqnarray*}
N(k)&=&N_{0}(k)+2(\Sigma_{1}+\Sigma_{2}+\Sigma_{3})\\
&=&\frac{k}{\pi}+O(1)+2\bigg(\Big(\nu_{0}-\frac{1}{2}\Big)\frac{k}{\pi}+O(k^{\frac{2}{3}})+\int_{\nu_{0}}^{\nu_{1}}{\Big(\frac{\eta}{\pi}-\frac{1}{4}\Big)d\nu}+\Sigma_{3}\bigg)\\
&=&\frac{k}{\pi}+2\nu_{0}\frac{k}{\pi}-\frac{k}{\pi}+2\int_{\nu_{0}}^{\nu_{1}}{\frac{\eta}{\pi}d\nu}-\frac{2}{4}(\nu_{1}-\nu_{0})+2\Sigma_{3}+O(k^{\frac{2}{3}})\\
&=&\frac{2\nu_{0}k}{\pi}+2\int_{\nu_{0}}^{\nu_{1}}{\frac{\eta}{\pi}d\nu}-\frac{1}{2}(k-k^{\frac{4}{9}}-k^{\frac{1}{3}})+2\Sigma_{3}+O(k^{\frac{2}{3}})\\
&=&\frac{2}{\pi}\nu_{0}k+\frac{2}{\pi}\int_{\nu_{0}}^{\nu_{1}}{\eta(k,\nu)d\nu}-\frac{k}{2}+2\Sigma_{3}+O(k^{\frac{2}{3}})\\
&=&\frac{2}{\pi}(k\nu_{0}+\nu_{0}O(k^{\frac{1}{3}}))+\frac{2}{\pi}\int_{\nu_{0}}^{\nu_{1}}{\eta(k,\nu)d\nu}-\frac{k}{2}+\frac{2}{\pi}\int_{\nu_{1}}^{k}{\eta(k,\nu)d\nu}+O(k^{\frac{2}{3}})\\
&=&\frac{2}{\pi}\int_{0}^{\nu_{0}}{\eta(k,\nu)d\nu}+\frac{2}{\pi}\int_{\nu_{0}}^{\nu_{1}}{\eta(k,\nu)d\nu}-\frac{k}{2}+\frac{2}{\pi}\int_{\nu_{1}}^{k}{\eta(k,\nu)d\nu}+O(k^{\frac{2}{3}})\\
&=&\frac{2}{\pi}\int_{0}^{k}{\eta(k,\nu)d\nu}-\frac{k}{2}+O(k^{\frac{2}{3}})\\
&=&\frac{2}{\pi}\Big(\frac{\pi{k^{2}}}{8}\Big)-\frac{k}{2}+O(k^{\frac{2}{3}})\\
&=&\frac{k^{2}}{4}-\frac{k}{2}+O(k^{\frac{2}{3}})
\end{eqnarray*}

\noindent Ce qui termine l'explication de \cite{KuFe_1965}.

\noindent \textbf{Nous rappelons les lecteurs que spécialement pour
la section 3.2 qui s'achève ici, nous avons utilisé $k$ au lieu de
$\lambda$ et étudier le problème $\triangle{u}+k^{2}{u}=0$ plutôt
que $\triangle{u}+k{u}=0$ et cela dans le but de faciliter les
références à l'article \cite{KuFe_1965}. Pour la section 3.3 qui
suit, nous revenons à notre lettre habituelle $\lambda$ et au
problème des valeurs propres exprimé sous la forme
$\triangle{u}+\lambda{u}=0$.}

\section{Algorithme pour calculer la fonction de compte \texorpdfstring{N($\lambda$)}{}}

\noindent Cette section est très courte et comporte une idée que
j'ai eue permettant de calculer la fonction de compte
\begin{displaymath}
N(\lambda)=\sum_{\lambda_{j}\leq\lambda}{1}
\end{displaymath}
efficacement \textbf{sans avoir à calculer toutes les valeurs
propres $\lambda_{j}$ telles que $\lambda_{j}\leq\lambda$}.

\begin{Rem}[Façon naïve de calculer $N(\lambda)$]
Étant donné l'algorithme de la section 2.2 qui détermine la suite
$\{\lambda_{j}\}_{j=1}^{l}$ pour un entier $l$ donné. La façon la
plus simple de calculer $N(\lambda)$ est de calculer
$\{\lambda_{j}\}_{j=1}^{l}$ pour un $l$ suffisamment grand et par la
suite de compter combien de valeurs propres $\lambda_{j}$ sont
inférieures à $\lambda$.
\end{Rem}

\noindent \textbf{\og Buts \fg{} des algorithmes}\\Étant donné un
entier strictement positif $l$, l'algorithme de la section 2.2 donne
$\{\lambda_{j}\}_{j=1}^{l}$. Étant donné un nombre réel strictement
positif $\lambda$, l'algorithme de la présente section donne
l'entier positif $N(\lambda)$. Toutefois de par les remarques
précédentes, nous pouvons voir qu'il est possible d'utiliser
l'algorithme de la section 2.2 pour évaluer $N(\lambda)$ et cela a
un sens en se rappelant les rappels théoriques de la section 1.1
concernant $N(\lambda)$.

\begin{Rem}
Les deux algorithmes, celui de la présente section et celui de la
section 2.2, déterminent un ensemble fini de solutions
$\mathrm{SOL}$ de paires $(n,k)$ avec $n=0,1,2,\ldots$ et
$k=1,2,\ldots$ où $n$ et $k$ sont les paramètres pour les zéros. Les
paramètres sont l'ordre ($n$) de la fonction de Bessel et le rang
($k$) du zéro. Dans le cas de l'algorithme de la section 2.2, un
ensemble fini $\mathrm{PA}$ contenant $\mathrm{SOL}$ est déterminé à
l'aide du théorème de Courant et $\mathrm{SOL}$ est obtenu en
ordonnant les zéros correspondant aux paires dans $\mathrm{PA}$.
Dans le cas de l'algorithme de la présente section, seul la
frontière de $\mathrm{SOL}$ est évaluée et il n'est pas nécessaire
de trouver tous les zéros qui sont plus petits que $\sqrt{\lambda}$
pour évaluer $N(\lambda)$.
\end{Rem}

\noindent L'idée utilise la monotonicité des zéros des fonctions de
Bessel $J_{\nu}(x)$ pour $x\geq 0$ et $\nu\geq 0$. Si $j_{k}(\nu)$
dénote le k\ieme{} zéro de la fonction de Bessel de $J_{\nu}(x)$,
alors ( voir \cite{Wa_1944} ),
\begin{eqnarray*}
j_{k}(\nu)&<&j_{k+1}(\nu)\\
j_{k}(\nu)&<&j_{k}(\nu+\nu')\phantom{1}\textrm{pour $\nu'>0$.}
\end{eqnarray*}

\noindent Les valeurs propres de $\triangle$ sur le disque $D$ de
rayon 1 correspondent aux carrées des zéros $j_{k}(n)$ avec
$k=1,2,\ldots$ et $n=0,1,\ldots$. En fixant une valeur de $\lambda$
et en dénotant par $k_{n,\lambda}$, la valeur entière telle que
\begin{displaymath}
j_{k_{n,\lambda}}(n)\leq \sqrt{\lambda}
\phantom{1}\textrm{et}\phantom{1}j_{k_{n,\lambda}+1}(n)>
\sqrt{\lambda}.
\end{displaymath}
alors, par la monotonicité des zéros,
\begin{displaymath}
j_{k_{n,\lambda}}(n) < j_{k_{n+1,\lambda}}(n+1).
\end{displaymath}

\noindent Il s'agit donc pour chaque ordre $n$ et $\lambda$ choisi a
priori d'évaluer $k_{n,\lambda}$ et ainsi en tenant compte de la
multiplicité
\begin{displaymath}
N(\lambda)=k_{0,\lambda}+2\sum_{n=1}^{N}{k_{n,\lambda}}
\end{displaymath}
où
\begin{displaymath}
N = \max\big\{n\phantom{1};\phantom{1}\textrm{tel
que}\phantom{1}k_{n,\lambda}=1\big\}.
\end{displaymath}

\begin{Rem}
L'ensemble des paires $(n,k)$ solution du problème dénoté auparavant
par $\mathrm{SOL}$ est ainsi donc donné par
\begin{displaymath}
\mathrm{SOL}=\big\{(n,k)\phantom{1},\phantom{1}k=1,\ldots,k_{n,\lambda},\phantom{1}n=0,1,\ldots,N\big\}
\end{displaymath}
et
\begin{eqnarray*}
N(\lambda)&=&2\cdot\mathrm{card}\big(\mathrm{SOL}\setminus\big\{(n,k)\phantom{1},\phantom{1}n=1,2,\ldots,N,\phantom{1}k=1,2,\ldots,k_{n,\lambda}\}\big)\\
&&\phantom{123}+\mathrm{card}\big(\big\{(0,k)\phantom{1},\phantom{1}k=1,2,\ldots,k_{0,\lambda}\}\big)\\
&=&k_{0,\lambda}+2\sum_{n=1}^{N}{k_{n,\lambda}}
\end{eqnarray*}
Dans le présent algorithme qui calcule $N(\lambda)$, il n'est donc
pas nécessaire d'évaluer tous les zéros associés aux paires dans
$\mathrm{SOL}$, mais seulement les valeurs $k_{n,\lambda}$ qui
déterminent la frontière de $\mathrm{SOL}$. Ce n'est pas le cas pour
l'algorithme de la section 2.2 où il faut tous les évaluer afin de
créer $\{\lambda_{j}\}_{j=1}^{l}$ et cela s'explique parce que les
buts poursuivis par les deux algorithmes sont différents quoique
mathématiquement équivalents par les rappels théoriques de la
section 1.1.
\end{Rem}

\noindent \textbf{L'algorithme en mots}\\Les valeurs $N$ et
$k_{0,\lambda}$ peuvent être trouvées par une méthode quelconque.
Soit $n=0,1,\ldots$, alors
$k_{n,\lambda}-k_{n+1,\lambda}\in\{0,1\}$. Donc pour trouver
$k_{n+1,\lambda}$ en possédant $k_{n,\lambda}$, il suffit de tester
la validité de l'inégalité suivante
\begin{displaymath}
j_{k_{n,\lambda}}(n+1)\leq\sqrt{\lambda}.
\end{displaymath}
Si l'inégalité est vraie, alors $k_{n+1,\lambda}=k_{n,\lambda}$
sinon $k_{n+1,\lambda}=k_{n,\lambda}-1$, etc. L'algorithme est
initialisé avec $k_{0,\lambda}$ et l'algorithme procède en reculant
et montant jusqu'à $N$.


\noindent L'algorithme en Maple est annexé à la fin de ce mémoire.

\begin{Rem}[Précision requise dans
l'algorithme]\label{remprecalgocompte} Si nous voulons calculer
$N(\lambda)$, alors nous voulons trouver tous les zéros
$j_{k}(n)\leq\sqrt{\lambda}$. Donc par la remarque \ref{remcapoute1}
avec $\lambda=9.00\cdot{10}^{08}$ qui est la valeur maximale de
$\lambda$ apparaissant dans le tableau \ref{bigTabTrueCount}, il
faut demander $\lceil\log(9.00\cdot{10}^{08})\rceil+1=10$ décimales
de précision. Également, le lecteur peut se référer à l'annexe A.5
pour plus d'information.
\end{Rem}


\subsection{Quelques résultats pour les valeurs de la fonction de compte \texorpdfstring{N($\lambda$)}{}}

\noindent Voici quelques résultats pour les valeurs de la fonction
de compte. Il est intéressant de constater que pour les valeurs de
$\lambda=1\cdot{10}^{2},\ldots,4\cdot{10}^{6}$ dans le tableau qui
suit, les valeurs de $N(\lambda)$ coïncident avec les valeurs
retouvées dans la liste des valeurs propres calculées par le premier
algorithme. Par exemple à l'aide de la liste donnée à l'hyperlien
précédent, la valeur propre $\lambda_{998978}$ est la dernière
valeur propre plus petite que $4.00\cdot{10}^{6}$ et le rang
correspond également à la valeur apparaissant dans le tableau qui
suit.

\noindent Voici donc quelques valeurs de la fonction de compte.
\begin{table}[ht]
\begin{center}
\begin{tabular}{|c|c||c|c||c|c|}
\hline $\lambda$ & $N(\lambda)$ & $\lambda$ & $N(\lambda)$ &
$\lambda$ & $N(\lambda)$\\
\hline 1.00E+02 & 21 &4.00E+04 & 9905 &7.00E+06 & 1748690\\
\hline 2.00E+02 & 42 &5.00E+04 & 12387 &8.00E+06 & 1998600\\
\hline 3.00E+02 & 67 &6.00E+04 & 14876 &9.00E+06 & 2248481\\
\hline 4.00E+02 & 92 &7.00E+04 & 17360 &1.00E+07 & 2498404\\
\hline 5.00E+02 & 115 &8.00E+04 & 19852 &2.00E+07 & 4997783\\
\hline 6.00E+02 & 142 &9.00E+04 & 22345 &3.00E+07 & 7497255\\
\hline 7.00E+02 & 160 &1.00E+05 & 24842 &4.00E+07 & 9996841\\
\hline 8.00E+02 & 187 &2.00E+05 & 49782 &5.00E+07 & 12496477\\
\hline 9.00E+02 & 209 &3.00E+05 & 74722 &6.00E+07 & 14996075\\
\hline 1.00E+03 & 232 &4.00E+05 & 99679 &7.00E+07 & 17495807\\
\hline 2.00E+03 & 478 &5.00E+05 & 124633 &8.00E+07 & 19995483\\
\hline 3.00E+03 & 725 &6.00E+05 & 149606 &9.00E+07 & 22495362\\
\hline 4.00E+03 & 972 &7.00E+05 & 174588 &1.00E+08 & 24994959\\
\hline 5.00E+03 & 1214 &8.00E+05 & 199546 &2.00E+08 & 49992941\\
\hline 6.00E+03 & 1458 &9.00E+05 & 224520 &3.00E+08 & 74991365\\
\hline 7.00E+03 & 1706 &1.00E+06 & 249494 &4.00E+08 & 99989990\\
\hline 8.00E+03 & 1952 &2.00E+06 & 499298 &5.00E+08 & 124988741\\
\hline 9.00E+03 & 2206 &3.00E+06 & 749151 &6.00E+08 & 149987791\\
\hline 1.00E+04 & 2456 &4.00E+06 & 998978 &7.00E+08 & 174986827\\
\hline 2.00E+04 & 4925 &5.00E+06 & 1248914 &8.00E+08 & 199985791\\
\hline 3.00E+04 & 7415 &6.00E+06 & 1498757 &9.00E+08 & 224984997\\
\hline
\end{tabular}
\vspace{5mm}\caption{Valeurs de $N(\lambda)$ pour $d\cdot{10}^{p}$
avec $d=1,\ldots,9$ et $p=2,\ldots,8$.}\label{bigTabTrueCount}
\end{center}
\end{table}
\newpage

\subsection{Temps de calcul pour les procédures naïve et efficace d'évaluer la fonction de compte \texorpdfstring{N($\lambda$)}{}}

\noindent Il est intéressant de constater que si notre seul but est
d'évaluer $N(\lambda)$, alors il est nettement mieux d'utiliser
l'algorithme de la présente section et de ne pas procéder en
utilisant l'algorithme de la section 2.2. Voici à cet effet quelques
comparaisons des temps de calculs.

\begin{table}[ht]
\begin{center}
\begin{tabular}{|c|c|c|}
\hline $\lambda$ & $N(\lambda)$ & Temps requis (secondes) \\
\hline 100 & 21 & 0.024 \\
\hline 1000 & 232 & 0.118\\
\hline 10000 & 2456 & 0.513\\
\hline 100000 & 24842 & 1.839\\
\hline 1000000 & 249494 & 7.704\\
\hline
\end{tabular}
\vspace{5mm} \caption{Temps requis pour évaluer $N(\lambda)$
rapidement avec l'algorithme de la section 3.3}
\label{TabTrueCount}\vspace{10mm}
\begin{tabular}{|c|c|c|c|}
\hline l (rang) & $\lambda_{l}$ & $\lambda_{l+1}$ & Temps requis
(secondes)\\
\hline 21 & 98.726272 & 103.499453 & 0.056\\
\hline 232 & 989.729086 & 1001.729979 & 0.928\\
\hline 2456 & 9998.868757 & 10038.582253 & 15.542\\
\hline 24842 & 99997.791736 & 100024.359657 & 392.066\\
\hline 249494 & 999994.341030 & 1000020.429162 & 16016.677\\
\hline
\end{tabular}
\vspace{5mm} \caption{Temps requis pour évaluer $N(\lambda)$
naïvement avec l'algorithme de la section 2.2}
\label{TabPremierAlgo}
\end{center}
\end{table}

\noindent Les valeurs de la table \ref{TabTrueCount} peuvent être
utlisées pour vérifier celles de la table \ref{TabPremierAlgo} et
vice-versa. En d'autres termes, nous pouvons utiliser l'algorithme
de la section 2.2 pour vérifier celui de la présente section et
vice-versa.

\chapter*{Conclusion}

\noindent Finalement, voici une liste résumant les points importants
de ce mémoire. Les points 1,2 et 4 sont des résultats et/ou idées
propres à l'auteur de ce mémoire.

\begin{enumerate}
\item La première ligne nodale de la deuxième fonction
propre $u_{2}(r,\theta)$ sur un secteur $S(\alpha)$ du disque de
rayon 1 et d'angle $\alpha\in(0,2\pi)$ est donnée par l'équation
\begin{eqnarray*}
r&=&\frac{j_{1}^{2}(\nu)}{j_{2}^{2}(\nu)}\phantom{1}\textrm{si}\phantom{1}\alpha<\alpha_{0}\\
\theta&=&\frac{\alpha}{2}\phantom{1}\textrm{si}\phantom{1}\alpha>\alpha_{0}
\end{eqnarray*}
où $\nu=\frac{\pi}{\alpha}$ et où $j_{k}(\nu)$ est le k\ieme{} zéro
de la fonction de Bessel d'ordre $\nu$. L'angle $\alpha_{0}$ est
l'angle critique où la ligne nodale n'est pas définie et il est
donné implicitement en résolvant l'équation
$j_{2}(\nu)-j_{1}(2\nu)=0$ où $\nu=\frac{\pi}{\alpha}$. La
multiplicité de la 2\ieme{} valeur propre vaut 2 lorsque l'angle
$\alpha=\alpha_{0}$.

\item Soit $m\in\mathbb{N}$. Ordonner les valeurs $m$ premières valeurs propres de
$\triangle$ sur le disque de rayon 1 ou sur un secteur est
équivalent à trouver la configuration des ensembles nodaux des $m$
premières fonctions propres qui y sont associées. En utilisant le
théorème de Courant, j'ai crée un algorithme permettant de créer la
liste des $m$ premières valeurs propres. Une analyse de son temps
d'exécution en fonction de $m$ pourrait être fait en connaissant une
borne sur le temps requis pour évaluer les zéros des fonctions de
Bessel.

\item Soit la fonction de compte $N(\lambda)=\sum_{\lambda_{j}<\lambda}{1}$ des valeurs propres du
problème $\triangle{u}+\lambda{u}=0$ avec des conditions de
Dirichlet sur $\partial{D}$. Une explication de l'article de
\cite{KuFe_1965} a permis de comprendre pourquoi asymptotiquement
\begin{displaymath}
N(\lambda)=\frac{\lambda}{4}-\frac{\sqrt{\lambda}}{2}+O\big(\lambda^{\frac{1}{3}}\big).
\end{displaymath}

\item En utilisant la monotonicité des zéros des fonctions de Bessel
ainsi qu'un principe de \og marche \fg{} sur un ensemble de pairs
d'entiers bien déterminé par le problème en soit, j'ai créé un
algorithme permettant d'évaluer exactement et efficacement
$N(\lambda)$. Si notre seul but est d'évaluer $N(\lambda)$, le
2\ieme{} algorithme est beaucoup plus efficace qui si nous
procédions naïvement en utilisant le 1\ier{} algorithme pour trouver
le grand index $l$ tel que $N(\lambda)\leq l$ et $N(\lambda)> l$.
\end{enumerate}

\begin{center}$\mathfrak{Fin}$\end{center}

\appendix
\chapter{Liste des programmes en Maple et en Matlab}

\noindent Les programmes donnés ci-dessous ont été ajusté pour le
format de ce mémoire. Ne faites pas \og copier et coller \fg{} sans
réajuster les lignes des programmes.

\section{Programme pour obtenir la suite des \textrm{m} premières
valeurs/fonctions propres pour $m\in\mathbb{N}$}

\begin{center}\line(1,0){300}\end{center}
\noindent 
\begin{verbatim}
# La variable "nbEV" est le nombre de valeurs propres
# désirées.
#
# Le temps requis à la fin de l'exécution est écrit
# dans le fichier "ProgState.txt".
#
# Thm. de Courant : nb. de comp. conn. <= 0.7*l où
# l est le rang.
#
# Lancer en arrière-plan: "maple -i Prog2BesFctEV.txt &"
###
# The variable "nbEV" indicates the number of wanted
# eigenvalues.
#
# Time required at the of the execution of the program
# is given in the file "ProgState.txt".
#
# Courant Nodal Thm: nb. of connected components <= 0.7*l
# where l is the rank.
#
# Send it in background mode : "maple -i Prog2BesFctEV.txt &"
#

Time0:=time();
nbEV:=1000000;

# nombre de décimal pour les calculs
# number of decimal for calculations
# with this accuracy is 1e-6 for all the zeros.
Digits:=ceil(log(nbEV)/log(10))+1+6;

# calcul des grandeurs de tableaux
# evaluate size of tables
sizeTAB:=ceil(0.7*nbEV);
for nn from 1 to ceil(0.5*0.7*nbEV)
 do sizeTAB:=sizeTAB+ceil(0.5*0.7*nbEV/nn)
end;

TABZero:=Array(1..sizeTAB);
TABInd:=Array(1..sizeTAB);

#cas nn = "ordre" = 0

foo:=1;

# la multiplicité est 1.
# multiplicity is 1
for kk from 1 to ceil(0.7*nbEV) do
    TABZero[foo]:=evalf(BesselJZeros(0,kk)^2);
    TABInd[foo]:=0;
    foo:=foo+1;
end do;

#STATE
fidWst:=fopen("ProgState.txt",WRITE);
fprintf(fidWst,"\nn = 0 done");
fclose(fidWst);


#cas nn = "ordre" > 0
# la multiplicité vaut 2.
# multiplicity is 2.
for nn from 1 to ceil(0.5*0.7*nbEV) do

    for kk from 1 to ceil(0.5*0.7*nbEV/nn) do
    TABZero[foo]:=evalf(BesselJZeros(nn,kk)^2);
    #printf("\n*** %.6f",TABZero[foo]);
    TABInd[foo]:=1;
    foo:=foo+1;
    end do;

    ##STATE
    #fidWst:=fopen("ProgState.txt",WRITE);
    #fprintf(fidWst,"\nComputation n = %d done",nn);
    #fclose(fidWst);

end do;

# tri sur place
# sort with respect to a categorical value in right column

C:=zip((x,y)->[x,y],TABZero,TABInd);
f:=(x,y)->evalb(x[1]<=y[1]);
E:=sort(C,f);
F:=map(x->x[1],E);
G:=map(x->x[2],E);

#STATE
fidWst:=fopen("ProgState.txt",WRITE);
fprintf(fidWst,"\nSorting done");
fclose(fidWst);

TAB_EV:=Array(1..sizeTAB);

ABC:=1;
for ii from 1 by 1 while ii<=nbEV do
    #if order=0 then multiplicity = 1;
    if (G[ABC]=0) then TAB_EV[ii]:=F[ABC];
    #else order not 0 then multiplicity = 2;
    else TAB_EV[ii]:=F[ABC]; TAB_EV[ii+1]:=F[ABC]; ii:=ii+1;
    end if;
    ABC:=ABC+1;
end do;

#STATE
fidWst:=fopen("ProgState.txt",WRITE);
fprintf(fidWst,"\nList of EVs created");
fclose(fidWst);

fidW:=fopen("ListEVs.txt",WRITE);
for ii from 1 to nbEV do
    fprintf(fidW,"%d %.6f\n",ii,TAB_EV[ii]);
end do;
fclose(fidW);

#STATE
fidWst:=fopen("ProgState.txt",WRITE);
fprintf(fidWst,"\nCPU time required = %.4f",time()-Time0);
fclose(fidWst);

## for computational checks
#fidW:=fopen("TMP_01.txt",WRITE);
#for ii from 1 to (foo-1) do
#    fprintf(fidW,"%d %.6f\n",G[ii],F[ii]);
#end do;
#fclose(fidW);

quit;
\end{verbatim}
\begin{center}\line(1,0){300}\end{center}

\section{Programme pour évaluer la fonction de compte}

\begin{center}\line(1,0){300}\end{center}
\noindent 
\begin{verbatim}
ppBY:=1;
ppMAX:=10;
aaBY:=4;
aaMAX:=9;
abc:=0;
for pp from 1 by ppBY to ppMAX do
    for aa from 1 by aaBY to aaMAX do
      abc:=abc+1;

    end do;
end do;
SizeSample:=abc;
SqRLamb:=Array(1..SizeSample);
TimeTab:=Array(1..SizeSample);


Digits:=10;

# Ajuster boucles en fonction nombre de cas voulus
# donc SizeSample...
abc:=0;
for pp from 1 by ppBY to ppMAX do
    for aa from 1 by aaBY to aaMAX do

      abc:=abc+1;
      SqRLamb[abc]:=evalf(sqrt(aa*(10^pp)));
    end do;
end do;

fidOut:=fopen("ResCount.txt",WRITE):
fclose(fidOut);
fidT:=fopen("TimeTrueCount.txt",WRITE);
fclose(fidT);


for aa from 1 to SizeSample do

    tmpT0:=time();

    EigMax:=SqRLamb[aa];
    # K0 = max(k, zero(ordre 0, index k) <= lambda)
    for K0 from 1 to infinity do
    tmp:=evalf(BesselJZeros(0,K0));
    if (tmp>EigMax) then kk:=K0-1; break; end;
    end do;

    # N1 = max(n, zero(ordre n, index 1) <= lambda)
    for N1 from 0 to infinity do
    tmp:=evalf(BesselJZeros(N1,1));
    if (tmp>EigMax) then sizetab:=N1; break; end;
    end do;

    Ord:=Array(1..N1+1);
    Ind:=Array(1..N1+1);
    Ord[1]:=0;
    Ind[1]:=kk;

    # on monte au ciel
    for nn from 1 to N1 do

    tmp:=evalf(BesselJZeros(nn,kk));
    # on recule au Moyen-Age
    for ll from 1 while tmp>EigMax to infinity do
        kk:=kk-1;
        tmp:=evalf(BesselJZeros(nn,kk));
    end;

    Ord[nn+1]:=nn;
    Ind[nn+1]:=kk;
    if (kk = 0) then break; end;

    end do;

    # Ind[1], c'est pour l'ordre zéro donc
    # multiplicité est 1
    Tcount:=Ind[1];
    for bb from 2 to (N1+1) do
    # multiplicité est 2
    Tcount:=Tcount+2*Ind[bb];
    end do;

    fidOut:=fopen("ResCount.txt",APPEND):
    fprintf(fidOut,"%.12E %d\n",round(SqRLamb[aa]^2),Tcount);
    fclose(fidOut);
    fidT:=fopen("TimeTrueCount.txt",APPEND);
    fprintf(fidT,"N(%.12E) = %d *** Temps requis = %.6f
     secondes.\n",round(SqRLamb[aa]^2),Tcount,time()-tmpT0);
    fclose(fidT);
end do;

quit;
\end{verbatim}
\begin{center}\line(1,0){300}\end{center}

\section{Évaluation des points stationnaires des fonctions
apparaissant dans le développement asymptotique d'Olver}

\noindent Rappelons que
\begin{eqnarray*}
\zeta(z)&=&\bigg(\frac{3}{2}\int_{z}^{1}{\frac{\sqrt{1-t^{2}}}{t}dz}\bigg)^{\frac{2}{3}}\\
&=&\bigg(\frac{3}{2}\log\Big(\frac{1+\sqrt{1-z^{2}}}{z}\Big)-\frac{3}{2}\sqrt{1-z^{2}}\bigg)^{\frac{2}{3}}\phantom{12}0<z<1\\
\zeta(z)&=&-\bigg(\frac{3}{2}\int_{1}^{z}{\frac{\sqrt{t^{2}-1}}{t}dt}\bigg)^{\frac{2}{3}}\\
&=&-\bigg(\frac{3}{2}\sqrt{z^{2}-1}-\frac{3}{2}\arccos\Big({\frac{1}{z}}\Big)\bigg)^{\frac{2}{3}}\phantom{12}z\geq1
\end{eqnarray*}
alors pour la fonction pour la fonction
$|\zeta(z)|^{\frac{1}{2}}B_{0}(\zeta(z))$ où $B_{0}$ est donnée par
\begin{eqnarray*}
B_{0}(\zeta)&=&-\frac{3t-5t^{3}}{24\zeta^{\frac{1}{2}}}-\frac{5}{48\zeta^{2}}\phantom{12}\zeta>0\phantom{12}\textrm{ou}\phantom{12}0\leq
z\leq1\\
B_{0}(\zeta)&=&\frac{3t_{1}+5t_{1}^{3}}{24\zeta^{\frac{1}{2}}}-\frac{5}{48\zeta^{2}}\phantom{12}\zeta<0\phantom{12}\textrm{ou}\phantom{12}z\geq1\\
\end{eqnarray*}
et où
\begin{displaymath}
t=\frac{1}{\sqrt{1-z^{2}}}\phantom{12}\textnormal{et}\phantom{12}t_{1}
= \frac{1}{\sqrt{z^{2}-1}}
\end{displaymath}
alors voici le code évaluant ses points stationnaires.
\begin{center}\line(1,0){300}\end{center}
\begin{verbatim}
restart;
Digits:=20;

#I put t below but it IS t1 in the text before!
t:=(z^2 - 1)^(-0.5);

K:=1.5*sqrt(z^2 - 1) - 1.5*arccos(1/z);

b0:=(3*t)/(24*(K^(1/3))) +(5*(t^3))/(24*(K^(1/3))) -
                                           5/(48*(K^(4/3)));

foo:=diff(((abs(K))^(1/3))*b0,z$1);

ptstat:=fsolve(foo,z=2.71);

xipt:=evalf(-(subs(z=ptstat,K))^(2/3));

totval:=(1/12)+2*evalf(subs(z=ptstat,((abs(K))^(1/3))*b0));

fidTerm:=fopen(terminal,WRITE);
fprintf(fidTerm,"\n\n***** Total variation = %.20f \n\n",totval);

fidFile:=fopen("VarTotResultB0.txt",WRITE);
fprintf(fidFile,"\n\n***** Variation totale = %.20f
                                              \n\n",totval);
fprintf(fidFile,"\n\n***** Point stationnaire autre que z=1 ;
                                      z = %.20f\n\n",ptstat);
fprintf(fidFile,"\n\n***** xi(point stationnaire) =
                                              %.20f\n\n",xipt);
fprintf(fidFile,"\n\n***** sqrt(|xi|)B0(xi) = %.20f\n\n",
                   evalf(subs(z=ptstat,((abs(K))^(1/3))*b0)));
fclose(fidFile);
\end{verbatim}
\begin{center}\line(1,0){300}\end{center}

\noindent Pour la fonction $|\zeta(z)|^{\frac{1}{2}}B_{1}(\zeta(z))$
où $B_{1}$ est donnée par
\begin{eqnarray*}
B_{1}(\zeta)&=&\frac{-(30375t^{3}-369603t^{5}+765765t^{7}-425425t^{9})}{414720\zeta^{\frac{1}{2}}}\\
&&\phantom{1}-\frac{(405t^{2}-2310t^{4}+1925t^{6})}{55296\zeta^{2}}\\
&&\phantom{1}-\frac{(1155t-1925t^{3})}{110592\zeta^{\frac{7}{2}}}-\frac{85085}{663552\zeta^{5}}\phantom{12}\zeta>0\phantom{12}\textrm{ou}\phantom{12}0\leq z\leq 1\\
B_{1}(\zeta)&=&\frac{(30375t_{1}^{3}+369603t_{1}^{5}+765765t_{1}^{7}+425425t_{1}^{9})}{414720\zeta^{\frac{1}{2}}}\\
&&\phantom{1}+\frac{(405t_{1}^{2}+2310t_{1}^{4}+1925t_{1}^{6})}{55296\zeta^{2}}\\
&&\phantom{1}+\frac{(1155t_{1}+1925t_{1}^{3})}{110592\zeta^{\frac{7}{2}}}+\frac{85085}{663552\zeta^{5}}\phantom{12}\zeta<0\phantom{12}\textrm{ou}\phantom{12}z\geq 1\\
\end{eqnarray*}
et où $t$ et $t_{1}$ sont comme précédemment, alors voici le code
évaluant ses points stationnaires.
\begin{center}\line(1,0){300}\end{center}
\begin{verbatim}
restart;
Digits:=20;

# 0<z<1 ou bien xi>0

t:=(1-z^2)^(-0.5);

K:=1.5*log((1+sqrt(1-z^2))/z) - 1.5*sqrt(1-z^2);

b11:=-(30375*t^3 - 369603*t^5 + 765765*t^7 -
                           425425*t^9)/(414720*(K^(1/3)));
b12:=-(405*t^2 - 2310*t^4 + 1925*t^6)/(55296*(K^(4/3)));
b13:=-(1155*t-1925*t^3)/(110592*(K^(7/3)));
b14:=-85085/(663552*(K^(10/3)));
b1:=b11+b12+b13+b14;

foo:=diff(((abs(K))^(1/3))*b1,z$1);

ptstat1:=fsolve(foo,z=0.5);

##xipt:=evalf((subs(z=ptstat1,K))^(2/3));

tmp0:=evalf(limit(((abs(K))^(1/3))*b1,z=0));
tmp1:=evalf(subs(z=ptstat1,((abs(K))^(1/3))*b1));


###################################################

# z>1 ou bien xi<0

t1:=(z^(2)-1)^(-0.5);

K1:=1.5*sqrt(z^2 - 1) - 1.5*arccos(1/z);

b111:=(30375*t1^3 + 369603*t1^5 + 765765*t1^7 +
                          425425*t1^9)/(414720*(-(K1)^(1/3)));
b121:=(405*t1^2 + 2310*t1^4 + 1925*t1^6)/(55296*(K1^(4/3)));
b131:=(1155*t1+1925*t1^3)/(110592*(-(K1)^(7/3)));
b141:=85085/(663552*(K1^(10/3)));
b11:=b111+b121+b131+b141;

foo1:=diff(((abs(K1))^(1/3))*b11,z$1);

ptstat3:=fsolve(foo1,z=1.5);

tmp3:=evalf(subs(z=ptstat3,((abs(K1))^(1/3))*b11));
tmp4:=evalf(limit(((abs(K1))^(1/3))*b11,z=infinity));

totval:=tmp0+tmp4-2*tmp1-2*tmp3;
fidFile:=fopen("VarTotResultB1.txt",WRITE);
fprintf(fidFile,"\n\n***** Variation totale = %.20f
                                          \n\n",totval);
fprintf(fidFile,"\n\n***** Point stationnaire no. 1
                   autre que z=1 ; z1 = %.20f\n\n",ptstat1);
fprintf(fidFile,"\n\n***** Point stationnaire no. 2
                   autre que z=1 ; z2 = %.20f\n\n",ptstat3);
fprintf(fidFile,"\n\n***** sqrt(|xi|)B1(xi) pour
                                  xi(z1) = %.20f\n\n",tmp1);
fprintf(fidFile,"\n\n***** sqrt(|xi|)B1(xi) pour
                                  xi(z2) = %.20f\n\n",tmp3);
fprintf(fidFile,"\n\n***** limit(|xi|)B1(xi)
                    (vers xi(0) = infty) = %.20f\n\n",tmp0);
fprintf(fidFile,"\n\n***** limit(|xi|)B1(xi)
               (vers xi(infty) = -infty) = %.20f\n\n",tmp4);
fclose(fidFile);
\end{verbatim}
\begin{center}\line(1,0){300}\end{center}

\section{Programme en Matlab pour résoudre les problèmes de valeurs
propres et visualiser les ensembles nodaux des fonctions propres}

\noindent Pour de l'information sur ces programmes, vous pouvez
taper par exemple
\begin{verbatim}
>>help SpDisque
\end{verbatim}
dans Matlab dans le répertoire où vous enregistrez les fichiers pour
savoir quoi faire.

\begin{center}\line(1,0){300}\end{center}
\begin{verbatim}
function [] = SpDisque(NbEV,NbCourbesNiveau);
% NbEV : nombre de valeurs propres désirées
% NbCourbesNiveau : quantité des courbes à nouveau pour les figures
%
% Conditions de Dirichlet aux frontières sont appliquées.
%
% DÉPENDANCE : DisqueGeo.m
%
% Example : SpSecteur(50,[0.25*pi 2*exp(1)],1000);

% Conditions aux frontières
BdyCond = [1 1 1 1  1   1 '0' '0' '1' '0';  %1 Dir
           1 1 1 1  1   1 '0' '0' '1' '0';  %2 Dir
           1 1 1 1  1   1 '0' '0' '1' '0';  %3 Dir
           1 1 1 1  1   1 '0' '0' '1' '0']'; %4 Dir

% Constantes pour la description de l'équation
c = 1; a = 0; d = 1;

% Efface fichiers information sur mesh et état du programme
fidIG = fopen('SpDisqueMeshQual.txt','wt');
fclose(fidIG);
fidEt = fopen('SpDisqueEtat.txt','wt');
fclose(fidEt);
fidEr = fopen('SpDisqueErreurs.txt','wt');
fclose(fidEr);

try

  % Mesh ici
  %[p,e,t]=initmesh('DisqueGeo','Hmax',2e6,'Init','off','Box','off'
   ,'Hgrad',1.05,'Jiggle','minimum','JiggleIter',120);
  [p,e,t]=initmesh('DisqueGeo'); %%% for testing
  [p,e,t]=refinemesh('DisqueGeo',p,e,t);
  [p,e,t]=refinemesh('DisqueGeo',p,e,t);
  [p,e,t]=refinemesh('DisqueGeo',p,e,t);
  %[p,e,t]=refinemesh('SecteurGeo',p,e,t);


  % stat sur qualité du mesh
  Mesh_Quality = pdetriq(p,t);
  fidIG = fopen('SpDisqueMeshQual.txt','at');
  fprintf(fidIG,'\nQualité moyenne des triangles = %.6f',
   mean(Mesh_Quality));
  fprintf(fidIG,'\nÉcart-type sur les triangles = %.6f',
   std(Mesh_Quality));
  fprintf(fidIG,'\nMinimum de qualité atteint = %.6f',
   min(Mesh_Quality));
  fprintf(fidIG,'\nMaximum de qualité atteint = %.6f\n',
   max(Mesh_Quality));
  fclose(fidIG);

  % résoudre ici
  MaxEig=1.0;
  while(1)

    fidEt = fopen('SpDisqueEtat.txt','wt');
    fprintf(fidEt,'\nAppel à pdeeig avec MaxEig = %.6f\n\n',
     MaxEig);
    fclose(fidEt);

    [EF_tab,EV_tab]=pdeeig(BdyCond,p,e,t,c,a,d,[0 MaxEig]);

    if length(EV_tab) >= NbEV
      break;
    else
      MaxEig = 2.00*MaxEig;
    end;
  end;


  % fait des beaux dessins ici :-) Sans-Dessin ha ha

  for jj = 1:1:NbEV
    fnamestr = char(strcat({'Disque'},{'EF'},{num2str(jj)}));
    titlestr = char(strcat({'Disque '},{'  \lambda_{'},
     {num2str(jj)},'}=',{num2str(EV_tab(jj))}));
    pdeplot(p,[],t,'xydata',EF_tab(:,jj),'xystyle','flat',
     'colormap','Jet','colorbar','off');
    title(titlestr);
    %print('-depsc2','-r0',fnamestr);
    print('-djpeg','-r0',fnamestr);

  end;


catch

  v = lasterror;
  fidEr = fopen('SpDisqueErreurs.txt','at');
  fprintf(fidEr,'\n%s\n',v.message);
  fclose(fidEr);

end;
\end{verbatim}
\begin{center}\line(1,0){300}\end{center}
\begin{verbatim}
function [] = SpSecteur(NbEV,angles,NbCourbesNiveau);
% NbEV : nombre de valeurs propres désirées
% angles : vecteur des angles de chaque secteur
% NbCourbesNiveau : quantité des courbes à nouveau pour les figures
%
% Conditions de Dirichlet aux frontières sont appliquées.
%
% DÉPENDANCE : SecteurGeo.m
%
% Exemple : SpSecteur(3,[1.0 1.11 1.112 1.1125 1.11253 1.112531
   1.1125316 1.112531602 1.1125316027 1.11253160273 1.112531602739
   1.112531602739001 1.2],1000);


% Conditions aux frontières
BdyCond = [1 1 1 1  1   1 '0' '0' '1' '0';  %1 Dir
           1 1 1 1  1   1 '0' '0' '1' '0';  %2 Dir
           1 1 1 1  1   1 '0' '0' '1' '0']'; %3 Dir

% Constantes pour la description de l'équation
c = 1; a = 0; d = 1;

% Efface fichiers information sur mesh et état du programme
fidIG = fopen('SpSecteurMeshQual.txt','wt');
fclose(fidIG);
fidEt = fopen('SpSecteurEtat.txt','wt');
fclose(fidEt);
fidEr = fopen('SpSecteurErreurs.txt','wt');
fclose(fidEr);

try
  % Itération sur les cas différents i.e. les angles
  for ii = 1:1:length(angles)

    fid = fopen('AngleGeoSecteur.txt','wt');
    fprintf(fid, '%.16f\n',angles(ii));
    fclose(fid);

    % Mesh ici
    %[p,e,t]=initmesh('SecteurGeo','Hmax',2e6,'Init','off','Box',
     'off','Hgrad',1.05,'Jiggle','minimum','JiggleIter',120);
    [p,e,t]=initmesh('SecteurGeo'); %%% for testing
    [p,e,t]=refinemesh('SecteurGeo',p,e,t);
    [p,e,t]=refinemesh('SecteurGeo',p,e,t);
    [p,e,t]=refinemesh('SecteurGeo',p,e,t);
    %[p,e,t]=refinemesh('SecteurGeo',p,e,t);


    % stat sur qualité du mesh
    Mesh_Quality = pdetriq(p,t);
    fidIG = fopen('SpSecteurMeshQual.txt','at');
    fprintf(fidIG,'\nSecteur d''angle %.16f',angles(ii));
    fprintf(fidIG,'\nQualité moyenne des triangles = %.6f',
     mean(Mesh_Quality));
    fprintf(fidIG,'\nÉcart-type sur les triangles = %.6f',
     std(Mesh_Quality));
    fprintf(fidIG,'\nMinimum de qualité atteint = %.6f',
     min(Mesh_Quality));
    fprintf(fidIG,'\nMaximum de qualité atteint = %.6f\n',
     max(Mesh_Quality));
    fclose(fidIG);

    % résoudre ici
    MaxEig=1.0;
    while(1)

      fidEt = fopen('SpSecteurEtat.txt','wt');
      fprintf(fidEt,'\nAppel à pdeeig avec MaxEig = %.6f\n\n'
       ,MaxEig);
      fclose(fidEt);

      [EF_tab,EV_tab]=pdeeig(BdyCond,p,e,t,c,a,d,[0 MaxEig]);

      if length(EV_tab) >= NbEV
        break;
      else
        MaxEig = 2.00*MaxEig;
      end;
    end;

    % fait des beaux dessins ici :-) Sans-Dessin ha ha
    for jj = 1:1:NbEV
      fnamestr = char(strcat({'Sect'},{num2str(ii)},{'EF'},
       {num2str(jj)}));
      titlestr = char(strcat({'Secteur d''angle '},
       {num2str(angles(ii),16)},{'  \lambda_{'},{num2str(jj)},'}=',
        {num2str(EV_tab(jj))}));

      pdeplot(p,[],t,'xydata',EF_tab(:,jj),'xystyle','flat',
       'colormap','Jet','colorbar','off');
      title(titlestr);
      %print('-depsc2','-r0',fnamestr);
      print('-djpeg','-r0',fnamestr);
    end;
  end;

catch
  v = lasterror;
  fidEr = fopen('SpSecteurErreurs.txt','at');
  fprintf(fidEr,'\n%s\n',v.message);
  fclose(fidEr);
end;

%exit;
\end{verbatim}
\begin{center}\line(1,0){300}\end{center}
\begin{verbatim}
function [x,y]=DisqueGeo(bs,s)
%   Géométrie du disque de rayon 1 pour résolution PDE
nbs=4;

if nargin==0,
  x=nbs; % number of boundary segments
  return
end;

d=[   0.0  0.5*pi  1  0;   %1
      0.5*pi pi    1  0;   %2
      pi 1.5*pi   1  0;   %3
      1.5*pi 2*pi 1  0]'; %4


bs1=bs(:)';

if find(bs1<1 | bs1>nbs),
  error('Invalid boundary segment.',
   'Non existent boundary segment number.')
end

if nargin==1,
  x=d(:,bs1);
  return
end

x=zeros(size(s));
y=zeros(size(s));
[m,n]=size(bs);
if m==1 && n==1,
  bs=bs*ones(size(s)); % expand bs
elseif m~=size(s,1) || n~=size(s,2),
  error('PDE:Sector Geometry : SizeBs',
   'bs must be scalar or of same size as s.');
end

if ~isempty(s),

% boundary segment 1
ii=find(bs==1);
x(ii)=cos(s(ii));
y(ii)=sin(s(ii));

% boundary segment 2
ii=find(bs==2);
x(ii)=cos(s(ii));
y(ii)=sin(s(ii));

% boundary segment 3
ii=find(bs==3);
x(ii)=cos(s(ii));
y(ii)=sin(s(ii));

% boundary segment 4
ii=find(bs==4);
x(ii)=cos(s(ii));
y(ii)=sin(s(ii));

end
\end{verbatim}
\begin{center}\line(1,0){300}\end{center}
\begin{verbatim}
function [x,y]=SecteurGeo(bs,s)
%   Géométrie du secteur de rayon 1 pour PDE solving
nbs=3;

if nargin==0,
  x=nbs; % number of boundary segments
  return
end;

[alpha] = textread('AngleGeoSecteur.txt');

if (alpha>=2*pi)
  error('BadParam:GeoFile:DesignerCheck',
  'L''angle pour un secteur < 2*pi');
end;

d=[   0.0  1.0     1  0;   %1
      0.0  alpha   1  0;   %2
      1.0  0.0     1  0]'; %3


bs1=bs(:)';

if find(bs1<1 | bs1>nbs),
  error('Invalid boundary segment.',
   'Non existent boundary segment number.')
end

if nargin==1,
  x=d(:,bs1);
  return
end

x=zeros(size(s));
y=zeros(size(s));
[m,n]=size(bs);
if m==1 && n==1,
  bs=bs*ones(size(s)); % expand bs
elseif m~=size(s,1) || n~=size(s,2),
  error('PDE:Sector Geometry : SizeBs',
   'bs must be scalar or of same size as s.');
end

if ~isempty(s),

% boundary segment 1
ii=find(bs==1);
x(ii)=s(ii);
y(ii)=0.0;

% boundary segment 2
ii=find(bs==2);
x(ii)=cos(s(ii));
y(ii)=sin(s(ii));

% boundary segment 3
ii=find(bs==3);
x(ii)=s(ii)*cos(alpha);
y(ii)=s(ii)*sin(alpha);

end
\end{verbatim}
\begin{center}\line(1,0){300}\end{center}

\section{Précision des calculs}

\begin{Rem}[Calcul des zéros et précision]\label{remcapoute1}
Par (\ref{capoute1}), la longueur décimale de $j_{\nu}(k)$ n'excède
pas la longueur décimale de $\nu$. Par conséquent pour calculer la
partie entière des zéros $j_{k}(\nu)$ sans perte de précision, il ne
faut que $\lceil\log_{10}(\nu)\rceil$ décimales. Pour être prudent,
nous pouvons exiger $\lceil\log_{10}(\nu)\rceil+1$ décimales de
précision. Si nous voulons évaluer les zéros à une précision absolue
de $10^{-p}$, alors il faut que les calculs se fassent avec
$\lceil\log_{10}(\nu)\rceil+1+p$ de décimales.
\end{Rem}

\noindent La remarque précédente explique pourquoi la ligne de code
\begin{verbatim}
Digits:=ceil(log(nbEV)/log(10))+1+6;
\end{verbatim}
apparaît dans le programme de la section 2.2. \textbf{En effet}, le
plus petit ordre entier non nul $n$ des fonctions de Bessel
impliquées dans l'obtention de la suite spectrale
$\{\lambda_{j}\}_{j=1}^{m}$ est 1. Si, encore, nous dénotons par $k$
le rang d'un zéro, alors, par le théorème de Courant, $2nk\leq m$
comme il a été expliqué au chapitre 1 en tenant compte de la
multiplicité puisque ce sont des valeurs propres de multiplicité 2.
Si l'ordre est nul, alors $k\leq m$ et c'est une valeur propre de
multiplicité 1 comme il a été expliqué au chapitre 1. Dans le cas
d'un secteur, $nk\leq{m}$. Dans les deux cas (le disque ou le
secteur), nous avons que l'ordre maximal impliqué est inférieur à
$m$. Puisque $m$ est le nombre de valeurs propres désirées
c'est-à-dire $m:=:nbEV$ dans le programme, alors nous voyons d'où
sort la formule. Le choix d'imposer une précision de $10^{-6}$
s'explique par le fait que nous voulons éviter de trouver des zéros
ayant la même représentation jusqu'à $10^{-3}$ comme il est possible
d'en trouver dans la liste de 30 MB des valeurs propres $\lambda_{j}
$ pour $j=1,\ldots,10^{6}$.

\noindent Dans le cas du programme de la section 3.3, il faut
évaluer tous les zéros $j_{k}(n)\leq\sqrt{\lambda}$. Puisque la
valeur maximale de $\lambda$ pour laquelle nous voulons évaluer
$N(\lambda)$ est $9.00\cdot{10}^{8}$, alors le nombre de décimales
requises $\lceil\log_{10}(9.00\cdot{10}^{8})\rceil+1+0=10$ d'où la
ligne de code
\begin{verbatim}
Digits:=10;
\end{verbatim}
dans le programme de la section 3.3 évaluant $N(\lambda)$.

\section{Méthodes utilisées par le logiciel Maple pour évaluer les zéros}

\noindent Ci-dessous, nous trouvons une bonne partie de
l'information nécessaire pour savoir comment Maple évalue les zéros
des fonctions de Bessel. C'est un courriel qui m'a été envoyé par
l'équipe technique de Maplesoft.

\begin{verbatim}
Comments extracted from the source file:

Rayleigh's formula for the first zero of BesselJ function can be
obtained in reference:

G.N. Watson - A Treatise on the Theory of Bessel Functions,
2nd ed., CUP, 1966, page 502,


Olver's asymptotic expansion for large order can be obtained in
references:

M. Abramowitz and I. A. Stegun, Handbook of Mathematical Functions,
Dover, NY (1972)

Royal Society Mathmatical Tables, vol. 7, part III. Ed. F.W.J.
Olver, CUP, 1960

Olver, F.W.J., Philosophical Transactions of the Royal Society, A,
vol 247, pag. 307-368 (1954)


McMahon's asymptotic expansion for large n can be obtained in
references:

M. Abramowitz and I. A. Stegun, Handbook of Mathematical Functions,
Dover, NY (1972)

G.N. Watson - A Treatise on the Theory of Bessel Functions,
2nd ed., CUP, 1966

To find the N_th ZERO OF BESSEL FUNCTION of order m,
the procedure uses the following results proved by
Schafheitlin see
G.N. Watson - A Treatise on the Theory of Bessel Functions,
2nd ed., CUP, 1966:

For m<=1/2 the n_th zero of BesselJ[m](x) is in the interval:
[(n+m/2-1/4)*Pi .. (n+m/4-1/8)*Pi]

For 1/2<m<5/2 the n_th zero of BesselJ[m](x) is in the interval:
[(n+m/4-1/8)*Pi .. (n+m/2-1/4)*Pi]

For m>1/2, the zeros of BesselJ[m](x) which exceed
(2*m+1)*(2*m+3)/Pi are in the intervals:    [(p+m/2+1/2)*Pi ..
(p+m/2+3/4)*Pi] where p are consecutive positive integers. The
relation of p with n has been taken as p=n-1.

The zeros less than (2*m+1)*(2*m+3)/Pi will be found by a recursive
algorithm.

We use the following assymptotic aproximations for the first and
second zeros[see Watson]:

first zero: m + 1.85*m^(1/3)) + O(m^(-1/3))

second zero: m + 3.24*m^(1/3)) + O(m^(-1/3))

To find the n-th zero j(m,n) we use the interval
      [j(m,n-1) .. j(m,n-1) + incr]
where the increment is
      incr = j(m,n-1) - j(m,n-2)

The other part of the solution can be found by setting a Maple
Interface variable to display the initial routine that is used.

interface(verboseproc=3); eval(BesselJZeros);

I hope this helps answer your questions.

Cheers,

Ian Maplesoft Technical Support www.maplesoft.com/support/

-----Original Message-----
From: Maplesoft Customer Service Sent: Monday, November 05, 2007
4:15 PM To: Claude Gravel; Maplesoft Technical Support Subject: RE:
Customer Service Request: (Web) Other
\end{verbatim}

\noindent Comme nous pouvons le constater, Maple utilise les
développements asymptotiques d'Olver pour représenter les fonctions
de Bessel pour des ordres grands. Les développements asymptotiques
d'Olver des zéros pour des ordres grands sont également utilisés
pour évaluer les zéros.

\end{document}